\documentstyle{amltd2004}
\begin{document}
\annalsline{158}{2003}
\received{April 6, 1999}
\startingpage{355}
\def\bye{\end{document}}
 \font\tenrm=cmr10
\def\ritem#1{\item[{\rm #1}]}
\input amssym.def
\input amssym.tex

\catcode`\@=11
\font\twelvemsb=msbm10 scaled 1100
\font\tenmsb=msbm10
\font\ninemsb=msbm10 scaled 800
\newfam\msbfam
\textfont\msbfam=\twelvemsb  \scriptfont\msbfam=\ninemsb
  \scriptscriptfont\msbfam=\ninemsb
\def\msb@{\hexnumber@\msbfam}
\def\Bbb{\relax\ifmmode\let\next\Bbb@\else
 \def\next{\errmessage{Use \string\Bbb\space only in math
mode}}\fi\next}
\def\Bbb@#1{{\Bbb@@{#1}}}
\def\Bbb@@#1{\fam\msbfam#1}
\catcode`\@=12

 \catcode`\@=11
\font\twelveeuf=eufm10 scaled 1100
\font\teneuf=eufm10
\font\nineeuf=eufm7 scaled 1100
\newfam\euffam
\textfont\euffam=\twelveeuf  \scriptfont\euffam=\teneuf
  \scriptscriptfont\euffam=\nineeuf
\def\euf@{\hexnumber@\euffam}
\def\frak{\relax\ifmmode\let\next\frak@\else
 \def\next{\errmessage{Use \string\frak\space only in math
mode}}\fi\next}
\def\frak@#1{{\frak@@{#1}}}
\def\frak@@#1{\fam\euffam#1}
\catcode`\@=12

\def\ZZ {{\Bbb Z}}
\def\NN {{\Bbb N}}
\def\QQ {{\Bbb Q}}
\def\RR {{\Bbb R}}
\def\CC {{\Bbb C}}

\def\Si{\Sigma}
\def\La{\Lambda}
\def\De{\Delta}
\def\Om{\Omega}
\def\Ga{\Gamma}
 \def\Th{\Theta}

\def\cA{{\cal A}}  \def\cG{{\cal G}} \def\cM{{\cal M}} \def\cS{{\cal S}}
\def\cB{{\cal B}}  \def\cH{{\cal H}} \def\cN{{\cal N}} \def\cT{{\cal T}}
\def\cC{{\cal C}}  \def\cI{{\cal I}} \def\cO{{\cal O}} \def\cU{{\cal U}}
\def\cD{{\cal D}}  \def\cJ{{\cal J}} \def\cP{{\cal P}} \def\cV{{\cal V}}
\def\cE{{\cal E}}  \def\cK{{\cal K}} \def\cQ{{\cal Q}} \def\cW{{\cal W}}
\def\cF{{\cal F}}  \def\cL{{\cal L}} \def\cR{{\cal R}} \def\cX{{\cal X}}
\def\cY{{\cal Y}}  \def\cZ{{\cal Z}}

\def\degre{{\rm deg}}
 \font\emi= cmmi10 scaled 1700 
\font\eightmr=cmr10
\title{A {\emi C}\hskip2pt\raise8pt\hbox{\eightmr 1}-generic dichotomy for\\ diffeomorphisms: Weak forms of\\ hyperbolicity
or infinitely many\\ sinks or sources}
  \def\titleheadline#1{\def\one{#1}\ifx\one\empty\else
\gdef\thetitle{{\frenchspacing%
\let\\ \relax
{#1}}}\fi}
\newif\ifshort
\def\shortname#1{\global\shorttrue\xdef
\theauthors{{\eightsc\uppercase{#1}}}}
\let\shorttitle\titleheadline
\shorttitle{ \eightsc\uppercase{A} {\eightpoint \it C} 
\hskip-2pt\raise4pt\hbox{{\tiny 1}}\eightsc\uppercase{-generic dichotomy for diffeomorphisms} }
 
 \acknowledgements{Partially supported by CNPQ, FAPERJ, IMPA, and PRONEX-Sistemas Din\^amicos,  Brazil, and Laboratoire de
Topologie (UMR 5584 CNRS) and Universit\'e de Bourgogne, France.}
 
 \twoauthors{C. Bonatti, L. J. D\'\i az,}{E. R. Pujals}
  
\institutions{
Institut de Math\'ematiques de Bourgogne, Universit\'e de Bourgogne,  Dijon, France\\
{\eightpoint {\it E-mail address\/}: bonatti@u-bourgogne.fr}\\
\vglue6pt
Pontif\'{\i}cia Universidade Cat\'olica do Rio de Janeiro (PUC-RJ), Rio de Janeiro, Brazil\\
{\eightpoint {\it E-mail address\/}: lodiaz@mat.puc-rio.br}\\
\vglue6pt
IMPA, Rio de Janeiro, Brazil\\
{\eightpoint {\it E-mail address\/}: enrique@impa.br}}
\vglue-8pt
\centerline{\it  A Ricardo Ma\~n{\rm \'{\it e}} {\rm (1948--1995),} por todo su trabajo}

\vglue12pt \centerline{\bf Abstract}
\vglue8pt
 We show that, for
every compact $n$-dimensional manifold, $n\geq 1$, there is a residual subset of Diff$^1(M)$ of diffeomorphisms for which the homoclinic class of any periodic saddle of $f$ verifies one of the following two possibilities: Either it is contained in the closure
of an infinite set of sinks or sources (Newhouse phenomenon), or it
presents some weak form of hyperbolicity called  {\em dominated
splitting\/} (this is a generalization of a bidimensional result of
Ma\~n\'e \cite{Ma3}). In particular, we show that any $C^1$-robustly transitive diffeomorphism admits a dominated splitting. 
 \vglue12pt

\centerline{{\bf R\'esum\'e}}
\vglue8pt
G\'en\'eralisant un r\'esultat de Ma\~n\'e sur les surfaces \cite{Ma3}, nous montrons que, en dimension quelconque, il existe un sous-ensemble r\'esiduel de Diff$^1(M)$ de diff\'eomor\-phismes pour lesquels la classe homocline de toute selle p\'eriodique hyperbolique poss\`ede deux comportements possibles : ou bien elle est incluse dans l'adh\'erence d'une infinit\'e de puits ou de sources (ph\'enom\`ene de Newhouse), ou bien elle pr\'esente une forme affaiblie d'hyperbolicit\'e appel\'ee  une {\em d\'ecomposition domin\'ee}.  En particulier nous montrons que
tout diff\'eomorphisme $C^1$-robustement transitif poss\`ede une d\'ecomposition domin\'ee.
\vglue12pt
\intro
 
\demo{Context}
The Anosov-Smale theory of uniformly hyperbolic systems has played a double role in the development of
the qualitative theory of dynamical systems.
On one hand, this theory shows that chaotic and random behavior can appear in a stable way for deterministic systems depending on a very small number of parameters. On the other hand, the chaotic systems admit  in this theory a quasi-complete description from the ergodic point of
view. Moreover the
hyperbolic attractors satisfy very simple statistical properties (see \cite{Si}, \cite{Ru}, and \cite{BoRu}): For Lebesgue almost every point in the 
topological basin
of the attractor, the time average 
of any continuous function along its orbit converges to the spatial average of the function by a probability measure whose support is the attractor.

However, since the end of the 60s,
one knows that this hyperbolic theory does not cover a dense
set of dynamics: There are examples of open sets of nonhyperbolic diffeomorphisms. More precisely,
\vglue5pt\noindent \hskip1em \hangindent =20pt \hangafter=1 $\bullet$ 
For every
 compact surface $S$ there exist nonempty open sets of Diff$^2(S)$ of diffeomorphisms whose nonwandering set is not hyperbolic (see \cite{N1}).
\vglue4pt\noindent \hskip1em \hangindent =20pt \hangafter=1 $\bullet$ 
Given any compact manifold $M$, with $\dim (M)\geq 3$, there are nonempty open subsets of Diff$^1(M)$ of diffeomorphisms whose nonwandering set is not hyperbolic (see, for instance, \cite{AS} and \cite{So} for the first examples).
\vglue5pt

In the $2$-dimensional case, at least in the $C^1$-topology, the heart of this phenomenon is 
closely related to the appearance of homoclinic tangencies:
For every compact surface $S$ the set of $C^1$-diffeomorphisms with homoclinic tangencies
is $C^1$-dense in the  complement in Diff$^1(S)$ 
of the closure of the Axiom A diffeomorphisms 
(that is a recent
result in \cite{PuSa}). 

Even if in this work we are concerned with the  
 $C^1$-topology, let us recall that 
Newhouse showed (see \cite{N1}) that 
generic unfoldings of homoclinic tangencies  create  
$C^2$-locally residual subsets of Diff$^2(S)$ of diffeomorphisms having an infinite set of  
sinks or sources. 
In this paper, by {\it $C^r$-Newhouse phenomenon\/} we mean  the
coexistence of infinitely many sinks or sources in a $C^r$-locally residual subset of Diff$^r(M)$.

The main motivation of  this article
comes from the following result of Ma\~n\'e (see \cite{Ma3} (1982)), which
gives, for $C^1$-generic diffeomorphisms of\break surfaces,
a dichotomy  between hyperbolic dynamics and the Newhouse\break phenomenon:
\enddemo

{\elevensc Theorem} {\rm (Ma\~n\'e)}.
{\it Let $S$ be a closed surface. Then there is a residual subset $\cR\subset$ {\rm Diff}$^1(S)${\rm ,}  
$\cR = \cR_1\coprod \cR_2${\rm ,} 
such that every $f\in \cR_1$ verifies the Axiom~{\rm A}
 and every $f\in \cR_2$ has infinitely many sinks or sources.}
\vglue8pt

Recall that a diffeomorphism of a manifold $M$ is {\it transitive\/} if it has a dense orbit in the whole
manifold. Such a diffeomorphism
is called {\it $C^r$-robustly transitive\/} if it belongs to the $C^r$-interior of the set of transitive diffeomorphisms.
Since transitive diffeomorphisms have neither sinks nor sources, a direct consequence from   Ma\~n\'e's result is 
the following:
\pagegoal=50pc
 \vskip4pt
{\it Every $C^1$\/{\rm -}\/robustly transitive  diffeomorphism on a  compact surface admits a
 hyperbolic structure on the whole manifold\/{\rm ;} i.e.{\rm ,} it is an Anosov diffeo\-morphism.}

Let us observe that 
Ma\~n\'e's result has no direct generalization to higher dimensions: For every $n\ge 3$ there are compact 
$n$-dimensional manifolds
supporting $C^1$-robustly transitive nonhyperbolic diffeomorphisms (in particular, without 
sources and sinks). 
All the examples of such diffeomorphisms, successively given by \cite{Sh} (1972) on the torus $T^4$, by \cite{Ma2} (1978) on $T^3$, by \cite{BD1} (1996) in many other manifolds (those supporting a transitive Anosov flow or of the form $M\times N$,
where $M$ is a manifold with an Anosov diffeomorphism and $N$ any compact manifold),
by \cite{B} (1996) and \cite{BoVi} (1998) 
examples in $T^3$  without any hyperbolic expanding direction 
and examples in $T^4$ without any hyperbolic direction, present some weak form of hyperbolicity, the newer the examples the weaker the form of hyperbolicity, but always 
exhibiting some remaining weak form of hyperbolicity. Let us be more precise.

Recall first that an invariant compact set $K$ of a diffeomorphism $f$
on a manifold $M$ is hyperbolic if the tangent bundle $TM|_K$ of $M$
over $K$ admits an $f_*$-invariant continuous splitting $TM|_K = E^s\oplus E^u$, such that $f_*$ uniformly contracts  the
vectors in $E^s$  and uniformly expands the vectors in  $E^u$. This means
that there is $n\in\NN$ such that $\|f_*^n(x)|_{E^s(x)}\|<1/2$ and $\|f_*^{-n}(x)|_{E^u(x)}\|\break <1/2$ for any $x\in K$
(where $||\cdot||$ denotes the norm).

The examples  of $C^1$-robustly transitive diffeomorphisms $f$  in \cite{Sh} and \cite{Ma2} let an invariant    splitting
$TM=E^s\oplus E^c\oplus E^u$,  where $f_*$  contracts
uniformly the vectors in $E^s$ and expands uniformly
 the vectors in $E^u$. Moreover, this splitting is dominated (roughly speaking, the contraction (resp.\ expansion) in
$E^s$ (resp.\ $E^u$) is stronger than the contraction
(resp.\ expansion)
in~$E^c$;
for details see Definition \ref{d.dominated} below), and the central bundle $E^c$ is one dimensional. The examples in \cite{BD1} admit also such a 
nonhyperbolic splitting, but the central bundle may have any dimension. The diffeomorphisms 
in \cite{B} have no stable bundle $E^s$ and admit a splitting $E^c\oplus E^u$, where the restriction of $f_*$ to $E^c$ is not uniformly contracting, but it
uniformly contracts  the area. Finally,  \cite{BoVi} gives examples of robustly transitive diffeomorphisms on $T^4$ without any uniformly stable or unstable bundles:
They leave invariant some dominated splitting $E^{cs}\oplus E^{cu}$, where the derivative of $f$ contracts uniformly
the area in $E^{cs}$ and expands uniformly
 the area in $E^{cu}$.

Roughly speaking, in this paper we see that,
if a transitive set does not admit a dominated splitting, 
then one can create as many sinks or sources as one  wants in any neighbourhood of this set. In particular,
$C^1$-robustly transitive diffeomorphisms always 
admit some dominated splitting.

Before stating our results  more precisely,  let us mention two previous
results on $3$-manifolds which  are the roots of this work:
\pagegoal=49pc

\vglue3pt\noindent \hskip1em \hangindent =20pt \hangafter=1 $\bullet$  \cite{DPU} shows that there is an
open and dense subset of
 $C^1$-robustly transitive  $3$-dimensional 
diffeomorphisms $f$  
admitting a dominated splitting $E^1\oplus E^2$ such
that at least 
one of the two bundles is uniformly hyperbolic (either stable or unstable). In that case, 
by terminology,   $f$ is {\it partially hyperbolic.\/} 
Moreover, \cite{DPU} also gives a semi-local version of
this result defining $C^1$-robustly transitive sets. Given a $C^1$-diffeomorphism
$f$, a compact set $K$ is {\it $C^1$-robustly transitive\/} if it is the
maximal $f$-invariant 
set in some neighbourhood $U$ of it and
if, for every $g$ $C^1$-close to $f$, the maximal 
invariant set $K_g = \bigcap_{\ZZ} g^n(U)$  is  also compact and  transitive. 

\vglue6pt\noindent \hskip1em \hangindent =20pt \hangafter=1 $\bullet$   \cite{BD2} gives examples of diffeomorphisms $f$ on
$3$-manifolds having  two saddles $P$ and $Q$ with a  pair of contracting and expanding complex (nonreal) eigenvalues,
respectively, which belong in a robust way to the same transitive set $\Lambda_f$. Clearly, this simultaneous presence of complex
contracting and expanding eigenvalues prevents the transitive set $\La_f$  from admitting a dominated splitting! Then \cite{BD2}
shows that, for a  $C^1$-residual subset of such diffeomorphisms,  the transitive set $\La_f$ is contained in the closure of the
(infinite) set of sources or sinks.
\vglue6pt
\pagegoal=48pc
 The results of these two papers seem to go in opposite directions, but here we show that they describe two sides of the
same phenomenon. In fact,we give here a framework which allows us to unify these results and generalize them in any dimension: In
the absence of weak hyperbolicity (more precisely, of a dominated splitting) one can create an arbitrarily large number of sinks or
sources.

In the nonhyperbolic context, 
the classical notion of basic pieces (of the  Smale  spectral 
decomposition theorem)  is not defined and an important problem is to understand what could be a good substitute for it. The
elementary pieces
 of nontrivial transitive dynamics are the homoclinic classes of
 hyperbolic periodic points, which are exactly the basic sets  in  Smale
 theory.\break Actually, \cite{BD2} shows that, $C^1$-generically, two periodic points belong to the same transitive set if and only if
their two homoclinic classes  are equal
\footnote{\label{f.1}Recently, some substantial  progress has been made in understanding the elementary pieces of dynamics of 
nonhyperbolic diffeomorphisms.
In \cite{Ar1} and \cite{CMP} it is shown that,
for $C^1$-generic diffeomorphims or flows, any  homoclinic class is a maximal 
transitive set. Moreover, any pair of homoclinic classes is either equal or
disjoint. On the other hand, \cite{BD3} constructs examples of  $C^1$-locally
generic $3$-dimensional diffeomorphisms with maximal transitive sets
without periodic orbits.}. 
The hyperbolic-like properties of these homoclinic classes are the main subject of this paper.

Finally, 
we also see that some  of our arguments can be adapted almost straightforwardly   in the volume-preserving setting.
Let us now  state our results in a precise way.
 
\vglue4pt  {\it Statement of the results}.
Our first  theorem asserts
 that given any hyperbolic saddle $P$ its 
homoclinic class either admits an invariant dominated splitting or 
can be approximated  (by $\cC^1$-perturbations)
  by arbitrarily many sources or sinks.

\numbereddemo{Definition} \label{d.dominated}
Let $f$ be a diffeomorphism defined on a compact manifold $M$, 
$K$ an $f$-in\-va\-riant subset of $M$,
and $TM|_K = E\oplus F$ an $f_*$-invariant splitting of $TM$ over $K$, where the fibers $E_x$ of $E$ have constant dimension. 
We say that $E\oplus F$ is a {\it dominated splitting\/} (of $f$ over $K$)
if there exists $n\in \NN$ such that 
$$
\|f_*^n(x)|_E\|\,\|f_*^{-n}(f^n(x))|_F\|<1/2.
$$ 
We write $E\prec F$, or $E\prec_n F$ if we want to emphasize the role of
$n$, and we speak of $n$-dominated splitting.
\enddemo

Let us make two comments on this definition. 
First,
the invariant set $K$ is not supposed to be compact and the splitting is not supposed to be continuous. However, 
if $K$ admits a dominated splitting, it 
is always continuous and
can be extended uniquely to the closure $\bar K$ of  $K$
(these are classical results; for
details see Lemma \ref{l.dense} and Corollary \ref{c.limsup}).
 Moreover, the existence of a dominated splitting is equivalent to the existence of  some  continuous strictly-invariant cone field over $\bar K$;  this cone field  can be extended to some neighbourhood  $U$ of $\bar K$ and 
persists by $C^1$-perturbations. Thus there is an open neighbourhood of $f$ of diffeomorphisms
for which the maximal invariant set in $U$ admits a dominated splitting. In that sense,
the existence of a dominated splitting is a $C^1$-robust property.

Given a hyperbolic saddle $P$ 
of a diffeomorphism $f$ we denote by $H(P,f)$ the {\it homoclinic class\/} of $P$, i.e.\ the closure of the transverse intersections of the invariant manifolds of $P$. This set 
is transitive and the set
$\Si$ of hyperbolic periodic points $Q\in H(P,f)$ of $f$, whose stable and unstable manifolds intersect transversally the invariant manifolds of $P$, is dense in $H(P,f)$.

\specialnumber{1}\proclaim{Theorem} \label{t.main}
 Let $P$ be a hyperbolic saddle of a diffeomorphism $f$ 
defined on a compact manifold $M$. 
Then 
\begin{itemize}
\item either the homoclinic class $H(P,f)$
of $P$ admits a dominated splitting{\rm ,} 
\item or given any  neighbourhood $U$ of $H(P,f)$  and
any  $k\in \NN$ there is $g$ arbitrarily $C^1$\/{\rm -}\/close to $f$ having $k$ sources or sinks arbitrarily close to $P${\rm ,} whose
orbits are included in $U$. 
\end{itemize}

\endproclaim

If $P$ is a hyperbolic periodic point of $f$ then,  
for every $g$ $C^1$-close to $f$,
there is a hyperbolic periodic point $P_g$ of $g$ close to $P$
(this point is given by the implicit function theorem),
called the {\em continuation\/} of $P$ for $g$. From  Theorem~\ref{t.main} we get the following two corollaries. 

\proclaim{{C}orollary} \label{c.locdense}
Under the hypotheses of Theorem~{\rm \ref{t.main},} one of the
following two possibilities holds\/{\rm :}
\begin{itemize}
\item 
either there are a 
$C^1$\/{\rm -}\/neighbourhood $\cU$ of $f$ and a 
dense open subset $\cV\subset \cU$  such that  $H(P_g,g)$ has a dominated splitting for any $g\in \cV${\rm ,}
\item 
or there exist diffeomorphisms $g$ arbitrarily $C^1$\/{\rm -}\/close to $f$ such that 
$H(P_g,g)$ is contained in the closure of infinitely many
sinks or sources.
\end{itemize}

\endproclaim
 
\phantom{hi}
\vglue-.65in
\phantom{by}

\proclaim{{C}orollary} \label{c.residual}
There exists a residual subset $\cR$ of {\rm Diff}$^1(M)$ such that{\rm ,}
for every $f\in \cR$ and
any hyperbolic periodic saddle $P$ of $f${\rm ,} 
 the homoclinic class $H(P,f)$ satisfies one of the following possibilities\/{\rm :} 
\begin{itemize}
\item either $H(P,f)$ has a dominated splitting{\rm ,}
\item or $H(P,f)$ is 
included in the closure of the infinite set of sinks and sources of $f$. 
\end{itemize}

\endproclaim

\vglue-12pt
{\it Problem}. Is there a residual subset of {\rm Diff}$^1(M)$ of diffeomorphisms $f$ such that the homoclinic class of any hyperbolic
periodic point $P$ is the maximal transitive set containing $P$ (i.e.\ every transitive set containing $P$ is included in $H(P,f)$)?
Moreover,  when is $H(P,f)$ locally maximal?\footnote{Observe that the first part of the problem was
positively answered in \cite{Ar1} and \cite{CMP}  (recall footnote~\ref{f.1}).
Using these results,
\cite{Ab} shows that there is a $C^1$-residual set of diffeomorphisms 
such that the number of homoclinic classes is well defined and  
locally constant.
Moreover, if this number is finite, the homoclinic classes are locally maximal sets and there is a filtration whose levels correspond to homoclinic classes.} 
\vglue4pt

Actually, we prove a quantitative version of Theorem~\ref{t.main} relating the strength of the domination with the size of the perturbations of $f$
that we consider to get sinks or sources (see Proposition \ref{p.quantitative}).
This quantitative version is one of the keys 
for the next two results.

Note first that the creation of sinks or sources is not compatible with
the $C^1$-robust transitivity of a diffeomorphism. 
We apply Hayashi's connecting lemma (see \cite{Ha} and
Section~\ref{t.robust})
to get, by small perturbations, a dense homoclinic class
in the ambient manifold. Then using the quantitative version 
of Theorem~\ref{t.main} we show:

\specialnumber{2}\proclaim{Theorem}\label{t.robust}
Every $C^1$\/{\rm -}\/robustly transitive set {\rm (}\/or diffeomorphism\/{\rm )} 
admits a dominated splitting.
\endproclaim

Ma\~n\'e's theorem for surface diffeomorphisms
mentioned before gives a\break $C^1$-generic dichotomy between  hyperbolicity
and  the $C^1$-Newhouse\break phenomenon. It is now natural to ask what happens, in any 
dimension, {\it far from\/}  the  $C^1$-Newhouse phenomenon. 

\specialnumber{3}\proclaim{Theorem} \label{t.finite}
Let $f$ be a diffeomorphism such that the cardinal of the
set of sinks and sources is bounded in a $\cC^1$\/{\rm -}\/neigh\-bour\-hood of $f$. Then there exist
$l\in \NN$ and a $\cC^1$-neigh\-bour\-hood $\cV$ of $f$
such that{\rm ,} for any $g\in\cV$ and every periodic orbit $P$ of $g$ whose homoclinic class $H(P,g)$ is not trivial{\rm ,} 
 $H(P,g)$ admits an $l$\/{\rm -}\/dominated splitting. 
\endproclaim

A long term objective is to get 
 a spectral decomposition theorem in the nonuniformly hyperbolic case
for
diffeomorphisms far from the  Newhouse phenomenon.
 Having this goal in mind, we can reformulate Theorem \ref{t.finite} as follows: 
\vglue9pt
 
{\it Under the hypotheses of  Theorem {\rm \ref{t.finite},} 
for every  $g$ sufficiently $C^1$\/{\rm -}\/close 
to $f$ there are compact invariant sets $\La_i(g)${\rm ,} 
$i\in\{1,\dots,\dim(M)-1\}${\rm ,} such that\/{\rm :}
\begin{itemize}
\item 
Every
$\La_i(g)$ admits an $l$\/{\rm -}\/dominated splitting $E_i(g)\prec_l F_i(g)$ with\break $\dim (E_i(g)) =i${\rm ,}
\item
every  nontrivial homoclinic class $H(Q,g)$ is contained in some  $\La_i(g)$. 
\end{itemize} 
}
\vglue9pt

Unfortunately,
this result has two disadvantages.
First, the $\La_i(g)$ are supposed to be 
neither transitive nor disjoint. Moreover, the nonwandering
 set $\Omega(g)$ is not {\em a priori} contained in 
the union of the $\La_i(g)$ (but every homoclinic class of a periodic orbit in $(\Omega(g)\setminus\bigcup_i\La_i(g))$ is trivial). So
we are still far away from a completely satisfactory 
spectral decomposition theorem\footnote{Fortunately, the results in 
  footnote~2 gave a spectral decomposition for generic 
diffeomorphisms with finitely many homoclinic classes.}. 
In view of these comments the following problem
arises in a natural way.

\demo{Problem} 
Let $\cU\subset$ {\rm Diff}$^1(M)$ be an open 
set of diffeomorphisms 
for which
 the number of sinks and sources is 
uniformly bounded. Is there some subset $\cU_0\subset \cU$,
either dense and open or residual
in $\cU$, of diffeomorphisms $g$ 
such that $\Omega(g)$ is the union of finitely many disjoint compact sets $\La_i(g)$ having a dominated splitting? 
\enddemo

Let us observe that the
 Newhouse phenomenon is  not incompatible with the existence of a dominated
splitting if we do not have 
any  additional\break information on the action of $f_*$ on the subbundles of the splitting. 
Actually, using Ma\~n\'e's ergodic closing lemma (see \cite{Ma3}) 
we will get some control of the\break
action of the  derivative $f_*$ on the volume induced on the subbundles. For that we  need to introduce
dominated splittings having more than two bundles.
An invariant splitting $TM|_K = E_1\oplus\cdots\oplus E_k$ is {\it dominated\/}
 if $\bigoplus_1^j E_i \prec \bigoplus_{j+1}^k E_i$
for every $j$. In this case we use the notation
$E_1\prec E_2 \prec \cdots \prec E_k$. 

By Proposition~\ref{p.thefinest},  there is a unique 
dominated splitting, called {\it finest dominated splitting,\/} such that any  dominated splitting is obtained by regrouping its subbundles by packages corresponding to intervals.

\specialnumber{4}\proclaim{Theorem} \label{t.jac}
Let $\La_f(U)$ 
be a $C^1$\/{\rm -}\/robustly transitive set
and $E_1\oplus\cdots\oplus E_k$, $E_1\prec\cdots \prec E_k${\rm ,}
 be its
finest dominated splitting. Then there exists $n\in \NN$ such that
$(f_*)^n$ contracts uniformly the volume in $E_1$ and expands uniformly the volume in $E_k$.
\endproclaim

This result synthesizes  previous results 
in lower dimensions of \cite{Ma3} and \cite{DPU} on robustly transitive diffeomorphisms (or sets) and it shows that, in the list of robustly 
transitive diffeomorphisms
above, each example corresponds to the worst pathological case in the corresponding dimension.
Observe  that if $E_1$ or $E_k$ has dimension one, then it is uniformly hyperbolic (contracting or expanding). Then, for robustly 
transitive diffeomorphisms, we have:
\begin{itemize}
\item  
In dimension $2$ the dominated splitting has necessarily two $1$-dimensional bundles, so  that the diffeomorphism is hyperbolic
and then Anosov\break (Ma\~n\'e's result above).
\item  
In dimension $3$ at least one of the bundles has dimension $1$ and so it is hyperbolic and the diffeomorphism is partially hyperbolic (see \cite{DPU}). In this dimension, the finest decomposition can contain {\em a priori\/} two or three bundles and in the list above there are examples 
of both of these possibilities.   
\item 
In higher dimensions the extremal subbundles may have dimensions  strictly bigger than one and so they are not necessarily
hyperbolic: This is exactly what happens in the examples in \cite{BoVi}.
\end{itemize}

Theorem~4 motivates us to introduce the notions of {\it volume hyperbolicity\/} and {\it volume partial
hyperbolicity,\/} as the existence of dominated splittings, say $E\prec G$ and $E\prec F\prec G$, respectively, 
such that the volume is uniformly contracted on the  bundle $E$ and expanded on  $G$. 
We think that this notion could be the best substitute for the hyperbolicity in a nonuniformly hyperbolic context.   

The volume partial hyperbolicity 
is clearly incompatible with the existence of sources or sinks. 
However, in the proof of  Theorem~4 ,
at least for the moment, we need  the robust transitivity to 
obtain the partial volume hyperbolicity. 
 Bearing in mind this comment and  our previous results, let us pose some questions:

\demo{Problems} 1.  In Theorem~\ref{t.main}, is it possible to replace the notion of dominated splitting by the notion of volume
partial hyperbolicity?\footnote{In this direction, using the techniques in
this paper, \cite{Ab} shows the volume hyperbolicity of  
the homoclinic classes of 
generic diffeomorphisms
having finitely many homoclinic classes.}

\vglue8pt 2. 
 Is the notion of volume hyperbolicity (or volume partial hyperbolicity) sufficient to assure the generic existence of 
finitely many Sinai-Ruelle-Bowen (SRB) measures whose basins cover a total Lebesgue measure set? More \pagebreak precisely:
\begin{quote}
Let $f$ be a $C^1$-robustly transitive diffeomorphism of class $C^2$ on a compact manifold $M$. Does there exist $g$
close to $f$ having  finitely many SRB measures 
such that the union 
of their basins has total Lebesgue measure in $M$? 
\end{quote}
\enddemo

\phantom{now}
\vglue-.45in 

For ergodic properties of partially hyperbolic systems 
(mainly existence of SRB measures) we refer the reader to
\cite{BP}, \cite{BoVi}, and \cite{ABV}.  

Let us observe that in the measure-preserving setting
(also volume-pre\-serving) the notion of stably ergodicity
 (at least in the case of $C^2$-diffeo\-morphisms)
seems to play the same role as  the robust transitivity
in the topological setting.
See the results in \cite{GPS} and \cite{PgSh} which, in
rough terms, show that weak forms
of hyperbolicity may be necessary
for stable ergodicity and 
go a long way in guaranteeing it. Actually, in \cite{PgSh}
it is conjectured that stably ergodic diffeomorphisms are
open and dense among the partially hyperbolic 
$C^2$-volume-preserving diffeomorphisms.
See \cite{BPSW} for a survey on stable ergodicity
and \cite{DW} for recent progress on the previous conjecture.
Our next objective is precisely the study of 
$C^1$-volume-preserving diffeomorphisms.

Although this paper is not devoted to conservative diffeomorphisms some of our results have a quite straightforward generalization into the conservative 
context. This means that the manifold is endowed with a smooth volume form~$\omega$; then 
we can speak of {\it conservative\/} (i.e.\ volume-preserving)
diffeomorphisms. We denote by Diff$^1_\omega(M)$ the set of $C^1$-conservative diffeomorphisms.  

A first challenge is to get a suitable version of the 
generic spectral decomposition theorem by Ma\~n\'e 
(dichotomy between hyperbolicity and the Newhouse phenomenon) for 
conservative diffeomorphisms. Obviously, since conservative diffeomorphisms have
 neither sinks nor sources, one needs to replace sinks and sources by {\it elliptic points\/} (i.e.\ periodic points whose derivatives
have some eigenvalue of modulus one).
 Very little is known in this context. First,
there is an unpublished result by Ma\~n\'e (see \cite{Ma1}) 
which says that $C^1$-generically, area-preserving diffeomorphisms
of compact surfaces are either Anosov or have Lyapunov exponents
equal to zero for almost every orbit 
(see also \cite{Ma4} for an outline of a 
possible
proof).\footnote{Recently, \cite{Bc} gave   a complete proof of this 
result. For a 
generalization to higher dimensions, see \cite{BcVi1}.}    
Ma\~n\'e also announced a version of his theorem in higher
dimensions for symplectic diffeomorphisms; see \cite{Ma1}.\footnote{See \cite{Ar2} for progress on this subject in dimension four.} 
Unfortunately, as far as we know, there are no available 
complete proofs of these results. See also the results by Newhouse in \cite{N2} where 
he states a dichotomy between hyperbolicity (Anosov diffeomorphisms) and existence of elliptic periodic points.

Related to  the announced  results of Ma\~n\'e,
there is the following question
in \cite{He}:

\nonumproclaim{Conjecture {\rm (Herman)}}
 Let $f\in$ {\rm Diff}$^1_\omega(M)$ be a conservative diffeomorphism of a compact manifold $M$. 
Assume that there is a neighbourhood $\cU$ of $f$ in {\rm Diff}$^1_\omega(M)$ 
such that for any $g\in \cU$ and  every periodic orbit $x$ of $g$ the matrix $g_*^n(x)$ {\rm (}\/where $n$ is the period of $x${\rm )}
has at least one eigenvalue of modulus different from one.  Then $f$ admits a dominated splitting.
\endproclaim

The following results give partial answers to this question:
 
\specialnumber{5}\proclaim{Theorem}\label{t.herman}
Let $f\in$ {\rm Diff}$^1_\omega(M)$ be a conservative diffeomorphism  of
an $N$\/{\rm -}\/dimensional manifold $M$. Then there is $l\in \NN$ such that{\rm ,}
\begin{itemize} 
\item  
either there is a conservative $\varepsilon-C^1$\/{\rm -}\/perturbation $g$ of $f$ having a periodic
 point $x$ of period $n\in\NN$ such that $g_*^n(x)={\rm Id}${\rm ,} 
\item 
 or for any conservative diffeomorphism $g$
 $\varepsilon-C^1$\/{\rm -}\/close to 
$f$ and every periodic saddle $x$ of $g$ the homoclinic class $H(x,g)$
admits an $l$\/{\rm -}\/dominated splitting.
\end{itemize}

\endproclaim

\specialnumber{6}\proclaim{Theorem}\label{t.con}
Let $f$ be a conservative diffeomorphism defined on a compact $N$\/{\rm -}\/dimensional manifold.
Then there are two possibilities\/{\rm :}\/ 
\begin{itemize}
\item Either given any $k\in\NN$ there is a conservative diffeomorphism $g$ arbitrarily
$C^1$\/{\rm -}\/close to $f$ having $k$ periodic orbits whose derivatives are the identity.
\item Or the manifold $M$ is the union of finitely many {\rm (}\/less than $N-1${\rm )} invariant compact
{\rm (}\/a priori  {\rm nondisjoint)} sets  having a dominated splitting. \end{itemize}

\endproclaim

Observe that if in  Theorem~6 above the diffeomorphism $f$ is transitive and the second possibility of the dichotomy occurs, then one of the invariant compact sets has to be the whole manifold (one of them contains a dense orbit). This means that, in the transitive case, Theorem~6 gives a complete positive answer to Herman's conjecture:

 \proclaim{{C}orollary}
\label{c.tran} Let $f\in$ {\rm Diff}$^1_\omega(M)$ be a conservative
transitive diffeomorphism of a manifold $M$. 
Assume that there is a neighbourhood $\cU$ of $f$ in {\rm Diff}$^1_\omega(M)$ 
such that for any $g\in \cU$ and  every periodic orbit $x$ of $g$ the matrix $g_*^n(x)$ {\rm (}\/where $n$ is the period of $x${\rm )}
has at least one eigenvalue of modulus different from one.  Then $f$ admits a dominated splitting.
\endproclaim

Let us observe that if $f$ is transitive and   there is some periodic
point $x$ of $f$ such that
$f_*^n(x)={\rm Id}$, $n$ is the period of $x$,
 then given any $\varepsilon>0$
there is a  $C^1$-perturbation $g\in$ {\rm Diff}$^1_\omega (M)$
of $f$ such that its totally elliptic points (derivative equal to the identity) are
$\varepsilon$-dense in $M$.

A conservative diffeomorphism $f\in$ Diff$^1_\omega(M)$
is {\em robustly transitive in  {\rm Diff}$^1_\omega(M)$\/}
if there is $\varepsilon >0$ such that every 
$\varepsilon-$ perturbation $g\in$ Diff$^1_\omega(M)$ of $f$
is transitive.
Observe that {\em a
priori\/} the robust transitivity in  Diff$^1_\omega(M)$ does
not imply the robust transitivity in Diff$^1(M)$.

\vglue8pt {\elevensc Conjecture}.
Let $f$ be a robustly transitive diffeomorphism in  {\rm Diff}$^1_\omega(M)$.
Then $f$ admits 
a nontrivial dominated splitting defined on the whole of  $M$.
\vglue8pt

In view of Corollary~\ref{c.tran}, to prove this conjecture one needs 
to show that a robustly transitive diffeomorphism $f$
cannot have  periodic points $x$
whose derivative $f_*^n(x)$ is the identity.

Finally, we  observe that the control of the volume in the subbundles is almost
 straightforward for conservative systems: 

\vglue8pt {\elevensc Proposition 0.5.}
{\it Let $f$ be a conservative diffeomorphism  and $E\oplus F${\rm ,} 
$E\prec F${\rm ,} be a 
dominated splitting of $\,TM$. Then $f_*$ contracts uniformly the volume
in $E$ and expands uniformly the volume
in $F$.
}
\vglue9pt

{\it Main ideas of the proofs}.
 Ma\~n\'e's paper \cite{Ma3} combines two main ingredients: systems
of matrices and the ergodic closing lemma. He first  considers
the linear maps induced by the derivative of a diffeomorphism $f$ over the orbits of its
 periodic points, thus obtaining a {\em system of matrices.\/}
 He shows that (in his context) a system of matrices admits a dominated splitting  if it is not possible to perturb it to get a matrix  with some eigenvalue of modulus one.
By a lemma of Franks, see \cite{F} and Section~\ref{s.ls},
 each 
perturbation of the system of matrices over a finite number of periodic orbits
 corresponds to a $C^1$-perturbation of $f$ and vice versa. Hence the existence of a
dominated splitting also holds for $C^1$-diffeomorphisms.
Finally, to get the uniform expansion and contraction on the subbundles of the
 splitting he uses his ergodic closing lemma (see \cite{Ma3}). 

 \pagegoal=50pc

Our proof uses these two tools introduced by Ma\~n\'e.
 Using  Franks' lemma we translate the problem of the existence of a dominated splitting for  diffeomorphisms into the same
problem for abstract linear systems. However,  the systems of matrices in \cite{Ma3} do not contain one relevant  dynamical
information  about $f$ that we  need. 
Actually, the solution of this difficulty
is probably the subtlest point of our arguments, so let us be somewhat more precise:

On one hand, in the context of \cite{Ma3}, all the periodic points have the same {\em index\/} (dimension of the stable bundle); thus
the system of matrices has a natural splitting (the one corresponding to the stable and unstable bundles of~$f_*$).
Then if this splitting is not dominated one gets a perturbation of it
having one eigenvalue of modulus $1$. On the other hand,
in our case there are points having different indices.
Moreover, points having eigenvalues of modulus~$1$ are not forbidden.
So we need some extra arguments to conclude our proof. In fact,
we need to control {\em all\/} the eigenvalues to create sources or   sinks. 

The additional argument, that comes from the dynamics, is a  property of our linear systems called {\em transitions.\/} Given two periodic points $P$ and $Q$ in the same 
{\em homoclinic class\/}
(i.e.\ their invariant manifolds intersect transversally) there are periodic orbits passing  first arbitrarily close to $P$, and thereafter 
arbitrarily  close to $Q$, and so on. These orbits can be chosen upon
  arbitrary sequences of times
 (the orbit spends $k_1$ iterates close to $P$ then, after a bounded number of iterates,
it becomes close to $Q$ and remains $k_2$ iterates close to $Q$, and so
on). So we define a structure we call transitions which translates this
dynamical behavior into the world of the abstract linear systems. This
property allows us to consider  the product of  matrices of the system
corresponding to different orbits as a matrix of the system. In fact,
the transitions endow the linear system with a ``semigroup-like" structure. Clearly,
this is not the case for general linear systems.
\pagegoal=48pc

\vglue4pt
Finally, after we introduce the linear systems with transitions, the proof of the existence of a dominated splitting involves
only arguments of linear algebra.  Precisely, this algebraic approach has allowed us to improve previous results
by stating them in higher dimensions and by eliminating the robust transitivity hypothesis.

 \vglue2pt
The problem of the  existence of points with different indices  already
appeared in \cite{DPU}, where  it was solved by considering only robustly transitive sets;
thus any perturbation of the dynamics  remains transitive.
This additional hypothesis in \cite{DPU} allows us to jump from 
the dynamical world to the abstract linear world, here  
do some perturbation, and then jump back to the dynamical world to do a new perturbation,
and so on.  
In our context we have no control of the variation of a homoclinic class
after dynamical perturbations. So it is crucial for Theorem~\ref{t.main}
that all the perturbations we do ``live in the world of
abstract  linear systems" and do not modify the underlying dynamics (that is here possible because Ma\~n\'e's linear systems have been enriched with the transitions).
\vglue2pt
In our proof, assuming that there is no
 dominated splitting, we perform 
a series of perturbations of the linear system; as
a   final result of such perturbations we get a
linear system  having a homothety. It is only then that we realize this linear system as a diffeomorphism using  Franks' lemma,
and the point corresponding to the homothety becomes a sink or a source of the diffeomorphism.

\vglue4pt
Finally, for the control of the volume in the extremal subbundles  (Theorem~4) we use the ergodic closing lemma, which gives a
dynamical perturbation having a periodic point reflecting the lack of volume expansion or contraction of the bundles. Unfortunately,
without  any additional hypothesis, this point has {\em a priori\/} nothing to do with the initial homoclinic class. This explains why
Theorem~4 only holds for robustly transitive systems.         
\pagebreak

{\it Acknowledgments.} We thank IMPA (Rio de Janeiro) and
the  Laboratoire de Topologie (Dijon) for their warm hospitality 
during our visits while preparing this paper. We also want to express our gratitude to Marcelo Viana for his encouragement and
conversations  on this subject, to Floris Takens and Marco Brunella for their enlightening explanations about perturbations of
conservative systems, to Flavio Abedenur, Marie-Claude Arnaud, Jairo Bochi, Michel Herman and Gioia Vago,  
for their careful reading of this paper, and to the 
students of the Dynamical Systems Seminar of IMPA for many comments on the 
first version of this paper. 
Last, but not least, we should like to thank 
the late Michel Herman 
for his interest in the present work and his continuous encouragement.
 
\vglue12pt \centerline{\bf Contents}
\vglue9pt
\def\sni#1{\smallbreak\noindent{#1}. }
\def\ssni#1{\vglue-1pt\noindent\hskip18pt {#1}.}
\def\vni#1{\vglue-1pt\noindent \hskip40pt {#1}.}

\sni{1} Linear systems with transitions
\ssni{1.1} Linear systems: Topology and linear changes of coordinates
\ssni{1.2} Special linear systems
\ssni{1.3} Dominated splittings
\ssni{1.4} Periodic linear systems with transitions

\sni{2} Quantitative results: Proofs of the theorems 
\ssni{2.2} Proofs of the theorems
\vni{2.2.1} Proofs of Theorem 1 and Proposition 2.6
\vni{2.2.2} Proof of Corollary 0.2
\vni{2.2.3} Proof of Corollary 0.3
\vni{2.2.4} Proof of Theorem 2
\vni{2.2.5} Proof of Theorem 3
\sni{3} Two-dimensional linear systems
\ssni{3.1} Proof of Proposition 3.1
\sni{4} Invariant subbundles: Reduction of the dimension and the finest\hfill\break
\phantom{4. }\hskip5pt dominated splitting
\ssni{4.1} Quotient of linear systems and restriction to subbundles
\ssni{4.2} The finest dominated splitting
\ssni{4.3} Transitions and invariant spaces
\ssni{4.4} Diagonalizable systems
\sni{5} Dominated splittings, complex eigenvalues of rank $(i,i+1)$, and\hfill\break \phantom{4. }\hskip5pt homotheties
\ssni{5.1} Getting complex eigenvalues of any rank
\ssni{5.2} End of the proof of Proposition 2.4
\ssni{5.3} Proof of Proposition 2.5
\vni{5.3.1} End of the proof of Proposition 2.5
\vni{5.3.2} Proof of Lemma 5.4 (and Remark 5.5)
\sni{6} Finest dominated splitting and control of the jacobian in the extremal\hfill\break \phantom{4. }\hskip5pt  bundles: Proof
of Theorem 4
\ssni{6.1} Control of the jacobian over periodic points
\ssni{6.2} Ma\~n\'e's ergodic closing lemma: Proof of Proposition 6.2

\sni{7} The conservative case
\ssni{7.1} Proof of Theorem~6
\ssni{7.2} Volume properties of dominated splittings of conservative systems

\section{Linear systems with transitions }
\label{s.ls}

Let  $f$ be a diffeomorphism. By Franks' lemma below (see for instance
\cite{F}), to any perturbation $A$ of the derivative $f_*$ along the
orbits of finitely many periodic points corresponds a diffeomorphism
$g$, $C^1$-close to $f$, such that $g_*=A$ along these orbits. 
This lemma allows us to consider perturbations of the derivative $f_*$ keeping unchanged the dynamics of $f$, in order to get a suitable derivative along some periodic orbits.  The aim of this section is to define
 in details the framework ({\it periodic linear systems\/}) which gives a precise meaning of this kind of perturbations, and to
translate  into this language the dynamical properties that we will need (specially the notion of {\it transitions,\/} see 
Definition~\ref{d.epsilontransition}). Finally, we prove that the homoclinic classes define a  periodic linear system with transitions
(Lemma~\ref{l.Df}) and we state an easy (but typical) consequence of the existence of transitions (Lemma~\ref{l.linearsources}).

Before beginning this section let us state precisely   Franks' lemma: 
\nonumproclaim{Lemma {\rm (Franks)}}
Suppose the  $E$ is a finite set and $B$ is an
$\varepsilon$\/{\rm -}\/perturbation of $f_*$ along $E${\rm .} Then there is 
a diffeomorphism $g$ $\varepsilon$-$\cC^1$\/{\rm -}\/close to $f${\rm ,} 
coinciding with $f$
out of an arbitrarily  small neighbourhood of $E${\rm ,} 
equal to  $f$ in $E${\rm ,} and such that $g_*$ coincides with $B$ in $E$.
\endproclaim

Let us point out that Franks' lemma is the key which allows  us to translate results on linear systems to the  dynamical context and
it will  often  be used in this paper.

\demo{{\rm 1.1.} Linear systems\/{\rm :} Topology and linear changes of coordinates}
Let $\Si$ be a topological space and $f$ a
homeomorphism defined on $\Si$. Consider 
a locally trivial vector bundle (of finite dimension) $\cE$ over $\Si$. Denote by $\cE_x$ the fiber of $\cE$ at $x\in \Si$.  
We assume that the dimension of the fibers $\cE_x$, $\dim (\cE_x)$, 
does not depend on $x\in\Si$. In what follows, we call
this number {\it dimension of the bundle $\cE$,\/} denoted by $\dim(\cE)$. 

A  {\it euclidian metric $|\cdot|$ on the bundle $\cE$} is a collection of euclidian metrics on
the fibers $\cE_x$, $x\in \Si$, {\em a priori} not depending continuously on $x$. 

We denote by $\cG\cL(\Si,f,\cE)$ the set of maps $A\colon \cE\to \cE$ such that for every $x\in \Si$ 
the induced map $A(x,\cdot)$ is a linear isomorphism from $\cE_x\to\cE_{f(x)}$, thus $A(x,\cdot)$ belongs to $\cL(\cE_x,\cE_{f(x)})$ and is invertible. For each
$x\in \Si$ the euclidian metrics on $\cE_x$ and $\cE_{f(x)}$ induce a norm (always  denoted by $|\cdot|$) on 
$\cL(\cE_x,\cE_{f(x)})$:
$$
|B(x,\cdot)| = \sup\{|B(x,v)|,   v\in\cE_x, |v|=1  \}.
$$
 
Let now $A\in\cG\cL(\Si,f,\cE)$ and define 
$|A| = \sup_{x\in\Si}|A(x,\cdot)|$. 
Observe that, for any $A\in\cG\cL(\Si,f,\cE)$, its inverse  $A^{-1}$ belongs to $\cG\cL(\Si,f^{-1},\cE)$. So we can define $|A^{-1}|$ in
the same way. Finally,  the {\em norm\/} of a $A\in\cG\cL(\Si,f,\cE)$
 is  $\|A\| = \sup\{|A|,|A^{-1}|\}$.
\enddemo

\numbereddemo{Definition} 
A {\it linear system\/}\footnote{After writing this paper, we realized 
that this notion corresponds to the classical concept of 
{\it linear cocycle\/} over the homeomorphism $f$.} is a $4$-uple $(\Si,f,\cE,A)$ where $\Si$ is a topological space, $f$ is a
homeomorphism of $\Si$,  $\cE$ is a  euclidean bundle over $\Si$,  $A$ belongs to $\cG\cL(\Si,f,\cE)$, and 
$\|A\| <\infty$. 
\enddemo

In what follows, for the sake of simplicity,
we  sometimes  denote by $A$ a linear system $(\Si,f,\cE,A)$
if there is no ambiguity on $\Si$, $f$, and $\cE$. 

\demo{Example {\rm 1}}
Let $f$ be a diffeomorphism defined on a riemannian manifold $M$
and $\Si\subset M$ an $f$-invariant subset. Consider the restriction  to $\Si$  of the tangent bundle,
$\cE= TM|_{\Si}$.
The riemannian metric on $M$ induces a euclidean structure on $\cE$. 
Then $(\Si,f|_{\Si},\cE,f_*|_{\cE})$ is the 
{\em natural linear system induced by $f$ over $\Si$.\/}
\enddemo

We denote by $\cG\cL^{\infty}(\Si,f,\cE)$ the space of linear systems over
$(\Si,f,\cE)$ such that $||A||<\infty$ is endowed with the distance defined by
$$
d(A,B)= \sup\{|A-B|,|A^{-1}-B^{-1}|\},
\quad
A,\,B\in\cG\cL^{\infty}(\Si,f,\cE).
$$
We can now define an
{\em $\varepsilon$-perturbation} of $A$ as a linear system $\tilde A$, defined over
$(\Si,f,\cE)$, such that $d(A,\tilde A)<\varepsilon$.

Very elementary arguments of linear algebra show that any perturbation of a linear system can be obtained by composing it with linear maps close to 
the identity. More precisely,
let $ A\in\cG\cL^{\infty}(\Si,f,\cE)$ and consider
some linear system $E\in \cG\cL^{\infty}(\Si,{\rm Id}_{\Si},\cE)$. Then $E\circ A$ and $A\circ E$
(defined in the obvious way) belong to
$\cG\cL^{\infty}(\Si,f,\cE)$. Moreover, if $E$ is close to 
the identity linear system $(\Si, {\rm Id}_{\Si},\cE, {\rm Id}_{\cE})$, then $E\circ A$ and $A\circ E$ are also close to $A$. 

Consider now some change of the euclidean metrics on the fibers. Assume  that the matrices
of the changes of coordinates (from an orthonormal basis of the initial metric to an orthonormal basis of the new metric) and their inverses are uniformly bounded on $\Si$. Then 
every linear system in the initial metric
induces a new system (for the new metric). Moreover,
this change of  metrics keeps invariant the
 topology of the set of linear systems. Let us be a little bit more precise.

Let $\cE$ denote a euclidean bundle on a topological space $\Si$ endowed with the euclidean metric $|\cdot|$.
Denote  by $\cE_1$ the same bundle, but now endowed with a
different euclidean metric $|\cdot|_1$. Denote by $P\colon \cE\to\cE_1$ the identity map considered as a morphism of bundles.
Using the metrics $|\cdot|$ and $|\cdot|_1$ we can define the norms $|P|$ and $|P^{-1}|$. 
Write $\|P\| = \sup\{|P|,|P^{-1}|\}$. If $\|P\|<\infty$, the 
canonical bijection ${\rm Id} \colon \cG\cL(\Si,f,\cE)\to \cG\cL(\Si,f,\cE_1)$ induces a  homeomorphism
 from $\cG\cL^{\infty}(\Si,f,\cE)$ (with the distance $d$) to $\cG\cL^{\infty}(\Si,f,\cE_1)$
(with the corresponding distance $d_1$). 
These two simple facts are put together in the following lemma. 

\proclaim{Lemma} \label{l.changes}
\hskip-10pt
{\rm 1.} 
Given $K>0$ and $\varepsilon>0$ 
there is $\delta>0$ such that for any linear system 
$(\Si,f,\cE, A)${\rm ,} $A\in\cG\cL^{\infty}(\Si,f,\cE)$ and
$\|A\| < K${\rm , } 
and every\break
$\delta$\/{\rm -}\/perturbation of the identity $(\Si,{\rm id}_\Si,\cE,E)${\rm ,} 
  $E\circ A$ and $A\circ E$ are 
$\varepsilon$\/{\rm -}\/perturbations of $A$.
\vglue8pt
  {\rm 2.} For every  $K>0${\rm ,} $K_0>0${\rm ,} and $\varepsilon>0$
there are  $K_1>0${\rm ,}  $\delta>0$ satisfying
the following property\/{\rm :}
\vglue8pt
 Consider a pair of euclidean bundles
$\cE$ and $\cE_1$   over $\Si$ endowed with the metrics $|\cdot|$ and $|\cdot|_1${\rm ,} and the isomorphism of bundles
$P\colon \cE\to \cE_1$  
induced by the identity on $\Si$
{\rm (}\/i.e.\
given $x\in\Si$ the map $P(x,\cdot)$
is a linear isomorphism from $\cE_x$ to $(\cE_1)_x${\rm ).} 
Assume that $\|P\|< K_0$. 
Let $(\Si,f,\cE,A)$ be a linear system  such that $\|A\|$ 
is bounded by $K$.
 \vglue8pt
 Let $B=P\circ A\circ P^{-1}${\rm ,}
then
$(\Si, f, \cE_1,B )$ is a linear system bounded by $K_1$.
Moreover{\rm ,}  any $\delta$\/{\rm -}\/perturbation of $B$ is conjugate
by $P$ to some $\varepsilon$\/{\rm -}\/perturbation of~$A$.  
\endproclaim

Let $(\Si,f, \cE, A)$ be a linear system and $n\in \NN$. 
The {\em $n$-th iterate of $A$,\/} denoted by $A^{(n)}$,
is the linear system over $(\Si,f^n,\cE)$
defined by $A^{(n)}(x) = A(f^{n-1}(x))\circ\cdots\circ A(f(x))\circ A(x)$.

Consider an $f$-invariant subset $\Si'$ of $\Si$ and the
restriction of the linear bundle $\cE$ to $\Si'$, then  $A$ induces canonically a linear system over $(\Si',f|_{\Si'},\cE|_{\Si'})$ called
the {\em linear subsystem induced by $A$ over $\Si'$.\/}

\vglue12pt 1.2. {\it Special linear systems}.
Along this work, the linear systems we  consider will  often  be endowed with
some additional structures: In some
cases they are {\em continuous,\/}
and most of them
are {\em periodic.\/}
We also consider {\em systems of matrices.\/} 
Finally, the most important additional structure
we will introduce 
is the notion of {\em transitions.\/} 
Let us now present the three first quite natural structures.
Due to
its specific and  subtle nature we postpone to the next
paragraph the notion of transition.  This key definition will deserve
special attention and \pagebreak care. 
\enddemo

In the sequel,  $\Si$ is  a topological space, $f$ a
homeomorphism of $\Si$, and $\cE$ a locally trivial vector bundle over $\Si$ endowed with a euclidean metric 
$|\cdot|$ on the fibers.
       A linear system $(\Si,f,\cE,A)$
is  {\em continuous\/}
 if the euclidean structure on the fibers varies continuously and the function $A\colon \cE\to\cE$ is continuous.

The linear system $(\Si,f,\cE,A)$ is
 {\em periodic\/} if all the orbits of $f$ are periodic.
In this case we let 
$M_A(x)\colon \cE_x\to \cE_x$  be the product of the 
$A(f^i(x))$ along the orbit of $x$. More precisely, 
let $p(x)$ be the period of $x\in\Si$, then
$$
M_A(x) = A(f^{p(x)-1}(x))\circ\cdots\circ A(x)= A^{(p(x))}(x).
$$  

Finally, $(\Si,f,\cE,A)$ is a
{\em system of matrices\/} 
if the euclidean bundle $\cE$ is the trivial bundle 
$\Si\times \RR^N$, where $\RR^N$ is endowed with the canonical euclidean metric. 
In this case every linear map 
$A(x)$ is canonically identified with an element of ${\rm GL}(N,\RR)$.

Let $(\Si,f,\cE, A)$ be an ({\it a priori \/} noncontinuous) linear system.  It will  sometimes be useful to fix an orthonormal basis on
each fiber $\cE_x$ (this basis does not depend, in general, continuously on the point $x\in\Si$).
These bases give an (\/{\it a~priori} noncontinuous) trivialization
 of the Euclidean bundle 
$\cE$. So in these new coordinates $A$ can be considered as a system of matrices.
Two systems of matrices {\em define the same linear system\/} if
at each point 
 there exists an orthonormal change of coordinates conjugating the two
systems.

\demo{{\rm 1.3.} Dominated splittings}
The definition of dominated splitting for  an invariant set  
of a diffeomorphism (see Definition~\ref{d.dominated}) can be directly generalized for linear systems
as follows. 
Let $(\Si,f,\cE,A)$ be a linear system, an {\em invariant subbundle\/}  is a collection of linear subspaces $F(x)\subset \cE_x$
whose dimensions do not depend on $x$
and such that $A(F(x))=F(f(x))$. 
An {\em $A$-invariant splitting\/} $F\oplus G$ is given by two invariant subbundles such that $\cE_x=F(x)\oplus G(x)$ at each $x\in\Si$.  
\enddemo

\numbereddemo{Definition} 
\label{d.dominated2}
Let $(\Si,f,\cE,A)$ be a linear system 
and $\cE = F\oplus G$ an $A$-invariant splitting. 
We say that $F\oplus G$ is a {\em dominated splitting\/} if there exists $n\in \NN$ such that 
$$
\|A^{(n)}(x)|_F\|\,\|A^{(-n)}(f^n(x))|_G\|<1/2
$$
for every $x\in \Sigma$. We write $F\prec G$. 

If we want to emphasize the role of $n$ then
 we say that $F\oplus G$ is an {\em $n$-dominated splitting} and  write $F\prec_n G$ .

Finally, the {\em dimension of
the dominated splitting} is the dimension 
of the subbundle $F$.
\enddemo

Suppose now that $(\Si,f,\cE,A)$ is a {\it continuous} linear system, then any dominated splitting can be obtained 
by considering  subsystems induced by $A$ over dense subsets $\Si'\subset \Si$. More \pagebreak precisely,

\proclaim{Lemma} \label{l.dense}
Let $(\Si, f, \cE, A)$ be a continuous linear system such
that there is a dense $f$\/{\rm -}\/invariant subset $\Si_1\subset\Si$ whose corresponding linear subsystem admits an
$l$\/{\rm -}\/dominated splitting. Then $(\Si, f, \cE, A)$ admits an $l$\/{\rm -}\/dominated splitting.

More generally{\rm ,} suppose
that there is a sequence of {\rm (}\/not necessarily continuous\/{\rm )} systems $(\Si,f,\cE,A_k)$ converging to $(\Si,f,\cE,A)$ 
such that for every $k$ there is a dense invariant subset $\Si_k\subset
\Si$ where $A_k$ admits an $l$\/{\rm -}\/dominated splitting. Then 
$A$
admits an $l$\/{\rm -}\/dominated splitting in the whole $\Si$.

Finally{\rm ,} any dominated splitting of a continuous linear system is continuous.
\endproclaim

\demo{Proof} Given $x\in \Si$ 
consider a sequence $(x_k)$, $x_k\in\Si_k$, converging to $x$. 
For a fixed $k$ we have an
$l$-dominated splitting $E_k\oplus F_k$.  Taking a subsequence 
we can assume
 that the dimensions of these spaces are independent of $k$ and that
the sequences $E_k(x_k)$ and $F_k(x_k)$ converge to some subspaces $E(x)$ and $F(x)$. 

By definition of $l$-dominance, given  any $k$, 
$u_k\in E_k(x_k)$, and $v_k\in F_k(x_k)$, we have 
$$
2\, \frac{\|A^l_k(u_k)\|}{\|u_k\|} \leq \frac{\|A^l_k(v_k)\|}{\|v_k\|}.
$$ 
By the continuity of $A$ and the convergences of $A_k \to A$, $x_k \to x$, $E_k(x_k) \to E(x)$, and $F_k(x_k) \to F(x)$,
we get
$$
2\,\frac{\|A^l(u)\|}{\|u\|} \leq \frac{\|A^l(v)\|}{\|v\|}
$$
for every $u\in E(x)$ and $v\in F(x)$. So these two spaces are transverse.

Finally, it remains to check that these two spaces are uniquely defined and give
 an invariant splitting. Observe first that $A(E(x))$ and
$A(F(x))$ are the limits of the  (same) subsequences 
 before $E_k(f(x_k))$ and  $F_k(f(x_k))$. 
Then for any $m\in \ZZ$ we get
$$
2^m\,\frac{\|A^{ml}(u)\|}{\|u\|} \leq \frac{\|A^{ml}(v)\|}{\|v\|}
$$
for every $u\in E(x)$ and $v\in F(x)$.  
Now a standard dynamical argument asserts that the 
spaces $E(x)$ and $F(x)$ verifying this inequality are uniquely determined by their dimensions.

To complete the proof, observe that the unicity of the dominated splitting above gives the continuity.
\enddemo

\proclaim{{C}orollary}
\label{c.limsup}
 Let $f$ be a diffeomorphism defined on a compact manifold $M$ and
 $\Lambda$ an $f$\/{\rm -}\/invariant set.
 Assume that there are $l\in\NN${\rm ,}
$i\in 1,\ldots,\break \dim (M) -1${\rm ,} and
 a sequence of diffeomorphisms
$f_n$ converging to $f$ in  the $C^1$-topology
such that
\begin{itemize}
\item
 every $f_n$ has a periodic orbit $x_n$ such that $H(x_n,f_n)$ admits an
 $l$\/{\rm -}\/dominated splitting of dimension $i${\rm ,}
\item
the set  $\La$ is included in the topological upper limit set 
of the $H(x_n,f_n)${\rm ,} i.e.\ 
\vglue-8pt 
\hfil ${\displaystyle
\limsup_{n\to \infty} (H(x_n,f_n)) = \bigcap_{n=1}^{\infty} \mbox{\em closure}( \bigcup_{i>n} H(x_n,f_n)).
}$\hfill
\end{itemize}
Then
$\La$ admits an $l$-dominated splitting of dimension $i$.
\endproclaim

\demo{Proof} 
Consider the topological set 
$$
I=\{0\}\cup\{\frac1n , n\in\NN\setminus\{0\}\}.
$$
 In $M\times I$ we consider the union 
$$
(\La\times\{0\})\cup \bigcup_1^{+\infty}
(H(x_n,f_n)\times\{\frac1n\}).
$$
 The differentials of $f$ and $f_n$ define in a natural way a linear
 system on this set, which is continuous because the $f_n$ converge to
 $f$ in the $C^1$-topology. Moreover, $\bigcup_1^{+\infty }(H(x_n,f_n)\times\{\frac1n\}$ is a dense
subset (because $\La$ is contained in the topological upper limit set of
the $H(x_n,f_n)$) and the system over $\bigcup_1^{+\infty
}(H(x_n,f_n)\times\{\frac1n\}$ admits an $l$-dominated splitting. To finish the proof it is now enough to apply Lemma \ref{l.dense}. 
\enddemo

\demo{{\rm 1.4.} Periodic linear systems with transitions}
Saddles $P$ and $Q$ of the same index which are linked by transverse intersections of their invariant manifolds 
(i.e.\ they are {\em homoclinically related\/}) belong
to the same transitive hyperbolic set. So they are accumulated
by other periodic orbits which spend
an arbitrarily long time close to $P$, thereafter close to $Q$, and so on. 
In fact, the existence of  Markov partitions shows
that for any fixed  finite sequence of times
there is  a periodic  orbit expending
alternately the times of the sequence close to $P$ and $Q$,
respectively. Moreover, the transition time (between  a neighbourhood of $P$ and a neighbourhood of $Q$) can be chosen to be 
bounded. This property will allow us to scatter in the whole homoclinic class of $P$ 
some properties of the  periodic
points $Q$ of this  class.  

We aim  in this section  to translate this property into
 the language of  linear systems, introducing the
concept of {\em linear
system with  transitions.\/}
Then we shall deduce some direct consequences of the existence of such transitions.
Let us go into the details of our constructions.
We begin by giving some definitions.
 
Given a set $\cA$,
a {\em word with letters\/} in $\cA$ is a finite sequence of elements of $\cA$, its {\em length\/}
is the number of letters composing it. 
The set of 
words admits a natural   semi-group structure:
The product of the word $[a]=(a_1,\dots,a_n)$ by $[b]=(b_1,\dots,b_k)$ is $[a][b]=(a_1,\dots,a_n,b_1,\dots,b_k)$. 
We say that a word $[a]$ {\em is not a power\/}
if  
$[a]\ne[b]^k$ for every word $[b]$ and $k>1$.

In this section $(\Si,f,\cE,A)$ is
a periodic linear system of dimension $N$: Recall that every $x\in\Si$ is periodic for 
$f$, $p(x)$ denotes its period, and
 $M_A(x)$ denotes the product 
$A^{(p(x))}(x)$ of $A$ along the orbit of $x$.

If $(\Si,f,A)$ is a periodic system of matrices (in ${\rm GL}(N,\RR)$), then for any $x\in\Si$ we write
$[M]_A(x)=(A(f^{p(x)-1}(x)),\dots,A(x))$; which
is a word with letters in ${\rm GL}(N,\RR)$. 
Hence the matrix $M_A(x)$ is the product of the letters of the word $[M]_A(x)$.  

\numbereddemo{Definition} \label{d.epsilontransition}
Given $\varepsilon>0$, a periodic linear system 
$(\Si,f,\cE,A)$ admits {\it $\varepsilon$-transitions} if for every finite family of points
$x_1,\dots,x_n=x_1\in\Si$ there is an orthonormal system of  coordinates  of the linear bundle $\cE$ 
(so that $(\Si,f,\cE,A)$ can now be considered as a system of matrices $(\Si,f,A)$), and for any $(i,j)\in\{1,\dots,n\}^2$
there exist $k(i,j)\in \NN$ and a finite word  $[t^{i,j}]=(t_1^{i,j},\dots ,t_{k(i,j)}^{i,j})$ 
of matrices in ${\rm GL}(N,\RR)$, satisfying the following properties:
\vglue4pt
 1.  For every $m\in\NN$,  $\iota = (i_1,\dots,i_m)\in\{1,\dots,n\}^m$, and $a=(\alpha_1,\dots,\alpha_m)\break \in\NN^m$ consider
the word  
$$
[W(\iota,a)]= [t^{i_1,i_m}][M_A(x_{i_m})]^{\alpha_m}
[t^{i_m,i_{m-1}}][M_A(x_{i_{m-1}})]^{\alpha_{m-1}}\cdots [t^{i_2,i_1}][M_A(x_{i_1})]^{\alpha_1},
$$
where the word $w(\iota,a)=\left((x_{i_1},\alpha_1),
\dots,(x_{i_m},\alpha_m)\right)$ with letters in $M\times\NN$ {\it is not a power}.
Then there is $x(\iota,a)\in\Si$ such that
\begin{itemize}
\item[$\bullet$] 
The length of $[W(\iota,a)]$ is the period $p(x(\iota,a))$ of $x(\iota,a)$.
\item [$\bullet$]
The word $[M]_A(x(\iota,a))$ is $\varepsilon$\/{\rm -}\/close to $[W(\iota,a)]$
and there is an $\varepsilon$\/{\rm -}\/pertur\-bation $\tilde A$ of $A$ such that the word 
 $[M]_{\tilde A}(x(\iota,a))$ is $[W(\iota,a)]$.
\end{itemize}

2. One can choose $x(\iota,a)$ such that the distance between the orbit of $x(\iota,a)$ and any point $x_{i_k}$
is bounded by some function of $\alpha_k$ which
 tends to zero as $\alpha_k$ goes to infinity. 
\enddemo

Given  $\iota$ and $a$ as above, 
the word $[t^{i,j}]$ is an
{\em $\varepsilon$-transition from $x_j$ to $x_i$.\/}
We call {\em $\varepsilon$\/{\rm -}\/transition matrices\/} the matrices $T_{i,j}$ which are the product of the letters composing $[t^{i,j}]$.

\numbereddemo{{R}emark} \label{r.semigroup}
Consider points  $x_1,\dots,x_{n-1},x_n=x_1\in \Si$ and
$\varepsilon$-transitions $[t^{i,j}]$ 
from $x_j$ to $x_i$. Then 
\begin{itemize} 
\item[1.]   
for every positive 
$\alpha\geq 0$ and $\beta\geq 0$ the word $([M]_A(x_i))^{\alpha}\,[t^{i,j}]\,([M]_A(x_j))^{\beta}$
is also an $\varepsilon$-transition from $x_j$ to $x_i$,
\item[2.] 
for any $i,\,j,$ and $k$ the word $[t^{i,j}][t^{j,k}]$ is an
$\varepsilon$-transition from $x_k$ to $x_i$.
\item[3.]
As a consequence of the two items above, the words $W(\iota,\alpha)$ in
Definition~\ref{d.epsilontransition} are $\varepsilon$-transitions from
$x_{i_1}$ to itself. Moreover,  the set of such $\varepsilon$-transitions forms a semigroup.
\end{itemize}

\enddemo

\numbereddemo{Definition} \label{d.transitions}
We say that a periodic linear  system {\em admits transitions\/} 
if for any $\varepsilon>0$ it admits $\varepsilon$-transitions.
\enddemo

The following  lemma justifies the introduction of the notion of transition for  studying  homoclinic classes:

\proclaim{Lemma} \label{l.Df}
Let $P$ be a hyperbolic saddle of index $k$ {\rm (}\/dimension of its stable manifold\/{\rm ).} 
The derivative $f_*$ induces a  continuous periodic
linear system
with transitions on the set $\Si$ of hyperbolic saddles in $H(P,f)$ of index $k$ and homoclinically related
to $P$.
\endproclaim

\demo{Proof} 
Fix any $\varepsilon>0$ and a finite family $x_1,\dots,x_n$ in $\Si$. As the $x_i$ are homoclinically related to
$P$, there is a compact transitive hyperbolic subset $K$
of $H(P,f)$ containing all the $x_i$. 
So this set $K$ can be covered by a Markov partition with arbitrarily small rectangles. We can now choose orthonormal systems of coordinates in $T_x(M)$, $x\in K$, such that the orthonormal
bases depend continuously on $x$ when the points are in the same rectangle. 

Let $(K,f,A)$ be the system of matrices defined
on $K$ by writing $f_*$ in this system of coordinates. Now,
using the continuity of $f_*$, and by subdividing if necessary the rectangles of the Markov partition, we can assume that,
for any $x$ and $y$ in the same rectangle, 
 $$
\|A(x)-A(y)\|<\varepsilon
\quad \mbox{and}\quad
\|A^{-1}(x)-A^{-1}(y)\|<\varepsilon.
$$ 
The transitions from $x_i$ to $x_j$ are now obtained  by consideration of the derivative of $f$ along 
any orbit in $K$ going from the rectangle containing $x_i$ to the rectangle containing $x_j$.   
\enddemo

The next lemma shows how a property at one point of a system with transitions can   scatter to a dense subset: 

 \proclaimtitle{Scattering Property}
\proclaim{Lemma}
\label{l.linearsources}
Let $(\Si,f,\cE,A)$ be a 
periodic linear system with transitions. Fix $\varepsilon> \varepsilon_0 >0$ and assume that there exist
an\break $\varepsilon_0$-perturbation $\tilde A$ of $A$
and $x\in\Si$ such that $M_{\tilde A}(x)$ is
either a dilation {\rm (}\/i.e.\ all its eigenvalues have
modulus bigger than $1${\rm )} or a contraction {\rm (}\/i.e.\ all its eigenvalues have modulus less than $1${\rm ).}

Then there  are a dense $f$\/{\rm -}\/invariant subset $\tilde \Si$ of $\Si$ and
an $\varepsilon$\/{\rm -}\/perturbation $\hat A$ of $A$
such that for any $y\in\tilde\Si$ the linear map $M_{\hat A}(y)$
is either a dilation or a contraction {\rm (}\/according to the
choice before\/{\rm ).}
\endproclaim

\demo{Proof} Write $\varepsilon_1=\varepsilon-\varepsilon_0$, take
some point $z$ in $\Si$, and 
consider two $\varepsilon_1$-transitions $T_{x,z}$ (from $z$ to $x$) and $T_{z,x}$ (from $x$ to $z$).
For a fixed $\delta >0$,
by definition of transitions, 
there is $n(z,\delta)$ such that for any $n>0$ there are $y_n\in \Si$,
with $d(y_n,z)<\delta$, and an $\varepsilon_1$-deformation $A'$ of $A$ along the orbit of $y_n$ such that 
$$
M_{A'}(y_n)=T_{z,x}\circ M_A(x)^n\circ T_{x,z}\circ M(z)^{n(z,\delta)}.
$$ 
Define $\hat M_n$ by
$$
\hat M_n =T_{z,x}\circ M_{\tilde A}(x)^n\circ T_{x,z}\circ M(z)^{n(z,\delta)}.
$$
We can now choose $n$ big enough so that
$\hat M_n$   is either
a dilation or a contraction (according to $M_{\tilde A}(x)$). 
Thus by an  $\varepsilon_1$-perturbation $\hat A$ of $A'$
along the orbit of $y_n$ we can get $ M_{\hat A}(y) = \hat M_n$. 

Since we are not requiring the continuity of $\hat A$, we can build it as above, that is,  orbit by orbit 
considering points in a dense subset. 
This ends the proof of the lemma. 
\enddemo

\section{Quantitative results: Proofs of the theorems}
\label{s.qr}

In this section we state, in terms of linear systems (Proposition~\ref{p.main}) and in terms of diffeomorphisms (Proposition~\ref{p.quantitative}), quantitative results on the  existence of dominated splittings (giving the strength of the dominance). 

Proposition~\ref{p.main} gives a dichotomy between the existence of a dominated splitting for a linear 
system and the existence of perturbations
of the system with homotheties. This proposition is divided into two main steps: Proposition~\ref{p.rank}, asserting that the lack of dominance allows
us  to create complex eigenvalues, and Proposition~\ref{p.homothety}, which says that sufficiently many complex eigenvalues allow
us to get homotheties.   These propositions will be proved in the next two sections. 

In this section we deduce from Proposition~\ref{p.main} most of the results announced in the introduction.

\vglue12pt  2.1. {\it Reduction of the study of the dynamics
to a problem on linear systems}.
 
\proclaim{Proposition} \label{p.main}  For any $K>0${\rm ,} $N>0${\rm ,} and $\varepsilon>0$ there is $l>0$ such that
any continuous periodic $N$\/{\rm -}\/dimensional linear system $(\Si,f,\cE,A)$
bounded by $K$ {\rm (}\/i.e.\ $\|A\|<K${\rm )} and having  transitions satisfies
the following\/{\rm :}
\begin{itemize}
\item either $A$ admits an $l$\/{\rm -}\/dominated splitting{\rm ,} 
\item or there are an $\varepsilon$\/{\rm -}\/perturbation $\tilde A$ of $A$
and a point $x\in\Si$ such that $M_{\tilde A} (x)$ is an homothety.
\end{itemize}

\endproclaim

The proof of Proposition \ref{p.main} is divided in two main steps:
In the first one, we  show that, if $(\Si,f,\cE,A)$ is a  linear system with transitions
such that no dense subsystem of it 
admits an $l$-dominated splitting, then we can perturb $A$ to get
{\em a lot of complex eigenvalues.\/}
In the second step, we see that, if we can obtain sufficiently many complex eigenvalues, then 
we can perturb the system to get a homothety (which will be either a contraction or a dilation). 
Let us state precisely these two steps. We begin
with some definitions.

\numbereddemo{Definition} \label{d.rank}
Let $M\in {\rm GL}(N,\RR)$ be a linear isomorphism of $\RR^N$ such that $M$ has some complex  eigenvalue $\lambda$, 
i.e.\ $\lambda\in\CC\setminus\RR$. We  say
that $\lambda$ has {\em rank\/} $(i,i+1)$ if there is an $M$-invariant splitting of $\RR^N$, $F\oplus G\oplus H$, such that: 
\begin{itemize}
\item Every eigenvalue  $\sigma$ of $M|_F$ (resp.\ $M|_H$) has modulus  $|\sigma|<|\lambda|$ (resp.\ $|\sigma|>|\lambda|$),
\item $\dim (F) = i-1$ and $\dim (H) = N-i-1$,
\item the plane $G$ is the eigenspace of $\lambda$.
\end{itemize}
\enddemo

\numbereddemo{Definition} A periodic linear  system $(\Si,f,\cE,A)$ {\em has a complex eigenvalue of rank $(i,i+1)$\/} if there is 
$x\in \Si$ such that the matrix $M_A(x)$ has a complex (nonreal) eigenvalue of rank $(i,i+1)$. 
\label{d.complexe}
\enddemo

Proposition \ref{p.main} is a direct consequence of 
Propositions \ref{p.rank} and \ref{p.homothety} below:
 
\proclaim{Proposition} \label{p.rank}
For every $\varepsilon>0${\rm ,} $N\in\NN$, and $K>0$ there is $l\in \NN$ satisfying the following
 property\/{\rm :} 
\vglue9pt
Let $(\Si,f,\cE,A)$ be a continuous periodic $N$\/{\rm -}\/dimensional linear
system with transitions such that its norm $\|A\|$ is bounded by $K$.  Assume that  there exists\break $i\in\{1,\dots,N-1\}$ such
that every
$\varepsilon$\/{\rm -}\/perturbation $\tilde A$ of $A$ has no complex eigenvalues of rank $(i,i+1)${\rm .}  Then $(\Si,f,\cE,A)$
admits an 
$l$\/{\rm -}\/dominated splitting $F\oplus G${\rm ,}  $F\prec_l G${\rm ,} with $\dim (F)= i$.
\endproclaim 

\proclaim{Proposition} \label{p.homothety}
Let $(\Si,f,\cE,A)$  be a periodic  linear system with transitions. Given $\varepsilon>\varepsilon_0>0$
assume that{\rm ,} for any $i\in\{1,\dots,N-1\}${\rm ,}
there is an $\varepsilon_0$\/{\rm -}\/perturbation of $A$ having a complex
eigenvalue of rank $(i,i+1)$. Then there are an $\varepsilon$\/{\rm -}\/perturbation $\tilde A$
 of $A$ and $x\in \Si$ such that $M_{\tilde A}(x)$ is a homothety with ratio of modulus different from $1$. 
\endproclaim

The key of the proof of  Proposition \ref{p.rank} is a $2$-dimensional argument of Ma\~n\'e
that  we present in Section 3. The proof in higher dimensions consists of an inductive argument which allows
us to reduce the dimension of the linear space by 
considering some quotients (roughly speaking, considering
projections). Using this inductive procedure we
finally arrive at a two-dimensional space.
The lemmas in Section~4.1 allow us to make these successive reductions of dimension.
The proofs of Propositions~\ref{p.rank} and \ref{p.homothety}
are in Section~\ref{s.proofprop}.

Now using Proposition~\ref{p.main} we prove most of the results announced in the introduction.

\demo{{\rm 2.2.} Proofs of the theorems}
Let us first explain why 
Proposition \ref{p.main} implies Theorem \ref{t.main}.
Actually, this proposition implies the following 
 quantitative version of Theorem~\ref{t.main},
 which is our main (but a little bit technical) result:  
\enddemo

\proclaim{Proposition} \label{p.quantitative}
For every $K>0${\rm ,} $N>0${\rm ,} and $\varepsilon>0$ there is 
$l(\varepsilon, K,N)\in \NN$ such that 
for any diffeomorphism $f$ defined on a riemannian $N$\/{\rm -}\/dimensional manifold $M$ 
such that the derivatives $f_*$ and $f_*^{-1}$ are bounded by $K${\rm ,} 
and
any saddle $P$ of $f$ with a
nontrivial homoclinic class $H(P,f)${\rm ,}
the following holds\/{\rm :}\/  
\begin{itemize}
\item 
Either the homoclinic class $H(P,f)$ admits an $l(\varepsilon, K,N)$\/{\rm -}\/dominated splitting{\rm ,} 
\item 
or for every  neighbourhood $U$ of
 $H(P,f)$  and  $k\in \NN$ there is $g$ $\varepsilon$\/{\rm -}\/$C^1$\/{\rm -}\/close to $f$ having $k$ sources or sinks whose orbits are 
contained in $U$. 
\end{itemize}

\endproclaim

2.2.1. {\it Proofs of Theorem} \ref{t.main}  {\it and Proposition} \ref{p.quantitative}.
As Proposition \ref{p.quantitative} implies 
 Theorem \ref{t.main} directly, it remains to see that 
Proposition \ref{p.quantitative} follows from
Proposition \ref{p.main}.

For that, consider
a diffeomorphism  $f$, such that $\|f_*\|$ and $\|f_*^{-1}\|$
are  bounded by $K$, and a periodic saddle $P$ of $f$
with a nontrivial homoclinic class. 
 Let 
$$
\Si = H(P,f),\quad
\cE= TM|_{\Si},
\quad
\mbox{and}
\quad 
A = (f_*)|_{\Si}.
$$
Then $(\Si,f,\cE,A)$ is a continuous linear system.
Denote by $\Si'\subset \Si$ the set of saddles 
homoclinically related to  $P$ (in particular,
having the same index as $P$).  
Observe that $\Si'$ is a dense $f$-invariant subset of $\Si$. Moreover, by Lemma~\ref{l.Df}, 
 the subsystem induced by $A$
over $\Si'$ admits transitions.

If $A$ admits an $l-$dominated splitting over $\Si'$ then,
 by Lemma \ref{l.dense},
such a  splitting can be extended
to an $l$-dominated splitting on
the whole $\Si=H(P,f)$, and we are done.  

Now take the constant $l>0$ given by Proposition~\ref{p.main} corresponding to $K$, $N=\dim (M)$, and $\varepsilon/2$. 
If $A$ does not admit an $l$-dominated splitting over $\Si'$, then  Proposition~\ref{p.main} says that there is an
 $\varepsilon/2$-perturbation
$\tilde A$ of $A$ and a point $x\in \Si^\prime$ 
such that $ M_{\tilde A}(x)$ is a homothety.
We can suppose that (up to an arbitrarily small perturbation)
this homothety is either a dilation or a contraction.
Assume, for instance, the first possibility.

As the system admits transitions,  by 
Lemma \ref{l.linearsources},
there is a dense subset of $\Si'$ of points $y$ admitting
$\varepsilon$-deformations $\hat A$ along their orbits such that the corresponding linear map $M_{\hat A}(y)$
is a dilation. 
Choose now an arbitrarily large (but finite) number of 
such points $y$, and  denote by $E$ this set of periodic orbits. 

The proofs of Proposition  \ref{p.quantitative}
(thus of Theorem \ref{t.main}) follows now immediately from   Franks' lemma.
\hfill\qed
\vglue12pt
 
We now prove the  corollaries of Theorem~\ref{t.main} in the introduction.

\pagebreak 2.2.2. {\it Proof of Corollary} \ref{c.locdense}.
Let $P$ be a hyperbolic saddle of a diffeomorphism $f$ and $\cU$ a neighbourhood of $f$ where $P$ admits a continuation $P_g$ for every
$g\in \cU$.
 
 Denote by $\cD\cS$ the set of diffeomorphisms $g\in\cU$ for which $H(P_g,g)$ admits a dominated splitting. 
If the closure  of the interior
of $\cD\cS$, ${\rm cl}({\rm int}({\cD\cS}))$, is a neighbourhood of $f$ then the first possibility in the corollary holds and we are done.
Otherwise, for any $\varepsilon>0$ there is $g_1$, $\varepsilon/2$-close
to $f$ in the $C^1$-topology, in the complement of ${\rm cl}({\rm int}(\cD\cS))$.
Thus  $g_1$ has an open neighbourhood $\cU_1$ 
such that for any $g\in \cU_1$
 there is $h$ arbitrarily close to $g$ such that $H(P_h,h)$ does
not admit any dominated splitting.

Given a set $E\subset M$ and $\delta>0$, 
let $V(E,\delta)$ be the set of points of $M$
at distance strictly less than $\delta$ from $E$.
We now construct inductively sequences $\varepsilon_i>0$
and  $g_i\in\cU_1$ satisfying
the following properties:
\begin{itemize}
\item[1.] $H(P_{g_i},g_i)$ has no dominated splitting,
\item[2.] $g_{i+1}$ is $\varepsilon_i/2$-close to $g_i$ in the $C^1$-topology,
\item[3.] there is  a finite set $\cS_{i+1}$ of periodic sinks or sources of $g_{i+1}$ such that  $H(P_{g_i},g_i)\subset V(\cS_{i+1},
\varepsilon_i/2)$,
\item[4.] $\varepsilon_{i+1}<\varepsilon_i/2$,
\item[5.] for all $g$ $\varepsilon_{i+1}$-close to $g_{i+1}$
the set of sinks or sources $\cS_{i+1}$ has 
a continuation $\cS_{i+1}(g)$ such that $H(P_{g_i},g_i)\subset V(\cS_{i+1}(g), \varepsilon_i)$.
\end{itemize}

Let us first end the proof of the corollary using the sequences $(\varepsilon_i)$ and $(g_i)$ above.
The sequence $(g_{i})$ is a 
Cauchy sequence in Diff$^1(M)$,  so
it converges to some $C^1$-diffeomorphism $h$. Moreover, from
$\varepsilon_{i+1}<\varepsilon_i/2$
 and the $\varepsilon_i/2$-proximity of $g_{i+1}$
 to $g_i$, we get
that $h$ is $\varepsilon_i$-close to $g_i$ for all $i$. 
Therefore, by item (5), the set of sources or sinks $\cS_{i+1}(h)$ is well defined
and $H(P_{g_i},g_i)\subset V(\cS_{i+1}(h), \varepsilon_i)$ for every $i$. 

Consider now the set 
$\cS(h)=\bigcup_1^{\infty}(\cS_i(h))$ 
consisting of sinks or  sources. By construction, the closure of $\cS(h)$
contains
the topological upper limit set of the $H(P_{g_i},g_i)$; that is,
$$
\mbox{closure}\,(\cS(h)) \supset
\limsup_{i\to\infty}H(P_{g_i},g_i)= \bigcap_{i=1}^{\infty}\,\mbox{closure}\,(\bigcup_{j>i}H(P_{g_j},g_j)).
$$
Finally, by definition of homoclinic class and since the transverse intersections vary continuously, this upper limit set contains $H(P_h,h)$, 
so that $H(P_h,h)$ is contained in the closure of the set of sinks or sources of $h$. Thus  $h$ is the diffeomorphism in the statement
of the corollary.

To end the proof of the corollary it remains
to build the sequences 
$(\varepsilon_i)$ and $(g_i)$ above. 
We proceed inductively, assuming that $\varepsilon_j$ and $g_j$
are defined  for  every $j\leq i$.
 Consider some finite set $\Si_i\subset H(P_{g_i},g_i)$ of saddles such that $H(P_{g_i},g_i)\subset V(\Si_i,\varepsilon_i)$. 
By item (1), applying Theorem \ref{t.main} finitely many times, we can create
a sink or a source close to each point of $\Si_i$,
obtaining a diffeomorphism $g_{i+1}$ which is $\varepsilon_i/2$-close to $g_i$ 
and has
a set of sinks or sources $\cS_{i+1}$ containing $H(P_i,g_i)$ in its
$\varepsilon_i/2$-neighbourhood $V(\cS_{i+1},\varepsilon_i/2)$.  Thus
 $g_{i+1}$ satisfies items (2) and (3). 
Having in mind the definition of $\cU_1$,
we can suppose (after a new perturbation if necessary) that $H(P_{g_{i+1}},g_{i+1})$ has no dominated splitting, i.e.\ $g_{i+1}$ satisfies item (1). 
Then, using the continuous variation of the finite set $\cS_{i+1}(g)$ in a small neighbourhood of $g_{i+1}$, we can  choose
$\varepsilon_{i+1}<\varepsilon_i/2$ (item (4)) verifying item (5) above.
This ends the proof of the corollary. 
\hfill\qed

\vglue12pt 2.2.3. {\it Proof of Corollary} \ref{c.residual}.
Recall first that there are dense open subsets $\cO_n$ of Diff$^1(M)$ of diffeomorphisms $f$ for which the set  $\cP(f,n)$ of 
periodic points of
period less than $n$ is finite and hyperbolic.
Note also that the cardinal of $\cP(f,n)$ is locally constant in $\cO_n$ and that the set
$\cP(f,n)$ depends continuously on $f$. 

 Denote by $\Si(n,f)\subset \cP(f,n)$ the set of saddles 
with nontrivial homoclinic class.  Then there 
is a dense open subset $\tilde \cO_n\subset\cO_n$ of diffeomorphisms $f$ 
such that
$\Si(n,f)$ has locally constant cardinal and depends continuously on $f$.

\nonumproclaim{Claim} There is a residual subset $\cR_n\subset \tilde\cO_n$ of diffeomorphisms
 $f$ such that for any $P\in\Si(n,f)$
either $H(P,f)$ admits a dominated splitting or $P$
belongs to the closure of the set of sinks or sources.
\endproclaim

\demo{{P}roof of the claim} Consider any open subset $\cO\subset \tilde\cO_n$
where the periodic
points in  $\Si(n,f)$ are continuous functions of $f$.
So let us denote by $\Si_n$ the (finite) set of these functions:
 Given $P\in\Si_n$ and $f\in\cO$ we  denote by
$P_f$ the corresponding periodic point of $f$.

For a fixed $P\in\Si_n$   let $\cD\cS(P)$ be the set of $f\in\cO$ 
such that $H(P_f,f)$ admits a dominated splitting.
Let 
$\cU(P)$ be the complement in $\cO$ of the closure of the interior of $\cD\cS(P)$.
Let $\cU(P,i)$  be the set of $f\in\cU(P)$ for which there is a sink or a source $Q_f$ (of any period) with 
$d(P_f,Q_f)<1/i$: This set is open and, by Theorem \ref{t.main}, dense in $\cU(P)$. 
Therefore the intersection 
$$
\cR_n(P)=\bigcap_{i=1}^\infty (\cU(P,i)\cup \mbox{int}(\cD\cS(P)))
$$
is a residual subset of $\cO$.
We write
$$
\cR_n(\cO)=\bigcap_{P\in\Si_n} \cR_n(P),
$$
by construction, noting that  the set $\cR_n(\cO)$ is a residual subset of $\cO$ consisting of diffeomorphisms $f$ satisfying the
conclusion in the  claim. 
Thus to end the proof of the 
claim it suffices to 
consider the set $\cR_n$ obtained as the union
(over all the open sets $\cO\subset\tilde\cO_n$) of  the $\cR_n(\cO)$.\hfill\qed
\enddemo

We are now ready to end the proof of the corollary.
Consider $\cR=\bigcap_{n\in\NN} \cR_n$. By the claim,
 the set $\cR$ is a residual subset of Diff$^1(M)$ of diffeomorphisms $f$ such that  for any 
saddle $P$ of $f$ there are two possibilities,
either $H(P,f)$ has a dominated splitting, or $P$ is in
the closure of the set $\cS(f)$ of sinks or sources of $f$ 
(remark that if the homoclinic class of $P$ is trivial then it admits a dominated splitting because $P$ is hyperbolic).

In the first case we are done. In the second one, we need to check
that the whole homoclinic class of $P$ is contained in
the closure of the sinks or sources.
So assume that $H(P,f)$ does not admit any dominated splitting. 
Observe that for every saddle $Q$ homoclinically related to $P$ one has $H(Q,f)=H(P,f)$, thus $H(Q,f)$ has no dominated splitting.
As $f\in\cR$, we have just seen that $Q$ is in
the closure of $\cS(f)$. 
Since the set of saddles homoclinically related to $P$ is dense in $H(P,f)$ we have that $H(P,f)$ itself is contained in the closure  of
$\cS(f)$. So $\cR$ is the residual set announced in Corollary \ref{c.residual} and the proof is complete.  
\hfill\qed

\vglue6pt 2.2.4. {\it Proof of Theorem} \ref{t.robust}.
Let us first prove this theorem
in the case of robustly transitive diffeomorphisms. There are two
reasons for that. First,
the proof of this case is simpler than the proof in
 the case of transitive
sets (i.e.\ the general case).
Second, proceeding in this way can emphasize
the additional difficulties and subtleties of the 
proof for transitive sets.

\demo{{P}roof of Theorem {\rm \ref{t.robust}} for robustly transitive diffeomorphisms} 
Consider a $C^1$-ro\-bus\-tly transitive diffeomorphism $f$
and an open neighbourhood  $\cU$ of $f$ such that any $g\in \cU$ is transitive. Reducing the size of $\cU$ if necessary, 
we can assume that there are $K>0$ and $\varepsilon>0$ such that every
$\varepsilon$-perturbation $h$ of any $g\in\cU$ is transitive and  the differentials $h_*$ and $h_*^{-1}$ 
 are bounded by $K$. 

Recall that by Pugh's closing lemma \cite{P}  there is 
 a residual subset of Diff$^1(M)$ 
of diffeomorphisms whose  nonwandering set is the closure of the hyperbolic periodic points. So there is a residual subset 
$\cR_0$ of $\cU$ of diffeomorphisms $g$ having a dense set of hyperbolic saddles (note that due to
the transitivity the diffeomorphisms in $\cU$ have
neither sinks nor sources).

Moreover, \cite[Th.~B]{BD2} says that 
there is a residual set $\cR_1$ of Diff$^1(M)$
of diffeomorphisms $f$  
such that two periodic points 
of $f$ belong to the same transitive set
if and only if their homoclinic classes are equal. 
Thus for any $g\in\cR = \cR_0\cap\cR_1$ and every
periodic point $P_g$ of $g$  
the homoclinic class $H(P_g,g)$ is the whole manifold $M$. 

By the robust transitivity of the $g$, given by the choice of $\varepsilon$,
it is not possible 
to create a sink or a source by an $\varepsilon$-perturbation of
any $g\in\cR$.
So Proposition~\ref{p.quantitative} gives 
$l$ such that every $g\in\cR$ admits an $l$-dominated splitting on $M=H(P_g,g)$.
Finally, choosing a sequence
$g_n\in\cR$ converging to $f$, Corollary \ref{c.limsup} ensures that $f$
admits an $l$-dominated splitting, ending
the proof of the theorem for robustly transitive diffeomorphisms.
\hfill\qed
\enddemo

{\it Proof of Theorem} \ref{t.robust}, {\it general case.}
Let $\La_f = \bigcap_{-\infty}^{+\infty} f^i(\bar U)\subset U$ be a\break
$C^1$-robustly transitive set in some open set $U$. 
Assume that $\La_f$ is not reduced to a single hyperbolic orbit 
(in that case we have nothing to do).

The proof follows essentially along the arguments
in the robust transitive case above, 
but we need to pay special attention to the following fact:
In the transitive case  ($U=M$) all the
orbits we consider are automatically in $\Lambda_f=M$
(that is, a tautology),
but {\em a~priori\/} this does not happen when $U\ne M$. 
Let us go into the details of the proof
 of this case.

Let $\cU$ be a $C^1$-open neighbourhood of $f$ such that, for every  $g\in \cU$, the maximal invariant set 
$\La_g= \bigcap_{-\infty}^{+\infty} g^i(\bar U)$ is transitive. 
As above, using Pugh's closing lemma
and  \cite[Th.~B]{BD2},
we get
a residual subset $\cR_0$ of $\cU$ of diffeomorphisms
 $g$ such that the hyperbolic periodic points of $\La_g$
 are dense in $\La_g$ and  have the same  homoclinic classes. 
Thus, for every $g\in\cR_0$, the set $\La_g$
is included in the homoclinic class $H(P,g)$ of some periodic point $P$. However, that is the special difficulty of this
case; we do not know {\it a priori\/} if $H(P,g)$ is contained in $U$
(in the robust transitive case that is obvious: $U=M$~!)

To solve this problem
denote
by $H(P,g,U)$ the points of 
the closure of the transverse intersections of the invariant manifolds of $P$
whose orbits remain in~$U$. We call this set the {\em  homoclinic class of $P$ in $U$.\/} 
So $H(P,g,U)$ is a transitive compact subset of $\La_g$ and,
since $\La_f\subset U$, it is  far from the boundary of $U$. 

\proclaim{Lemma}
\label{l.ro}
There is a residual set $\cR_2$ of $\cU$ such that for any $g\in\cR_2$ there is a periodic point $P$ such that 
$H(P,g,U)=\La_g.$     
\endproclaim

{\it Proof}. The proof is identical to that in  \cite[Th.~B]{BD2}
by the following version of Hayashi's connecting  lemma. So we 
do not go into the details.

\nonumproclaim{Theorem {\rm (Hayashi's Connecting Lemma)}}
Let $M$ be a compact manifold{\rm ,} $U$ an open set of $M${\rm ,}
$V$ an open set relatively compact in $U${\rm ,} and $f$ a diffeomorphism defined on $M$. 

Assume that there are  periodic saddles 
$P$ and  $Q$ whose orbits are contained in $U${\rm ,} a
sequence of points $x_i$ converging to some point $x\in W_{\rm loc}^u(P,f)${\rm ,} and a sequence 
of positive integers $n_i$ such that $f^{n_i}(x_i)$ converges to
some  $y\in W_{\rm loc}^s(Q,f)$.
Suppose also that for any $i$ and every
 $m\in\{0,\dots,n_i\}$ one has $f^{m}(x_i)\in V$. 

Then there is $g$ arbitrarily $C^1$\/{\rm -}\/close to $f$ such that  $x\in W^u(P,g)\cap W^s(Q,g)$.
Moreover{\rm ,} the whole orbit of $x$ is contained in $U$ and $g^n(x)=y$
for some $n>0$. \phantom{morerain}\hfill\qed
  \endproclaim

To prove Theorem~2 we apply Proposition~\ref{p.quantitative}, so
we need to see that the sets
$\La_g= H(P_g,g,U)$ given by Lemma~\ref{l.ro} are not all reduced to the single  periodic orbit $P_g$.

\proclaim{Lemma}
\label{l.2pp}
There is $g$ arbitrarily $C^1$\/{\rm -}\/close to $f$ having at least two different hyperbolic periodic orbits
contained in $U$.
\endproclaim

By Lemma~\ref{l.ro}, this  implies that one can choose $g$ such that $\La_g=H(P_g,g,U)$ is not reduced to the periodic orbit $P_g$.
  
\demo{Proof}
We are assuming that $\La_f$ is not reduced to a 
single periodic orbit (this is just the trivial case).
Hence, if $f$ has at least two periodic orbits the result
 is immediate: After perturbation we can make
these orbits
hyperbolic ones. 

So it remains to consider the case where $f$ has no periodic
orbits.
We know that there are diffeomorphisms $g$ close to $f$ having 
periodic orbits. 
We argue by contradiction:
suppose that every  $g$ (with periodic points) close to $f$ has only
one periodic orbit $Q$. Then considering an isotopy from
$g$ to $f$ we get a bifurcation of this periodic orbit. 
After a
perturbation, we get a saddle-node, a flip, or a Hopf
bifurcation. In these three cases,  a new perturbation gives 
two periodic orbits: A saddle-node and a flip split into two
hyperbolic periodic points, and a Hopf point into a periodic
point and an invariant circle (in this case to get a new periodic point it is enough to
modify the rotation number of the restriction of the
map to the invariant circle). 
\enddemo

By definition of robust transitivity, there are no perturbations of $f$ having sinks or sources
whose orbits are
contained in $U$. Take a sequence $g_n\to f$ of diffeomorphisms such that $\La_{g_n} = H(P_{g_n},g_n,U)$ is nontrivial. Proposition \ref{p.quantitative} 
implies that there is $l\in\NN$ such that 
$\La_{g_n}$ admits an $l$-dominated splitting for any $n$ large enough.
Corollary \ref{c.limsup} now implies that these dominated splittings
induce an $l$-dominated splitting on 
$\limsup_{n\to \infty} (\La_{g_n})$,
 and so on $\La_f$.
Now the proof of the theorem is complete.
\hfill\qed

\vglue12pt 2.2.5. {\it Proof of Theorem} 3.
This theorem  is a direct consequence of Proposition \ref{p.quantitative}.
We argue as follows:  let $m\in \NN$, $K>0$, and  $\varepsilon>0$ such
that any diffeomorphism $g$ which is  $\varepsilon$-close to $f$ has less than $m$ sinks and sources, and $g_*$
 and $(g_*)^{-1}$ are both uniformly bounded by $K$. 
Let $l_0$ be the constant 
$l(K,\varepsilon/2, \dim (M))$ given by Proposition \ref{p.quantitative}. Then, 
for every $g$ $\varepsilon/2$-close to $f$ and any saddle $P$
of $g$ having a nontrivial homoclinic class $H(P,g)$,
we have that $H(P,g)$ has an $l_0$-dominated splitting.    
\hfill\qed

\vglue12pt
To prove  Theorem 4 we  need  new arguments of
a very
different nature:
 the notion of {\em finest dominated splitting\/} and the 
{\em ergodic closing lemma\/} of Ma\~n\'e;  so let us postpone its proof until  Section~\ref{s.volume} of this paper.

We also postpone until the end of the paper 
(Section~\ref{s.con}) the results
about conservative diffeomorphisms.

\section{Two-dimensional linear systems \label{s.2dim}}

In this section we give a version (Proposition~\ref{p.2dim})
of Proposition~\ref{p.main} for two-dimensional systems  (without requiring transitions) following an argument essentially due to
Ma\~n\'e in \cite{Ma3}.  In the next sections, using an argument of reduction of the dimension of the system (quotients and restrictions
of linear systems; see Section~4.1) we deduce from this two-dimensional result the general version of it (see
Section~\ref{s.proofprop}).

\proclaim{Proposition} \label{p.2dim}
Given any $K>0$ and $\varepsilon>0$ 
there is $l\in\NN$
such  that for every two\/{\rm -}\/dimensional linear system
$(\Si,f,\cE, A)${\rm ,} with norm  $\|A\|$ bounded by $K$ and 
such that the matrices $M_A(x)$ preserve
the orientation{\rm ,}  
\begin{itemize}
\item either $A$ admits an $l$\/{\rm -}\/dominated splitting{\rm ,}
\item or there are  an $\varepsilon$\/{\rm -}\/perturbation $\tilde A$
of $A$ and  $x\in\Si$ such that $ M_{\tilde A}(x)$
has a complex {\rm (}\/nonreal\/{\rm )} eigenvalue.
\end{itemize}

\endproclaim

The difference between Proposition~\ref{p.2dim} and Proposition~\ref{p.main} (in the case
of $2$-dimensional systems) is that here we get a complex eigenvalue instead of a homothety. In fact, if the system admits transitions then one can use this complex eigenvalue to get  homotheties (this will be done later in any dimension, see Proposition~\ref{p.homothety}).  

We begin the proof of Proposition \ref{p.2dim} by a very elementary lemma whose proof we omit: 
 
\proclaim{Lemma} \label{l.rotation}
For every  $\alpha>0$  and  every matrix $M\in {\rm GL}_+(2,\RR)$ having two different eigenspaces $E_1$
and $E_2$  whose angle is less than $\alpha,$  there is $s\in[-1,1]$ such that  $R_{s\,\alpha}\circ M$ has
a complex {\rm (}\/nonreal\/{\rm )} eigenvalue {\rm (}\/here $R_{t\,\alpha}$
 denotes the rotation of angle $t\,\alpha${\rm ).}
\endproclaim

In what follows, for notational simplicity, let us write
$I_\mu=
\left(\begin{array}{cc}
                                 1 & \mu\\
                                  0   & 1\\
                    \end{array}\right)$.

\proclaim{Lemma} \label{l.angle}
For every $\alpha>0$ and $\mu>0$ there is $c>1$ verifying the following property\/{\rm :} 
Consider the matrices
$$
B=\left(\begin{array}{cc}
                                 b_1 & 0\\
                                  0   & b_2\\
                    \end{array}\right)
\,
\mbox{and}
\,
C= \left(\begin{array}{cc}
                                 c_1 & 0\\
                                  0   & c_2 \\
                    \end{array}\right)
\,
\mbox{ such that }
\,
\frac{|b_1|}{|b_2|}>c
\,
\mbox{ and }
\,
\frac{|b_1\, c_1|}{|b_2\,c_2|}<1.
$$ 
Then the angle between the eigenvectors of the matrix
$D= B\circ I_\mu
\circ C$ is less than $\alpha$.
\endproclaim

\demo{Proof} Observe first that $(1,0)$ 
is an eigenvector of the matrix $D$. 
The heuristic idea of the proof is very simple:
Consider the vector $(1,\beta)$, for some 
{\em small\/} $\beta\leq 2/(c\,\mu)$ fixed. As $|b_1/b_2|$
and $|c_2/c_1|$ are large (i.e.\ greater than $c$) the vectors $B^{-1}(1,\beta)$ and $C(1,\beta)$ are almost vertical (angle with the vertical less than $\mu$). The role of the matrix 
$I_\mu$ now is to send the direction of $C(1,\beta)$ into the direction of $B^{-1}(1,\beta)$,
thus  $(1,\beta)$ is an eigenvector of $D$.  
 
The precise calculations are not more complicated: 
Let $(1, \beta)$, $\beta\ne 0$, be some eigenvector of $D$ 
not parallel  to $(1,0)$ (i.e.\ associated to
the eigenvalue
$b_2\, c_2$). 
Then $\beta$ satisfies
$$
|\beta| = \frac{|b_2 \, c_2 -b_1\, c_1|}
{|\mu \,  b_1 \, c_2|}
=\frac{|b_2/b_1 - c_1/c_2|}{\mu}<
\frac{2}{c\,\mu}.
$$
This completes the proof of the lemma.
\enddemo

Consider a  periodic system of matrices $(\Si,f,A)$ in ${\rm GL}_+(2,\RR)$ such that 
all the matrices of the system are diagonal. Thus
the canonical splitting $\RR^2 = \RR\oplus\RR$ is invariant. Given $x\in \Si$ denote by
$\sigma(x)$ and $\lambda(x)$
the eigenvalues of $M_A(x)$ associated with the vertical direction
($\{0\}\times\RR$) and the horizontal direction
($\{0\}\times\RR$), respectively. 
Up to a trivial change of coordinates, one can assume that for any $x\in \Si$, the eigenvalue $\sigma(x)$ of $M_A(x)$ 
 is bigger in modulus than the eigenvalue $\lambda(x)$.

\proclaim{Lemma} \label{l.domangle}
For any $\varepsilon>0${\rm ,} $\alpha>0${\rm ,}
 and $K>0$ there is $l\in \NN$ with the following property\/{\rm :}

Consider a periodic system $(\Si,f,A)$ of diagonal
matrices in ${\rm GL}_+(2,\RR)$ as above{\rm ,} bounded by $K$
such that $|\sigma(x)| \ge |\lambda(x)|$ for every
$x\in \Si$.

Suppose that the  splitting $\RR^2 = \RR\oplus\RR$
is not $l$\/{\rm -}\/dominated. Then there are an $\varepsilon$\/{\rm -}\/per\-tur\-ba\-tion 
$\tilde A$ of $A$ and $x\in\Si$ such that the angle between the eigenspaces of $M_{\tilde A}(x)$ is less than $\alpha$.
\endproclaim

\demo{Proof}
Let us write 
$$A(x) = \left(\begin{array}{cc}
                                 a(x) & 0\\
                                  0   & b(x)\\
                    \end{array}\right). 
$$
Observe that if the splitting
$\RR^2=\RR\oplus\RR$ is $l$-dominated  
then for every $x\in\Si$ one has
$$
 2\, |\prod_0^{l-1} a(f^i(x))|<
|\prod_0^{l-1} b(f^i(x))|.
$$
Recall that, by Lemma \ref{l.changes},
there is $\mu>0$  (depending on $\varepsilon$ and $K$) such that 
multiplying  matrices $A(x)$
by diagonal matrices $\mu$-close to the identity one gets $\varepsilon/3$-perturbations of $A$. 

Suppose first that 
there is $x$ in $\Si$ such that 
$$
|\sigma(x)|\leq (1+\mu)^{2\,p(x)}\,|\lambda(x)|, \quad \mbox{ where $p(x)$ is the period of $x$.}
$$
 Then  multiplying
the matrices $A(f^i(x))$ by some matrix of the form  
$$
 \left(\begin{array}{cc}
                                 1+\nu & 0\\
                                  0   & 1/(1+\nu)\\
      
              \end{array}\right), 
\quad \mbox{for some $\nu\in[0,\mu],$}
$$
 we get an $\varepsilon/3$-perturbation $A'$ of $A$ such that $M_{A'}(x)$ is an homothety.
So given any pair of  (different) directions of $\RR^2$, 
there is an arbitrarily small perturbation
$\tilde A$ of $A'$ such that $M_{\tilde A}(x)$ 
has two eigenvectors parallel to such directions.
This ends the proof of the lemma in this first case.

So we can now assume that  
$$
|\sigma(x)|>(1+\mu)^{2p(x)}\,|\lambda(x)| 
\quad\mbox{ for every $x\in\Si$.}
$$ 
Consider the constant $c$, given by Lemma \ref{l.angle}, associated to $\alpha$ and $\mu$,  and
$l$ such that 
$$
(1+\mu)^l> 2\,c.
$$
 We  show that, for any system
of matrices $(\Si,A,f)$ bounded by $K$, 
$l$ is the constant announced in the statement of
the lemma. 

Recall the observation in the beginning of the proof
of the lemma;  since
 the canonical splitting $\RR^2=\RR\oplus\RR$
is not $l$-dominated, there is $x\in\Si$ such that 
$$ 
2\,|\prod_0^{l-1} a(f^i(x))| \geq |\prod_0^{l-1} b(f^i(x))|.
$$
Assume first that $l<p(x)$. Given
 $y\in \Si$ let 
\begin{eqnarray*}
\tilde a(y)& =& (1+\mu)\,a(y)
\mbox{ if }
y= f^i(x),i\in\{0,\cdots,l-1\},
\\
 \tilde a(y)& =&a(y)
\mbox{ if }
i\in\{l, \dots, p(x)-1\}.
\end{eqnarray*}
 Consider now 
 the $\varepsilon/3$-perturbation $\tilde A$
of $A$ given by 
$$
\tilde A(y) = \left(\begin{array}{cc}
                                 \tilde a(y) & 0\\
                                  0   & b(y)\\
                    \end{array}\right) .
$$  
Let $B=\tilde A^{(l)}(x)$ (the product of the  matrices 
of the system $\tilde A$ along the orbit of $x$ from time $0$ to time $l-1$)
and  $C=\tilde A^{(p(x)-l)}(f^l(x))$. Then, we get 
$$
M_{\tilde A}(x)=C\circ B\quad\mbox{and}\quad M_{\tilde A}(f^l(x)) = B\circ C.
$$
Observe that $B$ and $C$
verify the hypotheses of Lemma \ref{l.angle}, so that the angle between the eigenvectors of 
the matrix 
$D= B\circ I_\mu \circ C$ is less than $\alpha$.

 Denote by $\bar A$ the $\varepsilon/3$-perturbation
of $\tilde A$ obtained  modifying only the matrix $\tilde A(f^{-1}(x))$ $=  A(f^{-1}(x))$ by
replacing  this matrix by 
$\bar A(f^{-1}(x))= 
I_\mu \circ A(f^{-1}(x))$.
Then the angle between the eigenvectors
of $M_{\bar A} (f^l(x))= D$ is less than 
$\alpha$.

To finish the proof of Lemma~3.4    the case $l\geq p(x)$ remains. Remark that $l$ cannot
be a multiple of $p(x)$,  hence $l= k \, p(x)+ l_0$,
for some $k \ge 1$ and  $1\le l_0< p(x)$.
By hypothesis,
 $$
 (1+\mu)^{2\,p(x)}\,
|\prod_0^{p(x)-1} a(f^i(x))|< |\prod_0^{p(x)-1} b(f^i(x))|.
$$
Since we are assuming that the splitting is
not $l$-dominated, we also have
$$
 2\,|\prod_0^{l-1} a(f^i(x))|\geq |\prod_0^{l-1} b(f^i(x))|.
$$
So we get that 
$$ 
|\prod_0^{l_0-1} a(f^i(x))|> \frac{(1+\mu)^{2\, k\, p(x)}}{2}
\,|\prod_0^{l_0-1} b(f^i(x))|.
$$
Finally, note that $2\,k\, p(x)>l$ and recall the choice of $l$ (i.e.\ $(1+\mu)^l>2\,c$); thus 
$$
\frac{(1+\mu)^{2\,k\,p(x)}}{2} >\frac{(1+\mu)^l}{2} >c.
$$
Now, as in the previous case, we can apply Lemma \ref{l.angle} to the matrices $B=A^{(l_0)}(x)$ and $C= A^{(k\,
p(x)-l_0)}(f^{l_0}(x))$.  This completes the proof of the lemma.
\phantom{coffee}\enddemo 

 3.1. {\it Proof of Proposition} \ref{p.2dim}.
The proof of Proposition \ref{p.2dim} follows almost directly from Lemmas \ref{l.changes}, \ref{l.rotation}, and
\ref{l.domangle}. Let us go into the details.  

For a fixed $\varepsilon>0$ and $K>0$, 
consider  a system of matrices $(\Si,f,A)$ in ${\rm GL}_+(2,\RR)$,
bounded by $K$, such that 
it is not possible to create a complex eigenvalue by an
$\varepsilon-$per\-tur\-ba\-tion of $A$. We prove that such a system is\break
$l$-dominated. This clearly implies the
proposition.

By Lemma \ref{l.changes},
there is $\alpha=\alpha(K,\varepsilon)>0$
such that the composition of the system with
a rotation of angle
less than $\alpha$ gives an $\varepsilon/2$-perturbation
of the system. 
So Lemma \ref{l.rotation}  ensures that 
for any $x\in\Si$ the angle between the eigenspaces of $M_A(x)$ is bigger than $\alpha$.
This means that
there are $K_0 = K_0(\alpha)$ and 
a family of matrices $P(x)$, $x\in \Si$,
with
$P(x)$ and $P^{-1}(x)$ bounded by $K_0$, 
such that  the eigenspaces  of $P(x)\circ M_A(x)\circ P(x)^{-1}$ are orthogonal. 

Let $B=P\circ A\circ P^{-1}$ be the system 
defined by $B(x)= P(f(x))\circ A(x)\circ P(x)^{-1}$.
By  Lemma \ref{l.changes},  there are $K_1 = K_1(K,K_0)$  and $\delta =
\delta(K,K_0,\varepsilon)$ such that $B$ and $B^{-1}$ are bounded by
$K_1$ and any $\delta$-perturbation of $B$ is obtained by conjugating by
$P$ some $\varepsilon$-perturbation of $A$.

By construction, all the eigenspaces of
the matrices $M_B(x)$ are orthogonal.
Thus, by an orthonormal change of coordinates,
we can assume that  $B$ left invariant the canonical splitting $\RR^2=\RR\oplus\RR$,
and that the eigenvalue of $M_B(x)$ corresponding to the vertical direction is bigger (in modulus) than the 
eigenvalue associated to the horizontal
direction. Finally, there is no possibility to create a complex eigenvalue
by a $\delta$-perturbation   of $B$.

Fix $K_2=K_2(K_1)$ such that any $\delta$-perturbation of $B$ is bounded by $K_2$. Consider now 
$\alpha_2=\alpha_2(K_2)$ such that the composition by a rotation of
angle at most $\alpha_2$ of  a linear system bounded by $K_2$ gives a
$\delta/3$-perturbation.
Then (following Lemma \ref{l.rotation}) given any $\delta/3$-perturbation $\tilde B$ of $B$ the angles
between the eigenspaces of $M_{\tilde B}(x), x\in\Si$, are bigger than $\alpha_2$. 
So, from Lemma \ref{l.domangle}, we get $l_0(\delta,K_1,\alpha_2)$ such
that $B$ is $l_0$-dominated.

Using the fact that $A= P^{-1}\circ B\circ P$, where $P$ is bounded by $K_0$, we
get $l(l_0,K_0)$ such that $A$ is $l$-dominated. 

Now, to conclude the proof it is enough
to remark that all the constants introduced 
in the proof are functions of
$K$ and $\varepsilon$.
 \hfill\qed

\vglue-6pt
\section{Invariant subbundles: Reduction of the dimension\\ and the finest dominated splitting}
\label{s.is}
\vglue-6pt

4.1. {\it Quotient of linear  systems
and restriction to subbundles}. The proof of Proposition~\ref{p.rank} uses successively
Proposition~\ref{p.2dim} on $2$-dimensional subbundles.  
Our construction also involves an argument (considering
quotient spaces) that allows us to reduce the dimension of
the ambient space.
We pay  special attention to the invariant subbundles of a linear system, and we will often need to compare the action of the system on them. This motivates the  introduction of the notions of {\em restriction\/} and  {\em quotient\/} of a linear system.

Let $(\Si,f,\cE,A)$ be a linear system and  
$F$  an invariant subbundle  of $\cE$
(with constant dimension). 
We denote by $A_F$ the {\em restriction of $A$ to $F$\/} and 
by $A/F$ the {\em quotient of $A$ along $F$\/}
 endowed with the metric of the orthogonal complement $F^\perp$ of $F$;
i.e.,  given a class $[v]$ we let 
$$
|[v]|= |v^\perp_F|,\quad
\mbox{where $v=v^\perp_F +v_F$, $v^\perp_F\in F^\perp$, 
and $v_F\in F$.}
$$

Write $A$ in blocks of the form 
$$
\left(\begin{array}{cc}
                                 A_F & B \\
                                 0 & C \\
                                
                    \end{array}\right).
$$
Since
$C=A/F=(P_{F^\perp} \circ A)/F$,
where $P_{F^\perp}$ is the orthogonal projection on $F^\perp$,
we have the following lemma: 

\proclaim{Lemma} \label{l.perturbations} Given 
any $\varepsilon>0,$
\vglue4pt\noindent \hskip1em \hangindent =20pt \hangafter=1 $\bullet$ 
every $\varepsilon$\/{\rm -}\/perturbation of $A_F$ 
is the restriction of an $\varepsilon$\/{\rm -}\/perturbation of $A$ keeping invariant the other eigenvalues {\rm (}\/but not
necessarily the eigenspaces{\rm ).} Actually{\rm ,} $A/F$ is not modified.
\vglue4pt\noindent \hskip1em \hangindent =20pt \hangafter=1 $\bullet$ 
Any $\varepsilon$\/{\rm -}\/perturbation of $A/F$ is the quotient of an
$\varepsilon$\/{\rm -}\/perturbation of $A$ with $A_F$   invariant.

\endproclaim

The definition of domination of a linear system 
(recall Definition~\ref{d.dominated2}) has a direct generalization for pairs of invariant subbundles. 
Suppose that $E$ and $F$ are two invariant subbundles of a linear system
$(\Si,f,\cE,A)$; we say that {\em  $E$ is $l$-dominated by $F$\/} if for every $x\in\Si$,
 $$
\|A_E^{(l)}(x)\|\,\|A_F^{(-l)}(f^l(x))\|<1/2.
$$
 In this case we write $E\prec F$ or $E\prec_l F$. It is easy to see that this dominance implies that for any $x\in\Si$ the intersection $E(x)\cap F(x)$ is reduced to the zero-vector. So this definition is equivalent to
saying that 
$E\oplus F$ is an $l$-dominated splitting of the system $A_{E\oplus F}$. 

In what follows we will often use  the next  very easy  lemma whose proof we omit:

\proclaim{Lemma} \label{l.domiangle}
For any $K>0$ and $l>0$ there are $K_0>0${\rm ,} 
$l_0>0${\rm ,} and $K_1>0$ satisfying the following property. 

Consider a linear system $(\Si, f, \cE, A)$  bounded by $K$ with an invariant splitting $E\prec_lF$. 
Then  there is a linear change of coordinates
$P${\rm ,} bounded by $K_0${\rm ,}
such that 
the bundles $P(E)$ and $P(F)$ are orthogonal
and invariant by 
the system defined by $B= P\circ A\circ P^{-1}$. Moreover{\rm ,} $P(E)\prec_{l_0}P(F)$ {\rm (}\/for $B${\rm )}
and $B$ is bounded by $K_1$.
\endproclaim

One of the main difficulties
in the proof of Proposition~\ref{p.rank}  comes from the fact that the dominance has not a good behavior if one consider sums of dominated subbundles. The  following remark illustrates this difficulty:

\numbereddemo{{R}emark} There exist linear 
maps and invariant bundles
$\cE=E\oplus F\oplus G$ such that $E\prec F$ and $E\prec G$, but the splitting $E\oplus (F\oplus G)$ is not dominated.
\enddemo
This difficulty comes from 
the relationship between dominance and 
{\em angles\/} (the angles between two bundles of a dominated
splitting cannot be {\em very\/} small), and
the following easy geometric observation:  One can have
simultaneously
{\em big\/} 
angles between $E$ and $F$, $\angle (E,F)$, and
$E$ and $G$, $\angle(E,G)$, but an arbitrarily small angle
$\angle(E,F\oplus G)$. See the next example, and Lemma
\ref{l.domiangle} and Remark~\ref{r.quotient} below. 

\vglue8pt {\it Example {\rm 2}}.
Consider the quotient $\Si$ of $\ZZ\times \NN$ by the relation $(n,m)=(n+3m,m)$ 
 and the map $f\colon \Si \to \Si$ given by
$(n,m)\mapsto (n+1,m)$. 
Denote by $\cE$ the trivial bundle of fiber $\RR^3$. 
Define the linear map $A$ (acting on $\RR^3$) by 

$$
A(n,m)= 
\left\{
\begin{array}{lll}
&\left(\begin{array}{ccc}
                                 1 & 0& 0 \\
                                 0 & 4& 0 \\
                                 0 & 0& 4 \\
                    \end{array}\right) 
\quad
\mbox{ if 
$-m\leq n<0$, },\\ \noalign{\vskip5pt}
&  
\left(\begin{array}{ccc}
                                 1 & 0& 0 \\
                                 0 & 1/4& 0 \\
                                 0 & 0& 4 \\
                    \end{array}\right)  
\quad
\mbox{ if  $0\leq n<m$,}\\ \noalign{\vskip5pt}
& 
\left(\begin{array}{ccc}
                                 1 & 0& 0 \\
                                 0 & 64& 0 \\
                                 0 & 0& 4 \\
                    \end{array}\right)  
\quad
\mbox{ if $m \leq n<2m$}. 
\end{array}
\right.
$$
Then  $(\Si,f,\cE,A)$ is a periodic linear system.

Observe  that the directions $e_1 \RR$, $(e_2+e_3)\RR$,
and $e_3\RR$ define a splitting of $\cE(0,m)$ by  eigenspaces of the matrix $M_A(0,m)$. We now define $E$, $F$, and $G$
at the point $(n,m)$ as the images of these spaces by $A^n$.

By construction, $\cE_0= E\oplus F\oplus G$ is an invariant splitting and
$A$ induces an isometry on $E$ and a
 dilation
of  ratio bigger than $2$
on $G$ and $F$.
So  we get $E\prec_1F$ and $E\prec_1 G$.
However, this system does not admit
any dominated splitting: There are vectors $e_2\in F\oplus G$
which are contracted during an arbitrarily large time.
\vglue8pt

The following lemma explains how we solve the difficulty above.

\proclaim{Lemma} \label{l.transitivity}   With any $K>0$ and $l\in \NN${\rm ,}
   there exists $L$ with the following property\/{\rm :} 
Given  any   linear system $(A,f,\cE,\Si)$ 
such that  $\|A\|$ is bounded by $K$
with an invariant splitting $E\oplus F\oplus G${\rm ,} one has  
\begin{itemize}
\ritem{1.} $(E\prec_l F$  and $E/F \prec_l G/F) \Rightarrow E\prec_L (F\oplus G)${\rm ,}
\ritem{2.} $(F \prec_l G$ and $E/F \prec_l G/F) \Rightarrow (E\oplus F)\prec_L G)$. 
\end{itemize}

\endproclaim

\demo{Proof}
Observe that Lemma \ref{l.domiangle} allows us to assume 
that
(up to a bounded change of coordinates) 
 $E$ is orthogonal to $F$ and  $E/F$ is orthogonal to $G/F$ 
(in the quotient space).
 This means that $E$ is orthogonal to
$P_{F^\perp}(G)$; thus $E$ is also orthogonal to $G$.

Hence 
at each point $x\in\Si$ we can choose
 an orthonormal basis
$(e_i(x))$ of $\cE_x$ such that $e_i(x)\in E(x)$ for  
$i\leq \dim (E(x))$, 
and $e_i(x)\in F(x)$ for $\dim (E(x))< i\leq \dim (E(x))+\dim (F(x))$. 
Moreover, $G(x)$ is contained in the 
subspace spanned by 
the $e_i(x)$ for $i> \dim(E(x))$.
So this subspace is $F(x)\oplus G(x)$ and it is invariant by $A$.
Using the basis $(e_i(x))$  
we can write the matrices of the
system  $A$ as follows,
$$
A(x)=\left(\begin{array}{ccc}
                                 A_E(x)& 0 & 0\\
                                  0  & A_F (x)& B(x)\\
                                  0 & 0 & C(x)
                    \end{array}\right),
$$ 
where $C$ is $A/(E\oplus F)$ and
$B$ is bounded by a constant $K_1$ depending only on
$K$ and $l$.
Write 
$$
D=A_{F\oplus G} = \left(\begin{array}{cc}
                                 A_F& B\\
                                  0  & C\\
                         \end{array}\right).
$$
We can now consider $D$ as a linear system over $F\oplus G$, 
which allows us to iterate $D$.
 We prove that there is $L$ (depending on $K$ and $l$) such that 
$$
\|A_E^L(x)\|\, \|D^{-L}(x)\|<1/2;
$$
that is exactly the first assertion in the lemma
($E\prec_L F\oplus G$).  
To prove this claim observe that 
$$
D^{-l}=\left(\begin{array}{cc}
                                 A_F^{-l}& H 
\\
                                  0  & C^{-l}\\
                    \end{array}\right),
\quad
\mbox{where}
\quad
H = -\sum_{j=1}^{l} A_F^{-j}\circ B \circ C^{j-l-1}.
$$
Therefore  
all the matrices in the definition of $D^{-l}$ are bounded by some constant $K_2$ depending on $K$ and  $K_1$, and so depending on $K$ and $l$.
An elementary calculation shows that for $i\geq 1$
$$
D^{-i\, l} = \left(\begin{array}{cc}
                                 A_F^{-i\, l}& \sum_0^{i-1} A_F^{-j\, l}\circ H\circ C^{(j+1-i)\, l}\\
                                  0  & C^{-i\, l}\\
                    \end{array}\right).
$$
By the hypotheses, 
$$
\|A_E^l\|\, \|A_F^{-l}\|<1/2
\quad 
\mbox{and}
\quad \|A_E^l\|\, \|C^{-l}\|<1/2,
$$ 
the last inequality 
follows immediately from $E/F\prec_l G/F$.
Thus 
$$
\sup\{\|A_F^{-i\,l}\|,\, \|C^{-i\,l}\|\}
<
\frac{1}{2^i \, \|A_E^{i\, l}\|}.
$$ 
An elementary  estimate shows that the norm 
of a matrix of the form
                   $\left(\begin{array}{cc}
                                 X& Y\\
                                  0  & Z\\
                         \end{array}\right)$ is bounded by $\sup\{\|X\|,\|Z\|\} +\|Y\|$. 
So we get that 
$$
\|D^{-i \, l}\|
\le
\frac{1}{2^i\, \|A_E^{i\, l}\|}\, 
\left(1 + 2\,i\, K_1 \right).
$$
It is now immediate that
$$
\|A_E^{i\,l}\| \, \|D^{-i\, l} \| < \frac1{2^i} (1+ (2\, i \,K_1)).
$$
Taking $i$ big enough, this product is clearly less than $1/4$, ending
the proof of our claim.
This completes the first part of the lemma, the second one follows analogously, so we  omit it.
\enddemo

\numbereddemo{{R}emarks} \label{r.quotient}
\begin{itemize}
\item The relations $\prec$ and $\prec_l$ are (strict and partial) order
relations on the set of $A$-invari\-ant subbundles of a linear system 
$(\Si,f,\cE, A)$: The transitivity and the strict
 antisymmetry of these relations are clear. 
\item 
Let $E$ and $G$ be two subbundles of $\cE$, 
 the {\it angle $\angle(E,G)$
between $E$ and $G$\/}
is the infimum of the angles $(u,v)$, where
$u\in E_x, v\in G_x, x\in\Si$.
By Lemma \ref{l.domiangle},
 if  $E\prec_l G$ then
 the angle $\angle(E,G)$ between $E$ and $G$ 
is  greater than some constant $\alpha$ depending on 
$K$ and $l$. 
\item 
Let  $E$ and $F$ be two subbundles such that $E\prec F$ for $A_{E\oplus F}$. Let $G$ be
an $A$-invariant subbundle such that 
$G\cap E= G\cap F = 0$ and that the angles
$\angle(E,G)$ and $\angle(F,G)$ are
 bigger than some $\alpha>0$. Then 
$E/G$ is dominated by $F/G$. 
Moreover,  the constant of domination 
depends only on the bound $K$ of $A$, the constant of 
domination of the splitting $E\prec F$, and the angle $\alpha$.  The easy idea of the proof is the following:
As the angle $\angle(E,G)$ is  greater than $\alpha$,
  the metrics of $E/G$ 
(i.e.\ the metric in $G^\perp$)
and  of $E$ are obtained by a bounded change of coordinates. The same is true for $F$. 
It is now enough  to go back to the definition of domination.
\end{itemize}

\enddemo

\vglue-12pt
As a corollary we now get, 

\proclaim{Lemma} 
\label{l.sum}
Let $K>0$ and $l\in \NN$. 
There is $L$ such that for any linear system $(\Si,f,\cE,A)$ bounded by $K$ and any
invariant subbundles $E${\rm ,} $F${\rm ,} and $G${\rm ,}
 such that $E\prec_l F$ and $F\prec_l G${\rm ,}  
 one has $E\prec_L F\oplus G$
and $E\oplus F\prec_L G$. 
\endproclaim

\demo{Proof} As   $E\prec_lF$ and $F\prec_l G$ there is a constant $\alpha$ (depending only on $K$ and $l$)
such that the angles $\angle(E,F)$ and $\angle(F,G)$ 
are greater than 
$\alpha$. By transitivity of the relation $\prec_l$ one has $E\prec_l G$. 
So, by the third part of the remark,  we get  $l_1$
 (depending on $K$, $l$, and $\alpha$) such that 
$E/F\prec_{l_1} G/F$. Now the corollary follows directly from Lemma~\ref{l.transitivity}. 
\enddemo

4.2. {\it The finest dominated splitting}.
Lemma \ref{l.transitivity} allows us to define the notions of {\em undecomponible\/} and {\em finest dominated splittings.\/} Let us
recall that an invariant splitting 
$\cE = E_1\oplus\cdots\oplus  E_k$ is  dominated if 
for any $1\leq j\leq k-1$ one has $\bigoplus_1^j E_i\prec\bigoplus_{j+1}^k E_i$. 

\proclaim{{C}orollary} 
\label{c.prec}
A splitting $\cE=\bigoplus_{i=1}^k E_i$ is
dominated if and only if $E_i\prec E_{i+1}$
for every $i\in\{1,\dots,k-1\}$ . 
\endproclaim

\demo{Proof}
Observe that if $\cE=\bigoplus_{i=1}^k E_i$ is dominated
then one has 
$E_i\prec E_{i+1}$ straightforwardly from the
definition.
To prove the converse we
argue inductively (on 
the number $k$ of subbundles of the splitting). For $k=3$ the 
corollary  is exactly Lemma~\ref{l.sum}. Assume now that the lemma is true for $k-1$.

We have to prove that if $E_i\prec E_j$ for every $i<j$
then the splitting is dominated; that 
is,  $\bigoplus_1^i E_j\prec \bigoplus_{i+1}^k E_j$
for all $0<i<k$. 
Assume first that $i>1$. Applying 
the induction hypothesis to the restrictions of $A$ to $\bigoplus_2^k E_j$ and $\bigoplus_1^i E_j$,
 we get $\bigoplus_2^i E_j\prec\bigoplus_{i+1}^k E_j$
and $E_1\prec \bigoplus_2^i E_j$.
Now the dominance $\bigoplus_1^i E_j\prec \bigoplus_{i+1}^k E_j$
follows when we apply  Lemma~\ref{l.sum} to the 
subbundles $E_1$, $\bigoplus_2^i E_j$,
and  $\bigoplus_{i+1}^k E_j$. 
\vglue4pt
Finally,
if $i=1$ we apply the induction hypothesis to $\bigoplus_1^{k-1} E_j$
and 
$\bigoplus_2^{k} E_j$.
This gives
$E_1\prec \bigoplus_2^{k-1} E_j$
 and $\bigoplus_2^{k-1}E_j\prec E_k$.
Now the conclusion follows as above
applying Lemma~\ref{l.sum} to these three subbundles.
\enddemo

Corollary~\ref{c.prec} above allows us to use the notation $E_1\prec E_2\prec \cdots\prec E_k$ to denote a dominated splitting $\bigoplus_1^{k} E_j$.   

\proclaim{Lemma} \label{l.subsplitting}
Let $E_1\prec \cdots \prec E_k\prec F\prec G_1\cdots\prec G_m$ be a 
dominated splitting such that $A_F$ admits a dominated splitting $F_1\prec F_2$.
 Then $E_1\oplus\cdots\oplus E_k\oplus F_1\oplus F_2\oplus G_1\oplus\cdots\oplus G_m$ is a dominated splitting.
\endproclaim

\demo{Proof} It suffices to verify that $E_k\prec F_1$ and  
$F_2\prec G_1$.
\enddemo

\numbereddemo{Definition} \label{d.finer}
Let  $\bigoplus_1^k E_i$ and $\bigoplus_1^m F_j$ be two dominated splittings of $\cE$, 
$\cE=\bigoplus_1^k E_i$ and $\bigoplus_1^m F_j$. 
We say that the splitting $\bigoplus_1^m F_j$ 
{\em is 
finer\/} than  
 $\bigoplus_i^k E_i$
 if every $E_i$ is the direct sum of some of the $F_j$,
i.e.\ $E_i=\bigoplus_{j_1}^{j_2} F_j$. 
In this case we write $\bigoplus_1^m F_j\sqsubseteq 
\bigoplus_1^k E_i$.
\enddemo
 
\numbereddemo{{R}emark}
Let $(\Si,f,\cE,A)$ be a linear system of dimension $N$. 

\begin{itemize}
\item
 The relation $\sqsubseteq$ is an order relation on the set of dominated splittings of a linear system $(\Si,f,\cE,A)$.
\item
  Given any finite sequence $n_i$ with
$\sum_i n_i=N$, there is at most one dominated splitting $E_1\prec\cdots\prec E_k$
such that $\dim E_i= n_i$ for every $i$ (actually this easy assertion is a direct consequence of the next lemma).
As a consequence,  the set of 
dominated splittings of $(\Si,f,\cE,A)$ is finite. So there 
are splittings which are minimal for the relation $\sqsubseteq$, these splittings are called 
{\em undecomponible\/}. 
Moreover, given any splitting 
we can subdivide it to get a minimal one, i.e.\ every splitting is bigger, $\sqsupseteq$, 
than at least one undecomponible splitting.
\item 
Following Lemma~\ref{l.subsplitting}, a dominated splitting $\bigoplus_1^k E_i$ is undecomponible if and only if 
none of the subsystems $A_{E_i}$ admits dominated decomposition. 
\end{itemize}

\enddemo 

\proclaim{Proposition} \label{p.thefinest}
For every linear system $(\Si,f,\cE,A)$
 there is a unique undecomponible dominated splitting{\rm ,} called {the finest dominated splitting}.  
\endproclaim  

{\it Proof}.
The first step to prove this proposition is the following
lemma. 

\proclaim{Lemma} 
\label{l.E1}
Let $E_1\prec\cdots\prec E_k$ and $F_1\prec\cdots\prec F_m$
be two dominated splittings of $\cE$.
Then either  $E_1\subset F_1$ and $\bigoplus_2^m F_i\subset 
\bigoplus_2^k E_i$,  or $F_1\subset E_1$ and $\bigoplus_2^k E_i\subset \bigoplus_2^m F_i$.
 As a consequence{\rm ,} if $E_1 \neq F_1$ then 
either $E_1$ or $F_1$ admits a dominated splitting.
\endproclaim 

\demo{Proof}  First we prove that for any $x\in\Si$ one of the linear spaces $E_1(x)$ and $F_1(x)$ is contained in the other
one. 
We argue by contradiction: suppose,
contrary to our hypotheses, that there are 
$$
x\in\Si,\quad
u\in (E_1(x)\setminus F_1(x)),
\quad
\mbox{and}
\quad v\in (F_1(x)\setminus E_1(x)).
$$
 As $\bigoplus_1^k E_i$ is dominated, and $u\in E_1$ and $v\notin E_1$, for every $i$ big
enough we get 
$$
\frac {\|A^{(i)}(x)(u)\|}{\|u\|} \leq 
\frac {\|A^{(i)}(x)(v)\|}{2\,\|v\|}. 
$$
Similarly,  the dominance of 
$\bigoplus_1^m F_i$, $u\notin F_1(x)$, and $v\in F_1(x)$ 
imply (for big $i$)
$$
\frac {\|A^{(i)}(x)(v)\|}{\|v\|}
 \leq \frac {\|A^{(i)}(x)(u)\|}{2\|u\|}.$$ 
These two inequalities give a contradiction. 
So we have proven that for any $x\in\Si$ one has $E_1(x)\subset F_1(x)$ or vice versa. As the dimensions of these spaces do not depend on $x$, we have that
either $E_1\subset F_1$ or $F_1\subset E_1$.
   
The same arguments give the other inclusion;
that is, either $\bigoplus_2^m F_i\subset \bigoplus_2^k E_i$,  or  $\bigoplus_2^k E_i\subset \bigoplus_2^m F_i$.

Finally, assume now that $E_1\neq F_1$ and, for instance, $E_1\subset F_1$. 
Then $F_1$ is transverse to $\bigoplus_2^k E_i$,  so  that
we get  a splitting $F_1= E_1\oplus (F_1\cap\bigoplus_2^k E_i)$. This splitting is clearly invariant and dominated.  
This finishes the proof.
\enddemo

We are now ready to prove Proposition~\ref{p.thefinest}:

Let $\bigoplus_1^k E_i$ and $\bigoplus_1^m F_j$
be two undecomponible dominated splittings. Using Lemma~\ref{l.E1} above we get that $E_1=F_1$ and that $\bigoplus_2^k E_i =\bigoplus_2^m F_j$.
Consider now  the restriction $B=A_{\bigoplus_2^m F_j}$. 
Now
$E_2\prec\cdots\prec E_k$ and $F_2\prec\cdots\prec F_m$ are two undecomponible dominated splittings of $B$, so $E_2= F_2$.
We  conclude the proof repeating $k$ times this argument.  
\hfill\qed

\vglue12pt 4.3. {\it Transitions and invariant spaces}.

\proclaim{Lemma} \label{l.goodtransition}
Let $(\Si,f,\cE,A)$ be a periodic linear system with 
$\varepsilon$\/{\rm -}\/transitions and $\varepsilon_0>\varepsilon$.  Consider an 
invariant subset $\Si_0\subset \Si$ such that there is 
a dominated splitting $E_1\prec\cdots\prec E_k$ defined over $\Si_0$.
Then for every finite subset $\Lambda=\{x_1,\dots,x_m\}$ of $\Si_0$
there are $\varepsilon_0$\/{\rm -}\/transitions such that the corresponding linear maps {\rm (}\/transitions{\rm )} $T_{x_j,x_i}$ map
$E_l(x_i)$ on $E_l(x_j)$ for any $l\in\{1,\dots,k\}$ and $x_j$, $x_i\in \Lambda$. 
\endproclaim

Before proving this lemma let us observe that as a direct consequence of it we get the following corollary.

\proclaim{{C}orollary} \label{c.restr}
Let $(\Si,f,\cE,A)$ be a periodic linear system with transition and $E\oplus F$ a dominated splitting. 
Then the induced systems $A_E$ and $A_F$  admit  \pagebreak transitions.
\endproclaim

\demo{Proof} 
By hypothesis, there are $\varepsilon$-transitions $[t^0_{x_j,x_i}]$ from $x_i$ to $x_j$ for every pair of points in $\Lambda$. 
By Remark~\ref{r.semigroup}, for any $n_1,$
$n_2\geq 0$
the word 
$$
[M]_A(x_j)^{n_2}\, [t^0_{x_j,x_i}]\, [M]_A(x_i)^{n_1}
$$
 is also a transition. Moreover, any 
$(\varepsilon_0-\varepsilon)$-perturbation of an $\varepsilon$-transition
is also an $\varepsilon_0$-transition. 

Taking an arbitrarily small perturbation of 
$[t^0_{x_j,x_i}]$ we can assume that 
its corresponding matrix $\tilde T^0_{j,i}$ maps
$E_k(x_i)$ transversally to $\bigoplus_1^{k-1} E_m(x_j)$.
Thus, by the dominance 
$\bigoplus_1^{k-1} E_m \prec E_k$, taking $n_2$ large enough, 
we have that 
$$
\left((M_A(x_j))^{n_2}\circ \tilde T^0_{j,i}\right)(E_k(x_i))\quad \mbox{is arbitrarily close to}\quad E_k(x_j).
$$ 
So there is a small perturbation 
$M_{\tilde{A}}(x_j)$ 
of
$M_{{A}}(x_j)$ 
such that 
$$
\left( M_{\tilde{A}}(x_j)\circ (M_A(x_j))^{n_2}\circ \tilde T^0_{j,i}\right)(E_k(x_i)) = E_k(x_j).
$$
Denote by $T^1_{j,i}$ the linear map corresponding to this new transition. 
Now we consider the pre-image of $E_k(x_j)$
by $T^1_{j,i}$  and observe that 
$$
\left(M_A(x_i)^{-n_1}\circ (T^1_{j,i})^{-1}\right)(E_k(x_j))= E_k(x_i).
$$
 As above, since 
$\bigoplus_1^{k-1} E_m \prec E_k$, 
for $n_1$ sufficiently large, 
$$
\left(M_A(x_i)^{-n_1}\circ (T^1_{j,i})^{-1}\right)(\bigoplus_1^{k-1} E_m(x_j))\quad\mbox{is arbitrarily close to}\quad \bigoplus_1^{k-1} E_m(x_i).
$$
 So, arguing as before,
the announced transition $T_{j,i} $ is obtained composing (at the right) $T^1_{j,i}\circ M_A(x_i)^{n_1}$ with a small
perturbation of $M_A(x_i)$ mapping $E_k(x_i)$ into
 itself and $\bigoplus_1^{k-1} E_m(x_i)$ into  $\left(M_A(x_i)^{-n_1}\circ (T^1_{j,i})^{-1}\right)(\bigoplus_1^{k-1} E_m(x_j))$.

Next we consider the restriction of $A$ to $\bigoplus_1^{k-1} E_m$ to get a new perturbation of the transition mapping $E_k(x_i)$, $E_{k-1}(x_i)$, and $\bigoplus_1^{k-2}E_m(x_i)$
into the corresponding spaces for $x_j$. The inductive pattern
of our construction is now clear. 
\enddemo

4.4. {\it Diagonalizable systems. }

\numbereddemo{Definition}
 A  periodic linear system  $( \Si,f,\cE,B)$ of dimension $N$ is {\em diagonalizable} if  
for every $x\in \Sigma$
the matrix $M_B(x)$ has only real positive eigenvalues with multiplicity $1$.

Denote by $\lambda_1(x)<\cdots<\lambda_N(x)$
the eigenvalues of $M_B(x)$ and by $E_i(x)$ the one-dimensional eigenspace corresponding to $\lambda_i(x)$. Then
$\cE=\bigoplus_1^N E_i$ is an {\it invariant splitting} of $B$.
\enddemo

The next result says that every periodic system
with transitions can be approximated by a diagonalizable
subsystem (defined on a dense subset). More precisely,

\proclaim{Lemma} \label{l.real}
Let $(\Si,f,\cE,A)$ be a 
periodic linear  system with transitions.
Then for any $\varepsilon>0$ there is a diagonalizable 
$\varepsilon$\/{\rm -}\/perturbation $\tilde A$ of $A$
defined on   a dense invariant subset $\tilde \Si$ of $\Si$. 
\endproclaim

\demo{Proof} 
Observe first that any matrix $M\in {\rm GL}(n,\RR)$ 
can be perturbed (by an arbitrarily small perturbation)
to get $\tilde M$ 
having only eigenvalues
(complex or real) with multiplicity one and such that
any pair of eigenvalues with the same modulus are
complex and conjugate, having rational argument.
So there is $k$ such that $\tilde M^k$ has only real positive eigenvalues, but some of them may have multiplicity $2$
(those coming from a complex eigenvalue of $\tilde M$). 
As above, by a small perturbation of $\tilde M^k$, we get a matrix $M_1$ having only real positive eigenvalues with multiplicity $1$.

Take $x\in \Si$ and write $M=M_A(x)$.
By the comments above, there is an $\varepsilon/10$-perturba\-tion $\tilde M$ 
of $M$ having only multiplicity one eigenvalues. Denote
by $M_1$ the perturbation of an appropriate 
power of $\tilde M$ as the one  built above
(i.e having only real eigenvalues with multiplicity one).
Let $T$ be an\break $\varepsilon/10$-transition matrix from 
$x$ into itself. So for any positive $n_1$
and $n_2$ there are $y\in \Sigma$ and a  $3\varepsilon/10$-perturbation $\tilde A$
of $A$ along the orbit of $y$ such that 
$$
M_{\tilde A}(y)= M_1^{n_2}\circ T\circ M_1^{n_1}.
$$
So, taking $n_1$ and $n_2$ sufficiently big, we can now repeat the proof of 
Lemma~\ref{l.goodtransition} 
to get a $4\varepsilon/10$-transition $T_1$ preserving 
all the eigenspaces of $M_1$. These spaces are $1$-dimensional, thus they are the eigenspaces of $T_1$ and
their corresponding  eigenvalues are real. 

Recall (see Remark~\ref{r.semigroup}) that if $T_1$ is a transition from $x$ into itself then $T_1^2$ is also a transition. Therefore,
replacing if necessary $T_1$ by $T_1^2$,
we can assume that $T_1$ has only real positive eigenvalues
with multiplicity one. 
Observe that
$T_1$ and $M_1$ have the same one-dimensional eigenspaces;
thus 
$T_1$ commutes with
$M_1$.
So, for $k$ big enough, $M_1^k\circ T_1$ has only positive real
eigenvalues with multiplicity $1$. Moreover, there are a point $z$
(whose orbit passes arbitrarily close to $x$) and an $\varepsilon/2$-perturbation $\hat A$ of $A$ along the orbit of $z$ such that
$$
M_{\hat A}(z)=M_1^k\circ T_1.
$$  
This proves that the set $\tilde \Si$
of points $z\in\Si$ 
admitting an $\varepsilon$-perturbation along its orbits such that $M_{\hat A}(z)$ has only positive real eigenvalues with
multiplicity
$1$ is dense in $\Si$, ending the proof of the lemma.
 \enddemo

In what follows, to prepare the proof of Theorem~4, we will give 
 some volume controlling versions of our lemmas, as we do in  the next remark:

\numbereddemo{{R}emark} \label{r.real}
Let $(\Si,f,\cE,A)$ be a linear system with transitions such that there is some point $x_0\in\Si$ such that the modulus
of the jacobian of the matrix $M_A(x_0)$, denoted by $J(M_A(x_0))$, is greater than one. 
Then we can choose the diagonalizable $\varepsilon$-perturbation $\tilde A$ of $A$ given by Lemma~\ref{l.real} such that the modulus of the jacobian $J(M_{\tilde A}(p))$ 
is  bigger than one for every $p$ in the dense subset 
$\tilde \Si$ of $\Si$.
\enddemo

\demo{Proof} First, by using the transitions and the existence of the point $x_0$ with $J(M_A(x_0))>1$, we get a dense subset $\hat
\Si$ of $\Si$ such that $J(M_A(x))>1$ for all $x\in\hat\Si$. Now we can repeat the proof of Lemma~\ref{l.real} above: Taking the
exponents $n_1$ and $n_2$ large enough we get  that (in the proof of the lemma) the moduli of the jacobian of the corresponding
matrices $M_{\tilde A}(y)$ are greater than one.
\enddemo

\section{Dominated splittings, complex eigenvalues\\  of rank $(i,i+1)$, and homotheties}
\label{s.proofprop}

In  this section  we prove Propositions~\ref{p.rank} and \ref{p.homothety}. To prove Proposition~\ref{p.rank}, we use
Lemma~\ref{l.transitivity} and consider successive quotients of linear systems in order to get complex eigenvalues
(obtained by using the arguments of the $2$-dimensional case, see
Section~3). Then the transitions and the existence of complex
eigenvalues of any rank allow us to ``mix all the eigenvalues" of some matrix, obtaining the 
homothety announced in Proposition~\ref{p.homothety}.

\vglue12pt 5.1. {\it  Getting complex eigenvalues of any rank}.

\proclaim{Lemma} \label{l.oplus} Given  $K>0$ and $\varepsilon>0$ there is $l\in \NN$ such that 
for any diagonalizable linear periodic system $(\Si,f,\cE,B)$ of dimension $N$ and bounded by $K${\rm ,}
 and any $1\leq i\leq N-1$ one has\/{\rm :}
\begin{itemize}
\item
Either there is an $\varepsilon$-perturbation of $B$ having a complex eigenvalue of rank $(i,i+1)${\rm ,}
\item
or 
$$
E_j/({E_{j+1}\oplus \cdots \oplus E_k}) \prec_l E_{k+1}/({E_{j+1}\oplus \cdots \oplus E_k})
$$
for every $j\le i\le k$.
\end{itemize}

\endproclaim

\demo{Proof}
Fix $\varepsilon>0$ and let $l$ be the dominance constant given by  Proposition~\ref{p.2dim}. 
If
$E_j/({E_{j+1}\oplus \cdots \oplus E_k})$
is $l$-dominated by  $E_{k+1}/({E_{j+1}\oplus \cdots \oplus E_k})$
we are done. 
Otherwise, by Proposition~\ref{p.2dim}, 
we can perturb the quotient to get eigenvalues
 $\tilde\lambda_j = \tilde\lambda_{k+1}=\alpha$.
Moreover,  
by Lemma \ref{l.perturbations}, this perturbation
of the quotient gives a perturbation $\tilde B$ of $B$ having a pair of
eigenvalues 
 $\tilde\lambda_j = \tilde\lambda_{k+1}=\alpha$
and preserving the  eigenvalues of the restriction of $B$ to
$E_{j+1}\oplus \cdots \oplus E_{k}$
(we denote such eigenvalues by $\lambda_i$).

Consider  a small isotopy $B_t$ producing this perturbation (i.e.\ $B_0=B$ and $B_1=\tilde B$) and denote by $\lambda_{j,t}$ and $\lambda_{k+1,t}$ 
the continuations of the eigenvalues $\lambda_j$  and $\lambda_{k+1}$ at time $t$, so that 
$$
\lambda_{j,0}=\lambda_j,
\quad
 \lambda_{j,1}=\tilde\lambda_j,
\quad
\lambda_{k+1,0}=\lambda_{k+1},
\quad \mbox{and}
\quad 
\lambda_{k+1,1}=\tilde\lambda_{k+1}.
$$
 Moreover, we can assume that for
every  $0\leq t<1$ one has $\lambda_{j,t}<\lambda_{k+1,t}$ (otherwise we stop the isotopy at the first $t$ such that $\lambda_{j,t}=\lambda_{k+1,t}$). Then, by continuity, following this isotopy there are three possibilities:

\begin{itemize}
\item 
at some  stage
$t$ of the isotopy, $\lambda_{j,t} = \lambda_{i+1}\le \lambda_{k+1,t}$, 
\item
at some stage $t$ of the isotopy, $\lambda_{j,t}\leq \lambda_{i} = \lambda_{k+1,t}$, 

\item  for $t=1$, $\lambda_i<\lambda_{j,1}=\alpha=\lambda_{k+1,1}<\lambda_{i+1}$.
\end{itemize}
In each of these cases a small perturbation gives a complex eigenvalue of rank $(i,i+1)$, ending the proof of the lemma.
\enddemo

To deduce Proposition~\ref{p.rank} from Lemma~\ref{l.oplus} above we will apply successively
Lemma~\ref{l.transitivity}. 
Even if the arguments in the
proof of the proposition  are rather simple 
and the main difficulty of it 
is of combinatorial type, for clearness
and to organize the combinatorics,
 we divide the proof into two steps: 

\proclaim{Lemma} \label{l.croissant} For any $K>0$, $N>0${\rm ,} and $\varepsilon>0$
 there is $L_0\in \NN$ such that 
for every  diagonalizable periodic linear system 
$(\Si, f, \cE, B)${\rm ,}
of dimension $N$ and bounded by $K${\rm ,}
and for any $1\leq i\leq N-1$\/{\rm :}
\begin{itemize}
\item 
Either $B$ admits an $\varepsilon$-perturbation having a complex eigenvalue of rank $(i,i+1)$
\item
or  $$E_j/_{\bigoplus_{j+1}^i E_{m}} \prec_{L_0}  \bigoplus_{i+1}^N E_{k}/_{\bigoplus_{j+1}^i E_m} $$
for every $j\le i$.
\end{itemize}

\endproclaim

Before proving the lemma
 observe that if $j=i$ the second item of the lemma means that $E_i \prec_{L_0}\bigoplus_{i+1}^N E_{k}$.

\demo{Proof} 
Assume 
that every $\varepsilon$-perturbation of $B$ has 
no  complex eigenvalues of rank $(i,i+1)$. From Lemma~\ref{l.oplus}, taking
$k=i$ and $k=i+1$   we get
$l$ such that for every $j<i$  one has
$$
E_j/_{\bigoplus_{j+1}^i E_m} \prec_l E_{i+1}/_{\bigoplus_{j+1}^i E_m}
\quad
\mbox{and}
\quad
E_j/_{\bigoplus_{j+1}^{i+1} E_m} \prec_l E_{i+2}/_{\bigoplus_{j+1}^{i+1} E_m}.
$$
We now apply 
Lemma~\ref{l.transitivity} 
taking
$$
E=E_j/_{\bigoplus_{j+1}^i E_m},
\quad
F=E_{i+1}/_{\bigoplus_{j+1}^i E_m},
\quad
\mbox{and}
\quad
G=E_{i+2}/_{\bigoplus_{j+1}^i E_m}.
$$
Observe first that
$$
E/F =E_j/_{\bigoplus_{j+1}^{i+1}E_m}
\quad
\mbox{and}
\quad 
G/F=E_{i+2}/_{\bigoplus_{j+1}^{i+1}E_m}.
$$
By the comments above, $E\prec_l F$ and $E/F \prec_l G/F$.
Thus by Lemma~\ref{l.transitivity}
there is  $l_1$ such that 
$$
E_j/_{\bigoplus_{j+1}^i E_m} \prec_{l_1} (E_{i+1}\oplus E_{i+2})/_{\bigoplus_{j+1}^i E_m}.
$$
Moreover, this holds for every $j\leq i$. 
Taking $k=i+2$ in  Lemma~\ref{l.oplus}, one  gets 
$$
E_j/_{\bigoplus_{j+1}^{i+2} E_m} \prec_l E_{i+3}/_{\bigoplus_{j+1}^{i+2} E_m}.
$$
As before, applying  Lemma~\ref{l.transitivity} to $$
E=E_j/_{\bigoplus_{j+1}^i E_m},
\quad
F=(E_{i+1}\oplus E_{i+2})/_{\bigoplus_{j+1}^i E_m},
\quad
\mbox{and}
\quad
G=E_{i+3}/_{\bigoplus_{j+1}^i E_m},
$$
one gets $l_2$ such
that 
$$
E_j/_{\bigoplus_{j+1}^i E_m} \prec_{l_2} (E_{i+1}\oplus E_{i+2}\oplus E_{i+3})/_{\bigoplus_{j+1}^i E_m}.
$$
The inductive procedure 
to prove the lemma is now clear.
\enddemo

\proclaim{Lemma} \label{l.ii+1}
For any $K>0${\rm ,} $N>0${\rm ,} and $\varepsilon>0$ 
there is $L\in \NN$ such that 
for every  diagonalizable    periodic system
$(\Si, f, \cE, B)${\rm ,}
 of dimension $N$ and bounded by $K${\rm ,} and  any 
$1\leq i\leq N-1$\/{\rm :}
\begin{itemize}
\item 
Either $B$ admits an $\varepsilon$\/{\rm -}\/perturbation having a complex eigenvalue of rank $(i,i+1)${\rm ,}
\item 
or  
$$
\bigoplus_1^i E_j 
\prec_L \bigoplus_{i+1}^N E_j.
$$
\end{itemize}

\endproclaim

\demo{Proof}
Assume that it is not possible 
to perturb $B$ to get a complex eigenvalue of rank $(i,i+1)$.
Once more, the proof is an inductive argument using alternately 
Lemmas~\ref{l.croissant} and \ref{l.transitivity}
above. 

By Lemma~\ref{l.croissant}, 
applied to $j=i$ and to $j=i-1$, one knows
$$
E_i\prec_{L_0} \bigoplus_{i+1}^N E_k
\quad
\mbox{and}
\quad 
E_{i-1}/E_i \prec_{L_0}\bigoplus_{i+1}^N E_k /E_i.
$$
So Lemma~\ref{l.transitivity} gives $L_1$ such that 
 $(E_{i-1}\oplus E_i) \prec_{L_1}\bigoplus_{i+1}^N E_k $.
By Lemma~\ref{l.croissant}, taking $j=i-2$,
 one gets 
$$
E_{i-2}/(E_{i-1}\oplus E_i) \prec_{L_0}\bigoplus_{i+1}^N E_k /(E_{i-1}\oplus E_i).
$$ 
Thus applying Lemma~\ref{l.transitivity} to $E_{i-2}$, $E_{i-1}\oplus E_i$,
and $\bigoplus_{i+1}^N E_k $ one obtains
$$
(E_{i-2}\oplus E_{i-1}\oplus E_i) \prec_{L_1}\bigoplus_{i+1}^N E_k .
$$
The lemma now follows
by a very simple induction.
\enddemo

5.2. {\it End of the proof of Proposition}~\ref{p.rank}.
Recall that Proposition~\ref{p.rank}
 says that if there is no perturbation  of the system with a complex
 eigenvalue of rank $(i,i+1)$ then the system has an $l$-dominated splitting of dimension $i$. 

Let $(\Si,f,\cE,A)$ be a continuous periodic linear system, of dimension $N$ and bounded by $K$, having transitions.
Assume that there are
$\varepsilon>0$ and  $i\in\{1, \dots,
N-1\}$ such that 
every 
$\varepsilon$-perturbation of $A$
has no  complex eigenvalues of rank $(i,i+1)$.

Choose a sequence $\varepsilon_n$ such that $\varepsilon/2> \varepsilon_n\to 0$.
As the system $(\Si,f,\cE,A)$ has transitions, using Lemma~\ref{l.real} 
we get  dense subsets $\Si_n$ of $\Si$ and diagonalizable
$\varepsilon_n$-perturbations $B_n$ of $A$
defined on $\Si_n$.
Note  that 
$B_n$ is a diagonalizable periodic linear system of dimension $N$
bounded by $K_1 = K+\varepsilon/2$.
Thus, by hypotheses,  
it is impossible to create complex
eigenvalues of rank $(i,i+1)$ by  $(\varepsilon/2)$-perturbations of it.

By Lemma~\ref{l.ii+1},  there is $L=L(K_1,\varepsilon/2,N)$ such that
every system $B_n$ above admits an $L$-dominated splitting 
$E_n\oplus F_n$, with
$E_n \prec_L F_n$ and\break $\dim (E_n) = i$. 
Finally, as the system $A$ is continuous, 
the sets $\Si_n$ are dense in $\Si$, and $\| B_n- A\|\to 0$, Lemma~\ref{l.dense}
ensures that $A$ admits an $L$-dominated splitting $E\prec_L F$ with $\dim(E) = i$. The proof of the proposition is now
complete. \phantom{whattodo}
\hfill\qed

\demo{{\rm 5.3.} Proof of Proposition {\rm \ref{p.homothety}}}
The naive idea of the proof of 
Proposition~\ref{p.homothety} is to use the transitions to multiply matrices corresponding to different points of $\Si$ having complex eigenvalues of different rank. 
In this way one distributes 
homogeneously the action of the eigenvalues
in the whole fibers, obtaining homotheties.
The main step of the proof is the following lemma: 
\enddemo

\proclaim{Lemma} \label{l.transposition}
Let $(\Si,f,\cE,A)$
be  a continuous periodic linear 
 system of dimension $N$ with 
transitions. Fix  $0<\varepsilon_0$ and assume that 
for any $i\in \{1,\dots, N-1\}$
there is an $\varepsilon_0$\/{\rm -}\/perturbation of $A$ having a complex eigenvalue of rank $(i,i+1)$.

Then for every $0<\varepsilon_1<\varepsilon_0$ 
there is a point $p\in \Si$ 
such that for every $1\leq i<N$ there is an $\varepsilon_1$\/{\rm -}\/transition  $[t^i]$ from $p$ to itself
with the following properties\/{\rm :}
\begin{itemize}
\item 
There is an $\varepsilon_1$\/{\rm -}\/perturbation $[M]_{\tilde{A}}(p)$
of the word $[M]_A(p)$ such that 
the corresponding matrix $M_{\tilde{A}}(p)$ has only real positive eigenvalues with multiplicity $1$.
Denote  by $\tilde\lambda_1<\cdots<\tilde\lambda_N$
such eigenvalues and by $E_i(p)$ their respective {\rm (}\/$1$\/{\rm -}\/dimensional\/{\rm )} eigenspaces. 
\item 
There is an $(\varepsilon_0+\varepsilon_1)$\/{\rm -}\/perturbation 
$[\tilde t^i]$ of the transition $[t^i]$ such that 
the corresponding matrix $\tilde T^i$ satisfies
\begin{itemize}   
\item 
$\tilde T^i(E_j(p)) = E_j(p)$ if $j\notin\{i,i+1\}$,
\item 
$\tilde T^i(E_i(p)) = E_{i+1}(p)$ and $\tilde T^i(E_{i+1}(p)) = E_i(p)$.
\end{itemize}
\end{itemize}

\endproclaim 
\vglue-8pt
In fact we have also  a stronger (volume controlling) version of this lemma:

\numbereddemo{{R}emark} 
\label{r.transposition}
Under the hypotheses of Lemma~\ref{l.transposition}, suppose in addition that there is a point $x\in\Si$ such that the modulus of the jacobian $J(M_A(x))$ is bigger than one. Then we can choose the point $p$ and the perturbation $\tilde A$ in the lemma such that $J(M_{\tilde A}(p))>1$.
\enddemo

Before proving  Lemma~\ref{l.transposition} and Remark~\ref{r.transposition} 
let us deduce Proposition~\ref{p.homothety} from them.

\vglue12pt 5.3.1. {\it End of the proof of Proposition}~\ref{p.homothety}.
The hypotheses of the proposition (existence of $\varepsilon_0$-perturbations with complex eigenvalues of any rank) imply that we 
can apply Lemma~\ref{l.transposition} to the
system $(\Si, f, \cE, A)$ for all $i\in \{1, \dots, N-1\}$. To prove the proposition we show 
that for any $\varepsilon>\varepsilon_0$
there are an $\varepsilon$-perturbation $\hat A$ of $A$ 
and a point $x\in \Sigma$
for which $M_{\hat A}(x)$ is a  homothety. 

Choosing $0<\varepsilon_1<(\varepsilon-\varepsilon_0)/10$, 
let $p$ and $[t^i]$,  $i\in \{1,\ldots , N-1\}$, be the point and the $\varepsilon_1$-transition, from $p$ to itself, given by
Lemma~\ref{l.transposition}.

Observe that the action of the perturbed 
transition $\tilde T^i$ 
(which are $(\varepsilon_0+\varepsilon_1)$-perturbations of $T_i$)
on the
finite set $\{E_i(p)\}_{1\leq j\leq N}$ of  eigenspaces of $M_{\tilde A}(p)$ is the transposition $(i,i+1)$ which
interchanges $E_i(p)$ and $E_{i+1}(p)$,
keeping invariant the others $E_j(p)$.
Recall that a {\em transposition\/} is an order $2$ permutation;
thus  it is equal to its inverse. Moreover,
the transpositions $(i,i+1)$, $i\in\{1,\dots,N-1\}$,
 generate the group of (all) permutations of the finite set $\{E_j(p)\}_{1\leq j\leq N}$. 
\vglue2pt
Given $0\leq k< N$ denote by $\sigma_k$ the cyclic permutation
defined by\break $\sigma_k(E_j(p))= E_{j+k}(p)$, 
where the sum $i+j$ is considered
in the cyclic group $\ZZ/ \pagebreak N\ZZ$.

As a direct consequence of the previous comments, 
we have that for every $0<k<N$ 
there exists an element $[\tilde S_k]$ in the semi-group generated by the transitions $[\tilde t^i]$
 such that its action on the finite set $\{E_j(p)\}_{1\leq j\leq N}$ is
the permutation $\sigma_k$, i.e.\ if $\tilde S_k$ is the  matrix corresponding to the
word $[\tilde S_k]$ then one has
$\tilde S_k(E_j(p))= \sigma_k(E_j(p))$. 

Let $[S_k]$ be the word of matrices corresponding
to the perturbation $[\tilde S_k]$  
in the semi-group generated by the initials      
$[t^j]$. 
As the $[t^j]$ are $\varepsilon_1$-transitions from $p$ to 
itself, any word in the semi-group generated by the 
$[t^j]$ is also an $\varepsilon_1$-transition (recall
Remark~\ref{r.semigroup}).  In particular, the $[S_k]$ are
$\varepsilon_1$-transitions from $p$ to
itself. 
(For completeness let us write $[S_0]=[S_N]$ the empty word
whose corresponding  matrix is by convention the identity.)

By definition of transitions, for any $n\in \NN$ there is a point
$x_n\in\Si$ such that the word $[M]_A(x_n)$ is $\varepsilon_1$-close to
the word $[W_n]$ defined by 
$$
\begin{array}{ll}
[W_n] =& [S_1]\,[M]_A^n(p)\,[S_{N-1}]\,[S_1]\,[S_{N-1}]\,[S_2]\,[M]_A^n(p)\,[S_{N-2}]\,[S_2]\,[S_{N-2}] \cdots \\ \noalign{\vskip4pt}
&\cdots
[S_{N-2}]\,[M]_A^n(p)\,[S_2]\,[S_{N-2}]\,[S_2]\,[S_{N-1}]\cdots \\ \noalign{\vskip4pt}
&\cdots [S_{N-1}][M]_A^n(p)\,[S_1]\,[S_{N-1}]\,[S_1] \, [M]_A^n(p).
\end{array}
$$
Let us state some properties  justifying the introduction
of this word:
\begin{itemize}
\item[i)]
Given any $i$, the matrix 
$\tilde S_{N-i}\circ \tilde S_i$  acts trivially on the set of spaces $\{ E_j(p)\}$. Let us denote by $\mu_{i,j}$ the
eigenvalue of $\tilde S_{N-i}\circ \tilde S_i$
corresponding to the eigenspace $E_j(p)$.
\item[ii)]
Recall that $\tilde S_i$ maps $E_j(p)$ into $E_{i+j}(p)$
and that $M_{\tilde A}(p)$ is diagonal in the basis
corresponding to the directions $E_k(p)$ and
denote by $\tilde \lambda_k$ the eigenvalue of
$M_{\tilde A}(p)$ corresponding to such a direction.
 Hence
every  $E_j(p)$ is an eigenspace of 
$\tilde S_{N-i}\circ (M_{\tilde{A}}(p))^n\circ \tilde S_i$
whose corresponding eigenvalue is 
$\mu_{i,j}\, \tilde\lambda_{j+i}^n$.
\item[iii)]
By the two items before, for every $j$ the space $E_j(p)$ is an eigenspace of the matrix 
$$
\tilde W_{i,n} = \tilde S_{N-i}
\circ (M_{\tilde{A}}(p))^n\circ \tilde S_i\circ 
\tilde S_{N-i}\circ \tilde S_i
$$  
whose corresponding eigenvalue is $\mu_{i,j}^2\,\tilde\lambda_{j+i}^n$. (Recall that
$\tilde W_{0,n} = M_{\tilde A}(p)^n$ by convention.)
\end{itemize}

Proposition~\ref{p.homothety} is now an immediate
consequence of the following claim:

\nonumproclaim{Claim} 
For every  $n>n_0$ sufficiently large  there is an $\varepsilon$\/{\rm -}\/perturbation
 $\hat A$ of $A$ along the orbit of $x_n$ such that $M_{\hat A}(x_n)$ is a homothety of ratio $\tilde\La^n${\rm ,} where 
$$
\tilde\La=\prod_1^N \tilde\lambda_i= J(M_{\tilde A}(p)). \pagebreak
$$
\endproclaim

\demo{Proof}
Denote by $[\tilde W_n]$ the word obtained from $[W_n]$ 
by putting a $\, \tilde\cdot\, $ above any letter $S$ and $A$; 
this word is an $(\varepsilon_0+\varepsilon_1)$-perturbation of the word $[W_n]$,
so it  is an $(\varepsilon_0+2\,\varepsilon_1)$-perturbation of the word $[M]_A(x_n)$. 
Moreover, using the notation in item (iii) above, we see that the corresponding matrix $\tilde W_n$ is the product
$$
\tilde W_n = \tilde W_{N-1,n}\circ\cdots\circ\tilde W_{1,n}\circ \tilde W_{0,n}.
$$ 
So, by item (iii) above, for every $j$,  the one-dimensional space $E_j(p)$ is an eigenspace of $\tilde W_n$ and its corresponding eigenvalue $\tilde\lambda_{j,n}$ is  
$$
\tilde\lambda_{j,n}=
\prod_{i=0}^{N-1}\tilde\lambda^n_{j+i}\, \prod_{i=0}^{N-1}\mu_{i,j}^2
=
\prod_{i=1}^{N}\tilde\lambda^n_{i}\,
\prod_{i=0}^{N-1}\mu_{i,j}^2.
$$
Writing
$$
C_j= \prod_{i=0}^{N-1}\mu_{i,j}^2 >0
\quad
\mbox{and}
\quad
\tilde\La=\prod_1^N \tilde\lambda_i ,
$$
we get that,
for any $n\in\NN$, the eigenvalue 
$\tilde\lambda_{j,n}$ is 
$$
\tilde\lambda_{j,n}= C_j\, \tilde\La^n.
$$ 
This means that
the matrix $\tilde W_n$ is the product of a 
homothety 
$(\tilde\La^n\cdot {\rm Id})$ with a matrix $B$ which does not depend on $n$ and leaves invariant every one-dimensional space
$E_i(p)$. So the matrix $B$ commutes with every $\tilde W_{i,n}$.
Finally, by construction,  all the  eigenvalues of $B$ are positive. 

Denote by 
$$
C_{n,j}= (C_{j})^{-\frac 1n}.
$$
Clearly,
when $n$ becomes very large the $C_{n,j}$
are arbitrarily close to $1$. Consider the
matrix $B_n$ 
having the  $E_j(p)$ as eigenspaces 
and the $C_{n,j}$ as the corresponding eigenvalues. 
Denote by $[M]_{\hat A}(p)$ the word obtained from 
$[M]_{\tilde A}(p)$ by replacing its first letter 
$\tilde A(p)$ (at the right) by $\tilde A(p)\circ B_n$. For 
$n$ large enough 
this new word is an $\varepsilon_1$-perturbation of 
$[M]_{\tilde A}(p)$, so by item (i) of Lemma~\ref{l.transposition}
it  is also 
a $2\,\varepsilon_1$-perturbation of  $[M]_A(p)$.
Now, the  matrix corresponding to $[M]_{\hat A}(p)$ is 
$M_{\tilde A}(p)\circ B_n$. 
As $B_n$ commutes with $M_{\tilde A}(p)$, and by
the definitions of $B_n$ and  $C_{n,j}$, we get 
that 
$$
(M_{\tilde A}(p)\circ B_n)^n = M_{\tilde A}^n(p)\circ B^{-1}.
$$
As a conclusion, the word $[\hat W_n]$ 
obtained by changing the initial subword $[M]_{\tilde A}^n(p)$
of $\tilde W_n$ by $[M]_{\hat A}(p)$ is $(\varepsilon_0+2\, \varepsilon_1)<\varepsilon$ close to the
word $[M]_A(x_n)$,  and 
its corresponding matrix $\hat W_n=\tilde W_n\circ B^{-1} = \tilde\La^n\cdot {\rm Id}$ is a
homothety.
 This ends the proof of the claim.  
\enddemo

After reading carefully the proof before, one has the following remark  which will play a  key role in controlling the volume in the
next section.

\numbereddemo{{R}emark}
\label{r.rema}
Under the hypotheses of Proposition~\ref{p.homothety}, suppose in addition that there is $x_0\in\Si$ such that the modulus of the jacobian $J(M_A(x_0))$ is bigger than one. Then we can choose the perturbation $\tilde A$ of $A$ and the point $x\in\Si$ in the proposition such that $M_{\tilde A} (x)$ is a homothety with ratio of modulus bigger
 than one. 
\enddemo

{\it Proof}.
It is enough to repeat the proof of Proposition~\ref{p.homothety} bearing in mind the volume controlling Remark~\ref{r.transposition}: So the point $p$ and the first perturbation $\tilde A$ can be chosen such that $1< J(M_{\tilde A}(p)) = \tilde \La$. To conclude the proof, it is now enough to recall that (with the notation of the proof of Proposition~\ref{p.homothety}) $J(M_{\hat A}(x_n))=\tilde \La^n$ (see the claim in this proof). 
\hfill\qed\vglue8pt

To end the proof of Proposition~\ref{p.homothety} 
it remains to prove Lemma~\ref{l.transposition}.
This is done in the next section.

\vglue6pt 5.3.2. {\it Proof of Lemma}~\ref{l.transposition} ({\it and Remark}~\ref{r.transposition}).
The proof of Lemma~\ref{l.transposition} follows  from the ideas of  the proof of Lemma~\ref{l.goodtransition},  based on the
following fact: Given a vector $v$, a pair of  matrices $T$ and $M$,
and the eigenvector $w$ 
associated to the largest (in modulus) eigenvalue of
$M$, 
it is very easy to map $v$ into $w$ by an
arbitrarily  small perturbation of $M^n\circ T$, if 
$n$ is  large enough.
So, using this fact, given a dominated splitting, a simple inductive argument allows us to get transitions  preserving it. 
The difficulty here is that we want
to get  a transition interchanging two spaces of a dominated splitting
(in our case the eigenspaces of a
diagonal matrix). For that we will
use the complex eigenvalues, which enable us 
  to map an arbitrary vector into the eigenvector
corresponding to the weaker eigenvalue. 
Let us explain all that in detail.

 Fix any $0<\varepsilon_1<\varepsilon_0$. We now build some
$\varepsilon_1/10$-perturbation $\tilde A$ of $A$, modifying the initial system along a finite number of orbits.
In the next paragraphs we describe
  this perturbation along each orbit.

First, by Lemma~\ref{l.real}, there are a point $p\in\Si$ 
and a perturbation $\tilde A$ of $A$ along the orbit of
$p$ such that the corresponding matrix $M_{\tilde A}(p)$ 
is diagonalizable and has only positive real eigenvalues with multiplicity $1$. 
Denote by $\lambda_1<\cdots<\lambda_N$ such eigenvalues and 
by $E_j(p)$ the corresponding eigenspaces. Moreover, by  Remark~\ref{r.real}, if the linear system $A$ satisfies the hypothesis of Remark~\ref{r.transposition}, i.e.\ existence of some point with jacobian greater than one, one can choose the point $p$ and the perturbation $\tilde A$ such that $J(M_{\tilde A}(p))=\prod_{i=1}^N\lambda_i $ is strictly bigger than one.
 
\pagegoal=50pc

By hypotheses, there is a point 
$p_i\in\Si$ (whose orbit is disjoint from the one of $p$) and an
$\varepsilon_0$-perturbation $\tilde A$  (along the orbit of $p_i$) such that the
matrix $M_{\tilde A}(p_i)$  has a complex eigenvalue of rank
$(i,i+1)$. 
Now  we fix  $\varepsilon_1/10$-transitions $[t_{i,0}]$ and $[t_{0,i}]$
from  $p$ to $p_i$ and from $p_i$ to $p$.

To simplify the notation, we write $[M]=[M]_A(p)$, 
$[M_i]=[M]_A(p_i)$, and $M$ and $M_i$ for the corresponding
matrices.
We also write $[\tilde M]$, $[\tilde M_i]$, $\tilde M$,
and $\tilde M_i$ for the corresponding
perturbations of the  words
and matrices.

Observe, once more by Remark~\ref{r.semigroup}, that
 for every $i$  and 
positive $n_1$, $n_2$, and  $n_3$ the word
$$
[t_i(n_1,n_2,n_3)] = ([M])^{n_3}[t_{0,i}]([M_i])^{n_2 }[t_{i,0}][M]^{n_1}
$$
is also an $\varepsilon_1/10$-transition from $p$ to itself
of the system $A$.

The proof of Lemma~\ref{l.transposition} (and of 
Remark~\ref{r.transposition}) now
follows immediately from the next result:
\pagegoal=48pc
\proclaim{Lemma} There are $n_1\geq 0${\rm ,} $n_2\geq 0${\rm ,} and $n_3\geq 0$
such that the word
 $[t_i(n_1,n_2,n_3)]$ defined above admits an
 $(\varepsilon_0+\varepsilon_1)$\/{\rm -}\/perturbation
 such that the corresponding matrix $\tilde T_i(n_1,n_2,n_3)$
left invariant
every $E_j$, $j\notin\{i,i+1\}${\rm ,}
and interchanges $E_i$ and $E_{i+1}$. 
\endproclaim

\demo{Proof} For a fixed $i\in\{1,\dots,N\}$
consider  the 
$\varepsilon_0$-perturbation $\tilde A$ of $A$
 along the orbit of $p_i$ constructed above and
denote by $E(p_i)\oplus F(p_i)\oplus G(p_i)$
 its invariant splitting (over the orbit of $p_i$) 
where $F(p_i)$ is the $2$-dimensional eigenspace corresponding to the complex (conjugate) eigenvalues of $\tilde M_i$ and
$E(p_i)\prec F(p_i)\prec G(p_i)$.
Since $M_{\tilde A}(p)$ 
is diagonalizable we have the splitting
$$
E(p)=\bigoplus_1^{i-1} E_j(p),
\quad
F(p)= E_i(p)\oplus E_{i+1}(p),
\quad
\mbox{and}
\quad
G(p) =\bigoplus_{i+2}^N E_j(p),
$$ 
where $E_j(p)$ is the eigenspace associated to $\lambda_j$.

As in Lemma~\ref{l.goodtransition}, replacing, if necessary,
 the transitions $[t_{0,i}]$ and
$[t_{i,0}]$ by  words
of the form $[M]^n[t_{0,i}][M_i]^n$ and  $[M_i]^n [t_{i,0}][M]^n$ for some big  
$n$, we can assume that 
$[t_{i,0}]$ admits an
 $(\varepsilon_0+\varepsilon_1/10)$-perturbation 
$[\tilde t_{i,0}]$  
such that the corresponding matrix
 $\tilde T_{i,0}$
maps the splitting $E(p)\prec F(p)\prec G(p)$ into 
$E(p_i)\prec F(p_i)\prec G(p_i)$.
Conversely,
we can also suppose that the
matrix 
$\tilde T_{0,i}$ of the
$(\varepsilon_0+\varepsilon_1/10)$-perturbation 
$[\tilde t_{0,i}]$ of   
$[t_{i,0}]$  
maps $E(p_i)\prec F(p_i)\prec G(p_i)$ into
$E(p)\prec F(p)\prec G(p)$.

Our next objective is to get two 
different one-dimensional subspaces of 
$F(p_i)$
and a perturbation of $\tilde M_i$ interchanging such
subspaces. For that
write 
$$
E_i(p_i)= \tilde T_{0,i}^{-1}(E_i(p))\quad\mbox{and}\quad E_{i+1}(p_i) = \tilde T_{i,0} (E_{i+1}(p)).
$$ 
As $E_i(p)$ is a subspace of $F(p)$, 
 $E_i(p_i)$ is a (noninvariant!) subspace of $F(p_i)$. The same argument shows that $E_{i+1}(p_i)$  is also a noninvariant subspace of $F(p_i)$. 

Recall that $\tilde M_i$ has a pair of complex
(nonreal) eigenvalues whose eigenspace is $F(p_i)$. 
So it is an exercise to get $m>0$ 
and an $\varepsilon_1/10$-perturbation $\hat M_i$ of $\tilde M_i$, preserving the splitting 
 $E(p_i)\oplus F(p_i)\oplus G(p_i)$,  such that
$\hat M_i^m(E_{i+1}(p_i))= E_i(p_i)$. 

Consider now the linear map
 $$
B_0= \tilde T_{0,i}\circ \hat M_i^m\circ\tilde T_{i,0}
\colon\cE_p\to\cE_p,
$$ 
and recall that $\cE_p$ denotes the fiber of $p$.
By construction 

\begin{itemize}
\item[i)] $B_0$ preserves the splitting 
$E(p)\oplus F(p)\oplus G(p)$,
\item[ii)]
 $B_0(E_{i+1}(p))=E_i(p)$,
\item[iii)]
 $B_0(E_i(p))$ is a straight line (in the plane $F(p)$) different from $E_i(p)$ and so transverse to $E_i(p)$.
\end{itemize}

Observe that these
three properties also hold for 
$\tilde M^k\circ B_0$ for every $k>0$.

Recalling that 
$E_i(p)\oplus E_{i+1}(p)$ is a dominated splitting over the orbit of~$p$ (for the perturbed system whose matrix is $\tilde M$), using
the transversality in item (iii),  we have that if  $k>0$ is large enough
then  $\tilde M^k(E_i(p))$ is arbitrarily close to $E_{i+1}(p)$.  So we
can choose $k>0$ and an $\varepsilon_1/10$-perturbation $\hat M$
of $\tilde M$ such that 
$B_1 = \hat M\circ \tilde M^k\circ B_0$ satisfies the following two properties:

\begin{itemize}
\item
$B_1$ preserves the splitting $E(p)\oplus F(p)\oplus G(p)$,
\item $B_1(E_i(p))=E_{i+1}(p)$ and $B_1(E_{i+1}(p))=B_1(E_i(p)$. 
\end{itemize}

Note  that these properties of
 $B_1$ are also verified by every map of the form
$\tilde M^{k_1}\circ B_1\circ \tilde M^{k_2}$ ($k_1$ and $k_2>0$). 
Applying the arguments of the proof of Lemma~\ref{l.goodtransition}
to the restrictions of 
$B_1$ to $E(p)$ and $G(p)$ we get $k_1$ and $k_2>0$, and
$\varepsilon_1/10$-perturbations $[N_1]$ of the word $[\tilde M]^{k_1}$, and $[N_2]$ of $[\tilde M]^{k_2}$, 
 coinciding with $[\tilde M]^{k_1}$ and $[\tilde M]^{k_2}$ on $F(p)$, such that 
$$
N_1\circ B_1\circ N_2(E_j(p))=E_j(p)
\quad
\mbox{for every $j\notin\{i,i+1\}$}.
$$

To finish the proof of the lemma it suffices to observe that,
by construction,
the matrix $N_1\circ B_1\circ N_2$ corresponds to a word which is an
$(\varepsilon_0 +3\, \varepsilon_1/10)$-perturbation 
($\varepsilon_0+3\varepsilon_1/10<\varepsilon_0 +\varepsilon_1$) of 
$$
[\tilde M]^{k_1}\,[\tilde M]\cdots[\tilde M]\,[t_{0,i}]\,
[\tilde M_i]^{k}\,[t_{i,0}]\,[\tilde M]\cdots [\tilde M]\,
[\tilde M]^{k_2}.
$$
Thus this word is an $\varepsilon_1$-perturbation of $[t_i(n_1,n_2,n_3)]$ for some $n_1$, $n_2$, and $n_3$. Now the proof of the lemma is complete.
\enddemo

\section{Finest dominated splitting and control of the jacobian\\ in the extremal bundles: Proof of Theorem~4}
\label{s.volume}

6.1. {\it Control of the jacobian over periodic points}.
Let $(\Si, f, \cE, A)$ be a periodic linear system with transitions.
Suppose that $F_1\oplus F_2 \oplus \cdots \oplus F_{k-1}
\oplus F_k$, 
$F_1\prec F_2 \prec \cdots \prec F_{k-1}
\prec F_k$, is the finest dominated splitting of this system.
We call \pagebreak  $F_1$ and $F_k$ 
{\em extremal bundles of the dominated splitting.\/}
Denote by $A_i$ the restriction of $A$ to the subbundle $F_i$.
The goal of this section is to prove some estimates
on the determinant of $A_i$.
Let us begin with the following result:

\proclaim{Lemma} \label{l.vol}
For any $K>0, N\in\NN, L\in \NN${\rm ,} and $\varepsilon>0$ there is $l>0$ with the following property\/{\rm :}

Consider a periodic linear  system $(\Si,f,\cE,A)${\rm ,} of dimension $N$ and
bounded by $K${\rm ,} having  an $L$\/{\rm -}\/dominated splitting $E\oplus F$ such that 
\begin{itemize}
\item
the subbundle
$E$ does not admit any nontrivial $l$\/{\rm -}\/dominated splitting{\rm ,} and
\item 
there is a point
$p\in \Si$ such that
$\mbox{{\rm det}}(M_{A_E}(p))>1.$
\end{itemize}
Then there are
an $\varepsilon$\/{\rm -}\/perturbation  $\tilde A$ of $A$ and $x\in \Sigma$
such that all the eigenvalues of $M_{\tilde A}(x)$ have modulus
strictly bigger than $1$.
\endproclaim

Clearly, there is a version of this lemma where the bundle $F$
has no\break $l$-dominated splitting
and $\mbox{det}(M_{A_F}(p))<1$; then
the  moduli of the eigenvalues of the perturbation $\tilde A$ are 
strictly smaller than $1$.

 We will apply  this lemma twice to the finest dominated splitting of a system,
first taking $E=F_1$ and $F=F_2 \oplus \cdots \oplus F_k$, and
second taking
$E=F_1\oplus \cdots \oplus F_{k-1}$ and $F=F_k$. 
Let us now prove the lemma.

\demo{Proof} Note first that 
there is $\delta_L>0$ such that, for every system $C$\break $\delta_L$-close to
$A$,  we can define the
continuation $E_C\oplus F_C$ of the splitting $E\oplus F$, which
is also dominated (recall the comments after
Definition \ref{d.dominated}).

Let $\tau=\inf\{\varepsilon,\delta_L\}$. 
Consider the $l_0>0$ given by Proposition~\ref{p.rank}
associated to $\tau/2$, and fix $l=2\,l_0$.

Now take a system $(\Si,f,\cE,A)$ satisfying the hypotheses of the lemma.
From Corollary~\ref{c.restr}, the  system $A_E$ induced by $A$ on the
subbundle $E$ admits transitions and, by hypothesis, does not admit any
$l$-dominated splitting.
In particular, by Proposition \ref{p.rank}, for every 
$0 \le i< \dim (E)$, there is a\break $\tau/2$-perturbation
$B_i$ of $A_E$ 
having a complex eigenvalue of rank $(i,i+1)$, i.e.\ there is  $x_i\in \Si$ such that $M_{B_i}(x_i)$ has a complex eigenvalue of rank $(i,i+1)$. 

Moreover, by hypothesis, there is a point $p\in\Si$ for which the modulus of the jacobian $J(M_{A_E})(p)$ is strictly bigger than one.

The previous comments mean that
 taking $\varepsilon_0=\tau/2$ , the system\break $(\Si,f,E,A_E)$ satisfies all the hypotheses of Proposition~\ref{p.homothety} and
Remark~\ref{r.rema}. So there is a $\tau$-perturbation $\tilde A_E$ of $A_E$ and a point $x\in\Si$ such that $M_{\tilde A_E}(x)$ is a homothety of ratio bigger than one. 

Finally, by Lemma~\ref{l.perturbations}, $\tilde A_E$ is the restriction
to $E$ of some $\tau$-perturbation $\tilde A$ of $A$ which coincides with $A$ over $F$. Since the splitting $E\oplus F$  is  dominated for $\tilde A$,  we get that all the eigenvalues of  $M_{\tilde A}(x)$ associated to $F$ 
 necessarily have modulus bigger than one.
This  ends the proof of the lemma.  
\enddemo

Theorem 4 follows from the next proposition 
that 
is a consequence of the ergodic closing lemma in 
\cite{Ma3}
and
whose proof
(very similar to Ma\~n\'e's argument in \cite{Ma3} to get hyperbolicity)
 we postpone
until the next subsection:

\proclaim{Proposition} \label{p.vol} Let $f$ be a diffeomorphism and $\La_f(U)$ an $f$\/{\rm -}\/invariant compact set which is 
maximally invariant 
 in some neighbourhood $U$ of it.
Suppose that  $E\oplus F${\rm ,} $E\prec F${\rm ,} is  a dominated splitting  of $T_{\La_f(U)} M$ for $f_*$.

By shrinking $U${\rm ,} if necessary{\rm ,} there is a $C^1$-neighbourhood $\cU$ of $f$ such that 
for every $g\in \cU$
the maximal invariant set $\La_g(U)$ of $g\in\cU$ has a dominated splitting $E_g\oplus F_g$ which is the 
continuation of $E\oplus F$. Then one has that
\begin{itemize}
\item either there are an arbitrarily small $C^1$\/{\rm -}\/perturbation $g$ of $f$ and a 
hyperbolic periodic point $p\in\La_g(U) $ such that $g_*$ expands the  volume on $E_g(p)${\rm ,} 
\item or $f_*|_E$  contracts  the volume uniformly.
\end{itemize}

\endproclaim

Let us now end the proof of Theorem~4.

\demo{Proof of Theorem~{\rm 4}}
 Let $\La_f(U)$ be a (nontrivial) robustly transitive set. Denote by $N$ the dimension of the ambient manifold and let $K$ be a strict upper bound of the norms of $f_*$ and $f_*^{-1}$. 
 Take $\varepsilon>0$ such that for every $g$ 
 $\varepsilon-C^1$-close to $f$ the set $\La_g(U)$ is transitive  and the norms of $g_*$ and $g_*^{-1}$ are bounded by $K$. 

The proof is by contradiction. Let $F_1\oplus F_2\oplus \cdots \oplus F_k$
be the finest dominated splitting of $f_*$. Write
$E=F_1$ and $F=F_2\oplus \cdots \oplus F_k$, and fix $L$ such that the
splitting $E\oplus F$ is $2L$-dominated. Now let  $l>0$
 be the constant associated to $K$, $N$, $2L$, and $\varepsilon/2$ in Lemma~\ref{l.vol}.   

If $f_*|E$ does not uniformly contract  the volume in $E$ then, by
Proposition~\ref{p.vol}, there are $g$ arbitrarily close to $f$ and a hyperbolic periodic point
$p\in \La_g(U)$ such that $g_*$ expands the volume in $E_g(p)$. 
Since, by hypothesis, $\La_g(U)$ is  
 $C^1$-robustly transitive, 
using Lemma~\ref{l.ro} we can assume (after a perturbation) that the
relative homoclinic class  $H(p,g,U)$ of $p$ is the whole $\La_g(U)$. 

Consider the dense subset 
 $\Si\subset\La_g(U)$ consisting  of all
the hyperbolic periodic points 
of $\La_g(U)$ homoclinically related to $p$. Then $g$ induces the periodic linear system $(\Sigma, g, T_{\Si}M, g_*)$, of dimension $N$ and  bounded by $K$,
with transitions, see Lemma~\ref{l.Df}.

 Moreover, if $g$ is close enough to $f$,
$E(g)\oplus F(g)$ is  an $L$-dominated
splitting of this system and the bundle $E(g)$ does not admit any $l$-dominated splitting (this last assertion follows from Lemma~\ref{l.dense}). So applying Lemma~\ref{l.vol},
we get an $\varepsilon/2$-perturbation $B$ of  $(\Sigma, g, T_{\Si}M,
g_*)$ and a periodic point $q\in \Si$ such that all the eigenvalues of
$M_B(q)$ are bigger than one (in modulus). Using Franks'
lemma  we get that the (nontrivial) maximal invariant set
in $U$ of some $\varepsilon$-perturbation $g$ of $f$ contains a repeller (precisely the point $q$).  This contradicts the choice of $\varepsilon$ (robust transitivity of $\La_g(U)$).

To end the proof of the theorem it remains to
get the uniform expansion of the volume in $F_k$, this follows as above by replacing $f$ by $f^{-1}$. 
\enddemo

6.2. {\it Ma{\rm \~{\it n}}{\rm \'{\it e}}\/{\rm '}\/s ergodic closing lemma\/{\rm :} Proof of 
Proposition}~\ref{p.vol}.
 
\numbereddemo{Definition} Let $f$ be a diffeomorphism defined on  a compact
 manifold $M$ endowed with a Riemannian metric $d$. 
A point $x$ is {\em well closable\/} (for $f$) if for 
every $\varepsilon >0$ there are $g$ $\varepsilon-C^1$-close to $f$ 
and a periodic point $y$ of $g$ such that $d(f^i(x),g^i(y))<\varepsilon$ for every
$0\leq i<k(y)$, where $k(y)$ is the period of $y$.
We denote by $\cW(f)$ the set of well closable points of $f$.  
\enddemo

We have the following result (see \cite{Ma3}),

\nonumproclaim{Theorem {\rm (ergodic closing lemma)}} Let $f$ be a diffeomorphism and $\mu$
an $f$\/{\rm -}\/invariant probability.
Then $\mu$\/{\rm -}\/almost every point is well closable{\rm ,} i.e.\ $\mu(\cW(f))=1. $
\endproclaim

Suppose  now  that $\La_f(U)$ is 
a locally maximal set in a neighborhood $U$ of it
 and that $E\oplus F$,
$E\prec F$, is a dominated splitting of $f_*$ over $T_{\La_f(U)} M$. Recall that 
the bundle $E$ is continuous (see Lemma~\ref{l.dense}); thus
 $(\La_f(U), f, E, f_*|_E)$ is a continuous linear system. 

Since $E$ is endowed with a
continuous Euclidian metric,  we can define
the modulus of the determinant of $f_*|_E$, the {\em  jacobian
of $f$ on $E$,\/} denoted by $|J(f,E)|\colon \La_f(U)\to \RR$, which is 
a continuous  (and so integrable) positive function.
Thus $\log(|J(f,E)|)\colon\La_f(U)\to \RR$ is well defined and continuous. 
Moreover, by shrinking $U$, if necessary,
we have that for every $g$ close enough to $f$ there is defined a 
dominated continuation $E_g\oplus F_g$ of the splitting $E\oplus F$.
Thus we can define the function 
$\log(|J(g,E_g)|)$ depending  continuously on $g$
(observe that the subbundles $E_g$ and $F_g$ 
depend  continuously on $g$). 
 
The first step to prove Proposition~\ref{p.vol} is 
the following lemma.

\proclaim{Lemma} 
\label{l.measure}
Assume there is an $f$\/{\rm -}\/invariant probability measure $\mu$ supported on $\La_f(U)$ such that 
$$
\int \log(|J(f,E)|)\, d\mu \geq 0.
$$
Then there are $g$ arbitrarily $C^1$\/{\rm -}\/close to $f$ and
 a periodic orbit $y\in\La_g(U)$ of $g$ where $g_*$ expands the volume of $E_g${\rm ,}
i.e.\ $$
|J(g^k,E_g)(y)| >1,
\quad
\mbox{where $k$ is the period of $y$}.
$$ 
\endproclaim

\demo{Proof}
Observe first that we can assume that $\mu$ is ergodic: Otherwise
it is enough to consider the  decomposition of $\mu$ into ergodic invariant measures; then  for at least one of them the integral of
$\log(|J(f,E)|)$  is also positive. 

By the ergodic closing lemma, there is a $\mu$-generic point $x$ which is well closable. If $x$ is periodic we have nothing to do. 
In the other case, there are sequences of diffeomorphisms
$g_n$ converging to $f$ in the $C^1$-topology,  
of periodic points
$y_n$
(of period $k_n$) of $g_n$, and of numbers $\varepsilon_n\to 0$, such that 
$$
d(f^i(x),g_n^i(y_n))<\varepsilon_n,
\quad 
\mbox{for every $0\leq i<k_n-1$}.
$$
In particular, if $\varepsilon_n$ is small enough, this ensures that the point $y_n$ belongs to $\La_{g_n}(U)$. Moreover, 
$$
\frac 1{k_n}\sum_0^{k_n-1} \log (|J(g_n,E_{g_n}|)(g_n^i(y_n))
\to \int \log(|J(f,E)|)\, d\mu \geq 0.
$$
 So given any $\delta>1$, taking $n$ large enough,
we have that  the linear system defined over the orbit of $y_n$ obtained by multiplying $(g_n)_*|E_{g_n}$ by the scalar $\delta$ expands the volume on $E_{g_n}(y_n)$.  The conclusion in the lemma now
 follows immediately  by
Franks' lemma (see the very beginning of Section~\ref{s.ls}). 
\enddemo

Using Lemma~\ref{l.measure} we have that Proposition~\ref{p.vol} (thus  Theorem~4) is  a direct consequence of:

\proclaim{Lemma} \label{l.buildmeasure}
Let $E\oplus F${\rm ,} $E\prec F${\rm ,}
be a dominated splitting of $f$ over $\La_f(U)$. 
Assume that for any $N\in \NN$ the jacobian $|J(f^N,E)|$ is not
bounded uniformly from above by one. 
Then there is an $f$\/{\rm -}\/invariant  measure $\mu$ such that 
$$
\int \log(|J(f,E)|)\, d\mu \geq 0.
$$
\endproclaim

{\it Proof}.
By hypothesis, given any $N>0$ there exists some point $x_N\in\La_f(U)$
such that $|J(f^N,E)(x)|\geq 1$. 
Write
$$
\mu_N = \frac1N\sum_0^{N-1}\delta(f^i(x_N)),
$$
where $\delta(z)$ is the Dirac measure at the point $z$.
As the space of probabilities is compact for the weak topology, 
there is a  subsequence $N_i$ such that $\mu_{N_i}$ converges weakly to some probability measure $\mu$. 

A classical elementary argument proves that $\mu$ is $f$-invariant:
$f_*(\mu)-\mu$ is the weak limit of 
$\frac1{N_i}\, (\delta(f^{N_i}(x_{N_i}))-\delta(x_{N_i}))$,
which converges to $0$. \pagebreak
Finally, observing that
\begin{eqnarray*}
\int  \log(|J(f,E)|)\, d\mu_N& = &
\frac1N \, \sum_0^{N-1} \log (|J(f,E)(f^i(x_N))|)\\
&= &
\frac{1}N \, \log (|J(f^N,E)(x_N)|) \geq 0,
\end{eqnarray*}
one deduces immediately  that $\int  \log(|J(f,E)|\, d\mu\geq 0$.    
\hfill\qed

\section{The conservative case}
\label{s.con}

\vglue-4pt

In this section we translate some of our constructions
into the conservative context (volume-preserving
$C^1$-diffeomorphisms). 

In what follows the compact manifold $M$ is endowed with a smooth volume
form $\omega$, and we denote by Diff$^1_\omega(M)$ the set of
$C^1$-diffeomorphisms preserving this volume form $\omega$ endowed with
the usual $C^1$-topology. Following the traditional terminology, a {\em conservative diffeomorphism\/} is an element of Diff$^1_\omega(M)$. Here we only consider manifolds of dimension $N$ strictly bigger than one.

Observe that, given a linear system $(\Si,f,\cE,A)$, using the Euclidian metric on the bundle $\cE$, we can define the modulus of the determinant of the linear maps $A(x)$, denoted by $J_A(x)=|\det A(x)|$.

\vglue5pt {\it Definition} 7.1.
A periodic linear system $(\Sigma, f, \cE, A)$ 
is {\em conservative\/} if
$J_A(x)=1$ for all $x\in \Si$.
\advance\theoremcount by 2

\vglue5pt {\it {R}emark} 	7.2.
Given a linear system $(\Si,f,\cE,A)$ of dimension $N$ we 
define its {\em conservative part\/}  $A^c$  by 
$$
A^c(x)= J_A(x)^{\frac{-1}N}\cdot A(x).
$$
Clearly, if $A$ is conservative then $A=A^c$.
The map $A\mapsto A^c$ is continuous. Therefore for any $K>0$ and
$\varepsilon>0$ there is $\varepsilon_1>0$ such that if $A$ is a
conservative system bounded by $K$ and $B$ is a (\/{\it a~priori}
nonconservative) $\varepsilon_1$-perturbation of it, then $B^c$ is
$\varepsilon$-close to $A$. 

Finally, if the initial system $A$ is continuous,  periodic, and with transitions, then the same holds for its conservative part $A^c$.
\label{r.con}
\vglue5pt

We begin with a
straightforward corollary of Proposition~\ref{p.main}:

\proclaim{Proposition} \label{p.conlin} 
For any $K>0$, $N>0${\rm ,} and $\varepsilon>0$ there is 
$l>0$ such that for every continuous periodic conservative 
 linear system $(\Si,f,\cE,A)${\rm ,} of dimension $N$ and bounded by $K${\rm ,} with transitions
one has that{\rm ,}
\begin{itemize}
\item either $A$ admits an $l$\/{\rm -}\/dominated splitting{\rm ,}
\item or there is a conservative  $\varepsilon$-perturbation $\tilde A$ of $A$ such that
$M_{\tilde A}(x)$ is the identity for some $x\in\Si$.
\end{itemize}

\endproclaim

{\it Proof}.
According to Proposition~\ref{p.main},
there is $l$ with the following property:

Let $A$ be a  continuous periodic conservative linear system, of
dimension~$N$ and bounded by $K$, and admitting transitions. If $A$ does
not admit any\break $l$-dominated splitting then
 there is an $\varepsilon_1$-perturbation $B$ of $A$ 
such that the matrix $M_B(x)$ is a homothety for some $x\in \Si$. 
Then, by Remark~\ref{r.con}, $B^c$ is a conservative $\varepsilon$-perturbation of $A$ and $M_{B^c}(x)$ is 
a homothety with determinant~$1$ (because the system is conservative); thus
it is the identity. This completes the proof of the 
proposition.
\hfill\qed\vglue5pt

To proceed with our proofs in the
conservative setting we need suitable versions
of our perturbation  lemmas (Franks' lemma and Hayashi's connecting lemma)
 for volume-preserving diffeomorphisms.

\vglue-16pt
\phantom{odd}
\proclaimtitle{conservative version of  Franks'   lemma}
\proclaim{Proposition}
\label{p.fc}
Consider a conservative diffeomorphism
$f$ and a finite $f$\/{\rm -}\/invariant set $E$. 
Assume that $B$ is a conservative $\varepsilon$\/{\rm -}\/perturbation of $f_*$ along $E$. 
Then for every neighbourhood $V$ of $E$ 
there is a conservative diffeomorphism $g$ arbitrarily $C^1$\/{\rm -}\/close to $f$  coinciding with $f$
on $E$ and out of $V${\rm ,}
and  such that $g_*$ is equal to $B$ on $E$.
\endproclaim

 \vglue-4pt

As we did not find any precise reference for this probably well-known result,   let us give here the sketch of its proof:

\vglue5pt {\it Proof}.
The proof is based on the following  elementary fact of linear\break algebra:

\vglue-16pt
\phantom{odd}

\proclaim{Lemma}
\label{l.rots}
For any $N>1$ and $\varepsilon>0$ there is a neighbourhood $\cG$ of the
identity in ${\rm SL}(N,\RR)$ such that any matrix $A\in\cG$ can be written as
a product $B_1\circ B_2 \circ\cdots\circ B_{4N-4}${\rm ,} where the $B_i=
P_i\circ R_i\circ P_i^{-1}${\rm ,}  $P_i$ and $R_i$ are $\varepsilon$\/{\rm -}\/close to the identity in ${\rm SL}(N,\RR)${\rm ,}
 and $R_i$ is a rotation.
\endproclaim
\vglue-3pt

{\it Proof}. 
We just give the main steps of the proof.
All the matrices we consider will be in ${\rm SL}(N,\RR)$.
We have the following properties:

\vglue2pt\noindent \hskip1em \hangindent =20pt \hangafter=1 $\bullet$ Every matrix $A$ close to identity can be written
in the form $L\circ L^{-1}\circ A$, where $L$ is diagonal  with real eigenvalues of multiplicity one and close to~$1$,  $L^{-1}\circ
A$ is diagonalizable and  has eigenvalues close to $1$, and  the matrix of change of coordinates is close to identity.
\vglue2pt\noindent \hskip1em \hangindent =20pt \hangafter=1  $\bullet$ Any diagonal matrix in ${\rm SL}(N,\RR)$ can be written as
the product of $N-1$ diagonal matrices whose eigenvalues are  equal to one, except for two of them (inverse one to the other).
\vglue2pt\noindent \hskip1em \hangindent =20pt \hangafter=1 $\bullet$  Any diagonal matrix $D$
close to the identity, $D\in {\rm SL}(2,\RR)$, can be written as $R\circ R^{-1}\circ   D$, where $R$ is a rotation, $R^{-1}\circ D$ is conjugate to a rotation by some matrix $P$, and 
the matrices $R$, $R^{-1}\circ D$, and $P$ are close to the  identity. 
\vglue2pt
\noindent The lemma follows immediately from these three properties.
\hfill\qed\pagebreak

Observe that Franks' lemma consists of local perturbations around finitely many points. Taking local charts at these points we can
consider that the volume is the Lebesgue measure (Leb), see for instance \cite{Mo}. So    Franks' conservative lemma  follows
from the next lemma:

\proclaim{Lemma}
\label{l.fc} For every $N\in \NN$ and $\varepsilon>0$ there is a neighbourhood $\cG$ of the identity in ${\rm SL}(N,\RR)$ 
such that for every $A\in \cG$ there is  $h\in$ 
{\rm Diff}$^1_{\rm Leb}(\RR^N)$ satisfying the following properties\/{\rm :}\/

\vglue4pt\noindent \hskip18pt \hangindent =20pt \hangafter=1  $\bullet$ $h$ coincides with the identity outside the unit ball at
the origin{\rm ,}

\vglue4pt\noindent \hskip18pt \hangindent =20pt \hangafter=1  $\bullet$
 $h(0)=0$ and $h_*(0)=A${\rm ,}

\vglue4pt\noindent \hskip18pt \hangindent =20pt \hangafter=1  $\bullet$
$\|h_*-{\rm Id}\|<\varepsilon$.
\endproclaim

To prove this lemma it is enough to see that its proof is very easy if $A$ is a rotation or conjugate to a rotation.
Since  Lemma~\ref{l.rots} allows us to write $A$ as the product of $4N-4$ of such maps, all them close to the identity,
the general case follows from the simple first case. \hfill\qed
\enddemo

Exactly as Proposition~\ref{p.quantitative} follows from Franks' lemma and Proposition~\ref{p.main}, we
 deduce from Proposition \ref{p.fc} (conservative version of Franks' lemma) and Propositions~\ref{p.conlin} 
the following conservative version of Proposition~\ref{p.quantitative}:

\proclaim{Proposition} 
\label{p.conservative} Given any $K>0$, $N>0${\rm ,} and $\varepsilon>0$ there is 
$l(\varepsilon, K)\break \in \NN$ such that for every
conservative diffeomorphism $f$ 
on a Riemannian\break $N$\/{\rm -}\/dimensional manifold $M${\rm ,} 
with derivatives $f_*$ and $f_*^{-1}$ bounded by $K${\rm ,} 
and any saddle $p$ of $f$ 
having a nontrivial homoclinic class $H(p,f)${\rm ,} 
one has that\/{\rm :}  
\begin{itemize}
\item 
Either the homoclinic class $H(p,f)$ admits
 an $l(\varepsilon, K)$\/{\rm -}\/dominated splitting{\rm ,} 
\item 
or for every  neighbourhood $U$ of $H(p,f)$  and 
 $k\in \NN$ there are a conservative diffeomorphism $g$
 $\varepsilon$-$C^1$\/{\rm -}\/close to $f$ and
$k$ periodic points $x_i$ of $g$ arbitrarily
close to $p${\rm ,} whose orbits are contained in $U${\rm ,}
 such that
the derivatives $g_*^{n_i}(x_i)$ are  the identity {\rm (}$n_i$
is the period of $x${\rm ).} 
\end{itemize}

\endproclaim

Observe that this proposition implies Theorem~\ref{t.herman} (in fact, it is a quantitative version of Theorem 5).

\vglue6pt 7.1. {\it Proof of Theorem}~6.
The first step  is the following lemma:

\proclaim{Lemma} 
\label{l.cres}
There is a residual subset $\cR\subset$ {\rm Diff}$^1_\omega(M)$ of diffeomorphisms $f$ such that the nontrivial homoclinic
classes of hyperbolic periodic points of $f$ are dense in $M$.
\endproclaim

\demo{Proof}
Let us begin by recalling that, for conservative diffeomorphisms, the recurrent points are dense in $M$.
Moreover, using the
$C^1$-closing lemma in the conservative case (see
\cite{Rb}), one has that the conservative diffeomorphisms  whose periodic orbits are hyperbolic
and dense in the ambient manifold form a dense subset of Diff$^1_\omega(M)$. Moreover, by the continuity of the hyperbolic periodic orbits, this dense subset is in fact residual.

Suppose that $p$ is a hyperbolic periodic point 
of  a diffeomorphism $g$ such that its periodic points are hyperbolic and dense in the manifold.
Then, given fundamental
domains $D^s$ and $D^u$ of $W^s(p,g)$ and $W^u(p,g)$,
there is a sequence of periodic points $q_i$
converging to some $z\in D^u$ such that
$g^{k_i}(q_i)\to y\in D^s$ for some sequence $k_i\to +\infty$.
This implies that the invariant manifolds of any periodic point $p$ of $g$
satisfy the hypotheses of the conservative connecting
lemma of Xia  \cite{X}:

\nonumproclaim{Theorem  {\rm (conservative connecting lemma)}}
Let $M$ be a compact manifold{\rm ,} $g$ a conservative
diffeomorphism,{\rm } and $p$ a hyperbolic periodic point
of $g\in$ {\rm Diff}$^1_{\omega}(M)$.

Suppose that there are sequences of points $(x_i)${\rm ,}
$x_i\to z\in W^u(p)${\rm ,} and integers $(n_i)$
such that $n_i\to \infty$ and $g^{n_i}(x_i)\to y\in W^s(p)$. 

Then there is $h\in$ {\rm Diff}$_\omega^1(M)$ arbitrarily $C^1$\/{\rm -}\/close
to $g$ such that $W^u(p_h,h)$ intersects transversely
$W^s(p_h,h)$ at $z$ and $h^k(z)=y$ for some $k>0$.
\endproclaim

Now, a classical argument (see, for instance, the proof of  \cite[Prop.~1.1]{BD2}) gives that 
there is a dense open subset $\cR_n$ of Diff$^1_\omega(M)$, a
diffeomorphism whose nontrivial homoclinic classes are $\frac1n$-dense. To end the proof of the lemma it is enough to take $\cR=\bigcap_{n>0}\cR_n$.
\enddemo

Now to end the proof of Theorem~6  consider a diffeomorphism
$f\in$ Diff$^1_{\omega}(M)$ and $\varepsilon>0$ such that there is no
$\varepsilon$-perturbation of $f$ having periodic points whose
derivative is the identity. Thus, by Proposition \ref{p.conservative},
there is $l$ such that every nontrivial homoclinic class of any
conservative $\varepsilon/2$-perturbation $g$ of $f$ admits an
$l$-dominated splitting. For such a $g$ we define $\La_i(g)$,
$i=1,\dots,N-1$, as the closure of the union of the nontrivial homoclinic
classes with an $l$-dominated splitting $E\prec F$
of dimension $i$ (i.e.\ $\dim(E)=i$). 
By Lemma~\ref{l.dense}, this invariant compact set $\La_i(g)$ has an
$l$-dominated splitting. 

Using Lemma~\ref{l.cres}, we can take a sequence of diffeomorphisms $g_n\in\cR$ as above converging to $f$. Then, for each $n$, the union of the $\La_i(g_n)$ is the whole manifold. We define $K_i(f)$ as the topological upper limit set of the $\La_i(g_n)$, i.e.\ $$
K_i(f)=\limsup_{n\to\infty}\La_i(g_n) = \bigcap_k\mbox{closure}(\bigcup_{n\geq k} \La_i(g_n)).
$$ 
Again by Lemma~\ref{l.dense}, the set 
$K_i(f)$ admits an $l$-dominated splitting. Finally, by construction, the manifold $M$ is the union of the $K_i(f)$.

The argument  above shows that, if $f$ does not admit any
$\varepsilon$-pertubation with a periodic orbit whose derivative is the
identity, then $M$ is the union of finitely many invariant compact sets
with $l$-dominated splittings. Otherwise, there are a conservative perturbation $g$ of $f$ and 
a homoclinic class of a periodic point of $g$ 
whose induced  periodic linear system can be perturbed to get one point such that its linear map is the identity. Using the transitions we can get an arbitrarily large number of such points. 
Now the result follows from the conservative version of 
 Franks'  lemma (Proposition~\ref{p.fc}). 
\hfill \qed

\vglue12pt 7.2. {\it Volume properties of dominated splittings of conservative systems}.
We end this paper by giving some volume properties of
dominated splitting of conservative linear
systems.

Proposition~0.5 is a direct consequence of the following lemma:

\proclaim{Lemma}
\label{l.cEF}
Let $(\Si,f, \cE, A)$ be a conservative linear
system with an\break $l$\/{\rm -}\/dominated splitting
 $E\oplus F$, $E\prec_l F$. Then  
$$
|\mbox{{\rm det}}(A^l_E)(x))|\leq \frac1{\sqrt 2}
\quad
\mbox{and}
\quad
|\mbox{{\rm det}}(A_F^l)(x))|\geq \sqrt 2,
$$ 
for every $x\in\Si${\rm ;} recall that $A^l(x)= A(f^{l-1}(x))\circ\cdots\circ A(x)$.
\endproclaim

\demo{Proof}
Since the system is conservative,
$$
|\mbox{det}(A^l_E)(x)
\,
\mbox{det}(A^l_F)(x)|=1.
$$
In particular, the second inequality in the lemma follows from the first one.

We argue by contradiction. If the conclusion in the lemma
does not hold then there is 
$x\in \Si$ such that $|\det((A^l_E)(x))|>\frac1{\sqrt 2}$. 
Thus the modulus of the determinant of the
 matrix $(\sqrt 2)^{\frac1{\dim(E)}}\cdot A_E^l(x)$ is bigger than $1$. So this matrix
expands at least one vector. 
Thus there is some unit vector $u\in E(x)$ such that 
$$
\|A^l_E(x)(u)\|> 2^{-\frac1{2\dim(E)}}.
$$
As the system is $l$-dominated,  given any unit vector $v\in F(x)$ we have 
$$
\|A^l_F(x)(v)\|> 2^{1-\frac1{2\dim(E)}}\geq\sqrt 2.
$$
 In particular, 
$$
\left(|\det A_F^l(x)|> \sqrt 2^{\dim(F)}\geq \sqrt 2\right) \Longrightarrow \left(|\det(A^l_E)(x)
\,
\det(A^l_F)(x)|> 1\right),
$$  
contradicting that the system is conservative.
\enddemo

The next lemma immediately implies  Proposition~0.5:

\proclaim{Lemma}
\label{l.cds}
Let $f$ be a conservative diffeomorphism
with a dominated splitting $E\oplus F${\rm ,} $E\prec F$. 
Then
there is $\ell$ such that
$$
|\mbox{{\rm det}}(f_*^\ell(x)|_E)|<\frac{1}{2}
\quad\mbox{and}\quad
|\mbox{{\rm det}}(f_*^{-\ell}(x))|_F|>\frac{1}{2}
\quad \mbox{for every $x\in M$.}
$$
\endproclaim

\demo{Proof}
Observe first that,
using the conservative 
version of the closing lemma,
after a perturbation, we can assume that for every
$k$ the
 periodic points of $f$ of period bigger
than $k$ are dense in $M$. 

The proof now is by contradiction. 
Suppose that the thesis is false;
then for every $t>0$
there is $g_t\in$ Diff$^1_\omega
(M)$ close to $f$  and  a periodic
point $P_t$ such that
$$
g_t^n(P_t)=P_t
\quad 
\mbox{and}
\quad
\mbox{det}((g^n_t)_*(P_t)|_E)>(1-t)^n.
$$
It is not hard to see that we can assume
that the periods of the points $P_t$ go to
infinity as $t\to 0$ (otherwise, 
after perturbation, one gets a linear 
system $A$ and periodic point $x$ such that
$M_A(x)$ has an eigenvalue of modulus bigger
than $1$ in $E$ and an eigenvalue of
modulus less than $1$ in $F$, contradicting the dominance of the splitting).

Taking $t=1/n$,
arguing as in Lemma~\ref{l.cds}, we  perturb each
$(g_{1/n})_*$ along the orbit of $P_{1/n}$ to get a linear system $B$ close to $f_*$ such that
$E\oplus F$ is a dominated splitting of $B$ which does
not satisfies Lemma \ref{l.cEF}. This gives a 
contradiction.
\enddemo

\AuthorRefNames [BeVi2]

\end{document}

\vspace{1cm}



\end{document}